\documentclass[reqno, noamsfonts]{amsart}
\usepackage[utf8]{inputenc}
\usepackage[charter, cal=cmcal, expert]{mathdesign}
\usepackage[scaled=0.96, theoremfont]{XCharter}
\usepackage{mathtools, shuffle, microtype}
\usepackage[margin={3.5cm, 2.5cm}, a4paper]{geometry}
\usepackage[scr=boondoxo, bb=boondox]{mathalpha}
\usepackage{subdepth}
\usepackage{cite, bookmark}
\usepackage[nameinlink]{cleveref}
\usepackage[dvipsnames]{xcolor}
\usepackage[ocgcolorlinks]{ocgx2}
\usepackage[lsep=1.0, ssep=0.9]{ntrees}
\usepackage{tikz-cd}
\usepackage{dirtytalk}
\hypersetup{citecolor=RoyalBlue, linkcolor=BrickRed}
\usepackage{mleftright}
\usepackage{scalerel}
\mleftright

\numberwithin{equation}{section}
\usepackage{thmtools}
\declaretheorem[numberwithin=section]{theorem}
\declaretheorem[sibling=theorem]{proposition}
\declaretheorem[sibling=theorem]{lemma}
\declaretheorem[sibling=theorem]{corollary}

\theoremstyle{definition}

\declaretheorem[sibling=theorem,qed=\(\blacklozenge\)]{example}
\declaretheorem[qed=\(\blacktriangle\), sibling=theorem]{definition}
\theoremstyle{remark}
\declaretheorem[sibling=theorem,qed=\(\blacklozenge\)]{remark}

\title{Branched It\^o Formula and natural It\^o-Stratonovich isomorphism}
\author[C.~Bellingeri]{Carlo Bellingeri}
\address{Université de Haute-Alsace, France}
\curraddr{Université de Haute-Alsace, France}
\email{carlo.bellingeri@uha.fr}
\urladdr{https://sites.google.com/view/carlobellingeri/home}
\author[E.~Ferrucci]{Emilio Ferrucci}
\address{Mathematical Institute, University of Oxford, United Kingdom}
\curraddr{Mathematical Institute, University of Oxford, United Kingdom}
\urladdr{https://people.maths.ox.ac.uk/rossiferrucc/}
\email{emilio.rossiferrucci@maths.ox.ac.uk}
\author[N.~Tapia]{Nikolas Tapia}
\address{Technische Universität Berlin, Germany}
\curraddr{Weierstrass Institute, Berlin, Germany}
\email{tapia@wias-berlin.de}
\urladdr{https://wias-berlin.de/people/tapia}
\hypersetup{pdfauthor={Bellingeri, Carlo and Ferrucci, Emilio R. and Tapia, Nikolas}, pdftitle={Branched It\^o Formula and It\^o-Stratonovich correction}, unicode=true}

\newcommand{\cH}{\mathcal{H}}
\newcommand{\cF}{\mathcal{F}}

\newcommand{\scrF}{\mathscr{F}}
\newcommand{\scrH}{\mathscr{H}}
\newcommand{\scrC}{\mathscr{C}}

\newcommand{\scrT}{\mathscr{T}}
\newcommand{\cT}{\mathcal{T}}
\newcommand{\tee}{\mathbin{\top}}
\newcommand{\bX}{\mathbf{X}}
\newcommand{\bH}{\mathbf{H}}
\newcommand{\bY}{\mathbf{Y}}

\newcommand{\dif}{\mathrm{d}}
\newcommand{\D}{\mathrm{D}}

\newcommand\p{{\lfloor \rho \rfloor}}

\newcommand{\lAngle}{\langle \mkern-4.5mu \langle}
\newcommand{\rAngle}{\rangle \mkern-4.5mu \rangle}

\newcommand{\bZ}{\mathbf{Z}}
\newcommand{\cP}{\mathcal{P}}
\newcommand{\cQ}{\mathcal{Q}}

\newcommand{\hck}{\mathcal H_\mathrm{CK}}
\newcommand{\pck}{\cP}
\newcommand{\hgl}{\mathcal H_\mathrm{GL}}

\newcommand{\scrf}{\mathscr f}
\newcommand{\scrg}{\mathscr g}
\newcommand{\scrh}{\mathscr h}
\newcommand{\scrr}{\mathscr r}
\newcommand{\scrl}{\mathscr l}

\def\Fbar{\overline{F}}

\mathchardef\mhyphen="2D
\newcommand{\cat}[1]{\underline{\mathrm{#1}}}
\newcommand\binf{\mathbf{B}_\infty}

\def\bnohat{\cat{comm\mathbf{B}}_\infty}
\def\argument{\, \cdot \,}
\def\bbR{\mathbb R}
\def\cQ{\mathcal{Q}}
\def\cN{\mathcal{N}}
\def\cE{\mathcal{E}}
\def\cB{\mathcal{B}}
\def\tens{{\textstyle\bigotimes}}

\def\pairing{\langle\argument,\argument\rangle}
\def\cL{\mathcal{L}}

\DeclareFontFamily{U}{bigshuffle}{}
\DeclareFontShape{U}{bigshuffle}{m}{n}{
	<5-8> s*[1.6] shuffle7
	<8->  s*[1.4] shuffle10
}{}
\DeclareSymbolFont{BigShuffle}{U}{bigshuffle}{m}{n}
\DeclareMathSymbol\bigshuffle{\mathop}{BigShuffle}{"001}
\DeclareMathSymbol\bigcshuffle{\mathop}{BigShuffle}{"002}

\def\qshuffle{\mathbin{\widetilde{\shuffle}}}
\def\bigqshuffle{\widetilde{\bigshuffle}}

\def\scrt{\mathscr{t}}

\def\dgl{\Delta_\mathrm{GL}}
\def\dck{\Delta_\mathrm{CK}}
\def\wdck{\widetilde{\Delta}_\mathrm{CK}}

\def\dckRed{\widetilde\Delta_\mathrm{CK}}

\def\e{\mathrm{e}}
\def\eck{\e_\mathrm{CK}}
\def\egl{\e_\mathrm{GL}}

\def\Log{\mathrm{Log}}
\def\Exp{\mathrm{Exp}}
\def\qLog{\widetilde{\mathrm{L}}\mathrm{og}}
\def\qExp{\widetilde{\mathrm{E}}\mathrm{xp}}
\def\a{\alpha}
\def\b{\beta}
\def\c{\gamma}

\newcommand{\doubleBr}[1]{(\!(#1)\!)}

\usepackage{euscript}

\usepackage{subfiles}

\begin{document}
	\maketitle
	\begin{abstract}
    Branched rough paths, defined as paths with values in the character group of the Connes-Kreimer Hopf algebra $\mathcal{H}_\mathrm{CK}$, constitute integration theories that may fail to satisfy the usual integration by parts identity.
    Using known results on the primitive elements of $\mathcal{H}_\mathrm{CK}$ we can view it as a commutative cofree Hopf algebra (i.e.\ a commutative $\mathbf{B}_\infty$-algebra) and thus write an explicit change-of-variable formula for solutions to rough differential equations. This formula restricts to the well-known It\^o formula in the very special case of semimartingales.
    In addition, we establish an isomorphism between $\mathcal{H}_\mathrm{CK}$  and the shuffle algebra over its primitives, which extends Hoffman's exponential for the quasi-shuffle algebra, and can therefore be viewed as a far-reaching generalisation of the usual Itô-Stratonovich correction formula for semimartingales.
    Indeed, this can be stated as a characterisation of the algebra structure of any commutative $\mathbf{B}_\infty$-algebra.
    Compared to previous approaches, this transformation has the key property of being natural in the decorating vector space.
    We study the one-dimensional case more closely, by introducing the branched analogue of the Kailath-Segall polynomials and Dol\'eans-Dade exponential, and conclude with some examples of branched rough path lifts of a stochastic process which are not quasi-geometric.
	\end{abstract}
	\tableofcontents

	\maketitle
\tableofcontents

\section*{Introduction}\label{sec:intro}
One of the cornerstones of stochastic analysis is Itô's change-of-variable formula \cite{Ito51b}. Given a continuous semimartingale \(Y\) and a smooth function \(\varphi\), it tells us how to express \(\varphi(Y)\) in terms of Itô integration against \(Y\) and Riemann-Stieltjes integration against the quadratic variation path \([Y]\):
\begin{equation}\label{eq:itoFormula}
\varphi(Y_t) - \varphi(Y_s) = \int_s^t \D \varphi(Y_u)\,\dif Y_u + \frac 12 \int_s^t \D^2 \varphi(Y_u)\,\dif [Y]_u
\end{equation}
This result is arguably what elevates the status of Itô's theory to that of a \say{calculus}, albeit one that does not satisfy the same identities as ordinary calculus, as exemplified by the above identity. The Itô formula provides the means to carry out many computations of interest in probability theory, such as those involving (conditional) expectations. For example, if \(Y\) satisfies the stochastic differential equation (SDE)
\begin{equation}\label{eq:SDE}
	\dif Y_t = \sigma(Y_t,t)\,\dif W_t + \mu(Y_t,t)\,\dif t
\end{equation}
where \(W\) is a (multidimensional) Brownian motion, the Itô formula can be written as
\begin{equation}\label{eq:generator}
\begin{split}
\varphi(Y_t) - \varphi(Y_s) &= \int_s^t \D\varphi(Y_u)\sigma(Y_u,u)\,\dif W_u + \int_s^t \mathcal L\varphi(Y_u,u)\,\dif u \\
\text{with}\quad \mathcal L\varphi(y,u) &\coloneqq \D\varphi(y) \mu(y,u) + \frac 12  \D^2\varphi(y)  \mathrm{Tr}(\sigma_\cdot \sigma_\cdot^\intercal)(y,u) \,.
\end{split}
\end{equation}

While the Itô integral is often preferred in probability and finance, thanks to its martingale-preserving properties, the Stratonovich integral is often preferred in physics and geometry, since it does satisfy the same identities as ordinary calculus (e.g.\ the Itô formula but without the second-order correction term). The two integrals are respectively defined by left endpoint and midpoint (or trapezoidal) Riemann-Stieltjes approximations:
\[
\int H \dif Y  \eqqcolon L^2\!\lim \sum_{[s,t]} H_s Y_{s,t},\qquad \int H \circ \dif Y  \eqqcolon L^2\!\lim \sum_{[s,t]} H_{\frac{s+t}{2}} Y_{s,t}
\]
where \(Y_{s,t} \coloneqq Y_t - Y_s\), \(H\) is an adapted and continuous integrand, and the limits are taken in \(L^2\) along a partition with vanishing mesh size.

\begin{figure}[!ht]
  \begin{center}
    \hfill
    \begin{minipage}{0.475\textwidth}
      \includegraphics[width=0.8\linewidth]{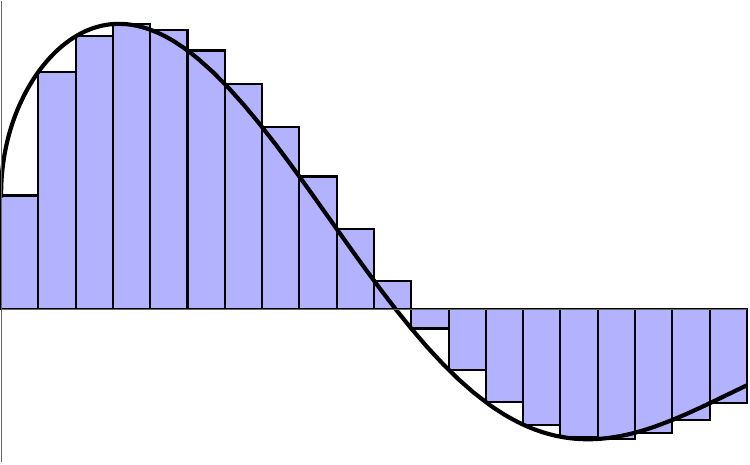}
    \end{minipage}
    \hfil
    \begin{minipage}{0.475\textwidth}
      \includegraphics[width=0.8\linewidth]{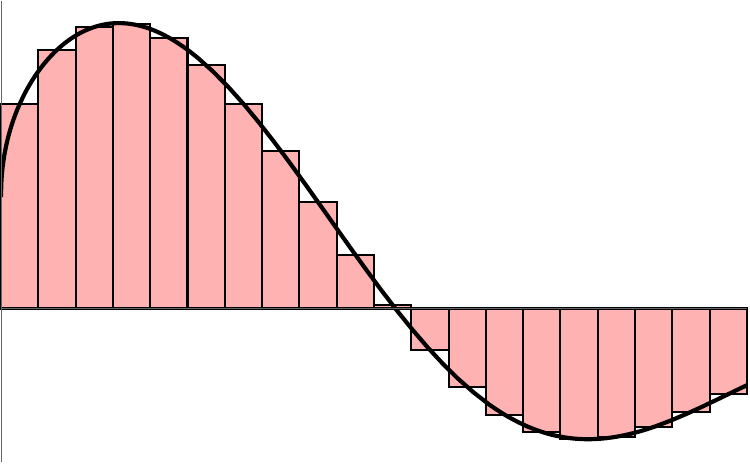}
    \end{minipage}
    \hfill
\end{center}
	\caption{Visual representation for left endpoint and midpoint Riemann-Stieltjes sum approximations.}
\end{figure}
Despite the Itô and Stratonovich integrals being different, they are related by the identity \(\int H \circ \dif Y = \int H \dif Y + \frac 12 [H,Y]\), the bracket denoting the quadratic covariation between \(H\) and \(Y\). In particular, this means that if \(Y\) satisfies the SDE \eqref{eq:SDE}, it equivalently satisfies the Stratonovich SDE
\begin{equation}\label{eq:itostrat}
  \dif Y_t = \sigma(Y_t,t)\,{\circ \dif W_t} +  \bigg[\mu(Y_t,t) - \frac 12 \D\sigma(Y_t,t) \cdot \sigma(Y_t,t) \bigg]\,\dif t
\end{equation}
This relates two a priori-distinct notions of stochastic differential equation, each with different advantages and drawbacks, at the sole cost of modifying the coefficients.

While Itô and Stratonovich SDEs constitute the main examples of random dynamical systems perturbed by instantaneous noise, they are by no means the only ones that can be conceived. A much more versatile notion of \emph{controlled differential equation} involves considering a driving process \(X \in C([0,T],V)\) (with \(V\) a finite-dimensional vector space) and an equation
\begin{equation}\label{eq:CDE}
\dif Y_t = F(Y_t)\,\dif X_t
\end{equation}
with \(F \in C^\infty(\mathcal L(V,W))\), \(\mathcal L\) denoting linear maps. When \(X\) is smooth, of course, one may substitute \(\dif X_t = \dot X \dif t\) above, but otherwise the equation must be understood differently. When \(X\) is of bounded \(\rho\)-variation with \(\rho<2\), such equation was first considered by Lyons \cite{Lyo94}, using Young's notion of integral \cite{You36}. When \(X\) is only of bounded \(2 < \rho\)-variation, however, Lyons's fundamental insight \cite{Lyo98} is that the path \(X\) must be augmented by additional information to give the equation unique, well-posed and robust meaning. Initially, this superstructure---a \emph{rough path}---involved postulating iterated integrals \(\int_{s < u_1 < \ldots < u_n < t} \dif X_{u_1} \otimes \cdots \otimes \dif X_{u_n}\) for all \(n \leq \p\): this is usually done, when \(X\) is random, via probabilistic notions of convergence, subject to certain constraints, including a first-order integration by parts identity. Rough path theory is able to handle a wide variety of signals, such as Gaussian and Markov processes, including the special case of Stratonovich SDEs, whilst not relying on probabilistic notions of convergence to give meaning to \eqref{eq:CDE} once the rough path is defined. This clean separation between probability theory and analysis dovetails nicely with the fact, famously observed by F\"ollmer \cite{Foll81}, that It\^o's change of variable formula is defined deterministically along a suitable sequence of partitions, and pathwise It\^o formulae have continued to attract interest  since, see e.g. \cite{CP19}.

Lyons's theory was extended by Gubinelli in \cite{Gub10}, who relaxed the algebraic requirements on the rough path,
making it possible for the resulting integral not to satisfy the same identities as ordinary calculus. This involved
requiring that the rough path not only contain the information representing linear iterated integrals but also
integrals of products (of integrals of products\ldots and so on). Such an object is readily viewed as a functional on
the Connes--Kreimer Hopf algebra \(\hck\), introduced at the end of the 20th century in the context of renormalization
in quantum field theory \cite{CK98}.
Unterberger \cite{MR2657813} independently considered rough paths defined in this way around the same time. \emph{Branched rough paths}
are general enough to reproduce Itô's theory of integration, and arguably provide the most general framework for what a
theory of path integration might look like. What was initially missing from Gubinelli's theory was a change-of-variable
formula akin to \eqref{eq:itoFormula}, and an explanation of the relationship between branched and the
previously-defined geometric rough paths. These two questions were addressed in Kelly's PhD thesis \cite{Kelly}, the
latter jointly with Hairer \cite{HK15}: for the former, an Itô-type formula is possible, at the cost of enlarging the
original branched rough path; for the latter, every equation driven by a branched rough path can be expressed in terms
of a geometric one, again, at the cost of enlarging (in a different way) the original path and choosing a geometric
rough path above it. These two results, while of great significance, are vulnerable to the objection that the extensions
they consider are chosen in a highly non-canonical fashion. In contrast, the correction needed
for the usual Itô formula, namely the quadratic variation \([Y]\), is defined intrinsically given the path \(Y\) (even algebraically, as the defect in the integration by parts formula: \([Y]^{\a\b} = Y^\a Y^\b - \int Y^\a \dif
Y^\b - \int Y^\b \dif Y^\a\)). This was addressed, for the second point above, by Boediharjo and Chevyrev: in
\cite{BC2019} they utilise results from Foissy \cite{Foi2002} and Chapoton \cite{Ch10} to observe that the mapping from
branched to geometric rough paths can be achieved in an entirely algebraic fashion, thus avoiding the troublesome choice
of the rough path lift. The isomorphism of which they show the existence is still, however, not explicitly identified
and computed, and many (indeed infinitely many) choices for it exist, each corresponding to a free algebraic basis of the
Grossman--Larson Hopf algebra \(\hgl\), introduced in \cite{GL89} to study compositions of differential operators. The relation between rooted trees and differential operators goes back to Cayley \cite{C1857}, and was extensively employed by Butcher (and earlier by Merson \cite{Merson57}) in relation to the analysis of Runge--Kutta schemes \cite{Butcher64} in what is now known as B-series.

Our contributions are as follows. In \Cref{sec:algebra} we view the decorated versions of $\hck$ and $\hgl$ as functors from  the category of finite-dimensional vector spaces \(\cat{Vec} \) to the category of graded, connected Hopf algebras \( \cat{gcHopf}\). We use Foissy's result \cite{Foi2002}, building on work of Broadhurst and Kreimer \cite{BK2000}, that \(\hck\) 
is a commutative \(\binf\)-algebra \cite{F2017}, and revisit Hoffman's duality result with \(\hgl\) in this framework. A thorough review of the Oudom-Guin result on the universal enveloping algebra of a pre-Lie algebra \cite{OG08} and its implications in the study of differential operators concludes the first section.

In \Cref{sec:brps} we consider branched rough paths defined on \(\hck\) and identify primitives in \(\hck\)---call the space of these \(\cP\)---as those elements which index paths: these form the \say{intrinsic correction terms} to the usual integration-by-parts identity; a similar point of view was taken in the context of renormalisation \cite{CQRV02}.
Integration against these is what makes the change-of-variable formula, \Cref{thm:ito} possible, the first of the two main results of this paper. Here we consider an equation similar in spirit to \eqref{eq:CDE}, but driven by \(X\) together with its correction terms, and thus strictly more general than those considered by previous authors. The Itô formula is achieved by lifting \(F\) to a Hopf algebra homomorphism \(\mathbf F\) from \(\hgl\) to the Hopf algebra of differential operators over the state space of the solution, which yields our \say{rough generator} (similar to \eqref{eq:generator}) to be integrated. When restricted to equations \say{without drift} (i.e.\ only driven by the original \(X\) and not its corrections), this yields a construction which can be compared to Kelly's, but that does not rely on non-canonical rough path lifts and admits a more precise description of its coefficients as the Oudom-Guin extension of a pre-Lie morphism. Our Itô formula is particularly tractable in the special instance of quasi-geometric rough paths---those rough paths defined on Hoffman's quasi-shuffle algebra \cite{Hof00}---in which case it extends previous work of one of the authors \cite{Bel20}.
Writing the quasi-shuffle algebra as a quotient of \(\hck\), and leveraging once again the algebra introduced previously to identify a section of this quotient map, leads to an explicit criterion for when a branched rough path is (quasi-)geometric, filling a gap in the literature.

In \Cref{sec:itostrat} we begin by introducing the Eulerian idempotent of a Hopf algebra, and use it to define an
isomorphism \(\Log \colon \hck \to \bigshuffle(\cP)\), the latter being the shuffle Hopf algebra over \(\cP\). This
constitutes our second main result, \Cref{thm:iso}. Since \(\Log\) is defined solely in terms of operations intrinsic to
the \(\binf\)-algebra \(\hck\), it can be viewed as a natural transformation between the functors \(\hck\) and
\(\bigshuffle(\cP)\). The definition and proof that this map is a isomorphism holds for general commutative $\binf$-algebras, and explicit
closed-form and recursive expressions for the inverse \(\Exp\) are identified. Indeed, one may view this result as a
structure theorem for (graded connected) cofree commutative Hopf algebras.
Since the shuffle Hopf algebra is the one
on which geometric rough paths are defined, this immediately yields a recipe for obtaining, given a branched rough path
\(\bX\), a geometric one \(\overline \bX\) which carries the same information. By this we mean that an equation driven
by \(\bX\) can be re-expressed as an equation driven by \(\overline \bX\) whose coefficients are explicitly computable
in terms of pre-Lie products of coefficients obtained through the Hopf morphism \(\mathbf F\) and the adjoint (under the
\(\hck\)-\(\hgl\) duality) of \(\Log\). We proceed by showing how \(\Exp\) extends Hoffman's exponential (after which it
is named) between the quasi-shuffle and shuffle algebras, originally introduced to study algebraic relations between multiple zeta
values \cite{Hof00}. This isomorphism can itself be regarded as an extension of the It\^o-Stratonovich correction formula for iterated integrals of Brownian motion, first written out by Hu and Meyer \cite{HM88}, and extended to more general processes several times since. We end the section with a comparison to the previous approaches of Hairer-Kelly and Boediharjo-Chevyrev, providing further motivation for our work. More precisely, the dual of our isomorphism may be viewed as a special case of the latter, the general case of which does not however guarantee naturality.

In the final section, \Cref{sec:one}, we restrict our attention to the more tractable case in which the underlying
vector space is one-dimensional (the \say{undecorated} setting). While there is not much to say about one-dimensional
geometric rough paths (there is only one over each trace path, given by powers \(X^n_{s,t}/n!\)), one-dimensional
branched rough paths are already interesting: their linear iterated integrals over a single primitive path component can
be expressed as polynomials in the correction terms. This fact, which is proved using Eulerian and Dynkin idempotents,
extends the well-known Kailath-Segall polynomials \cite{Segall1976OrthogonalFO} to the branched setting. Our proof of such identities is new even in the classical semimartingale setting, in which such polynomials are well known to be the orthogonal polynomials for Gaussian processes (Hermite polynomials) and Poisson processes (Charlier polynomials), of great significance in stochastic analysis. We end the paper by laying out a framework, leveraging the algebraic framework presented in the previous sections, through which branched rough paths with given correction terms may be constructed. Motivated by some results concerning Itô formulae in law \cite{Nourdin09}, we consider rough paths above scalar fractional Brownian motion with Hurst parameter \(1/4\): two out of four such examples are \say{truly branched}, i.e.\ not quasi-geometric.

We believe that many of the modern challenges presented by stochastic analysis can be met by a careful and parsimonious treatment of the algebraic structures involved, with an emphasis on natural structure and functoriality. We look forward to further developments of the theory in this spirit, particularly in the direction of SPDEs and regularity structures.

\subsection*{Acknowledgments}
CB and NT were partially supported by the DFG Research Unit FOR2402 at the start of this project. EF is supported by
UKRI EPSRC Programme Grant EP/S026347/1. NT is also funded by the Deutsche Forschungsgemeinschaft (DFG, German Research
Foundation) under Germany's Excellence Strategy – The Berlin Mathematics Research Center MATH+ (EXC-2046/1, project ID:
390685689) and CRC/TRR 388 "Rough Analysis, Stochastic Dynamics and Related Fields" – Project ID 516748464.
We would like thank Frederic Patras, Lo\"ic Foissy and Kurusch Ebrahimi-Fard for fruitful conversations which took place at the conference \say{Structural Aspects of Signatures and Rough Paths} held at the Centre for Advanced Study (CAS) of Oslo, where a preliminary version of these results was presented.
We are also grateful to the Mathematische Forschungsinstitut Oberwolfrach for their warm hospitality during the final stage of the writing of this manuscript.

\section{Algebraic preliminaries}\label{sec:algebra}
	\subsection{The Connes-Kreimer Hopf algebra over a vector space}\label{subsec:CK} 
	We begin by introducing the Connes-Kreimer Hopf algebra over a vector space as an algebra of tensors. We assume the reader is familiar with the original construction \cite{CK98} (see for example the survey \cite{Foi13}), as well as with the basics of Hopf algebras \cite{Man06,Und15}. 
	
	Let \(\scrH\) be the set of non-planar rooted forests (including the empty forest \(\mathbf{1}\); we denote \(\scrH_+ \coloneqq \scrH \setminus \{\mathbf{1}\}\)), and \(\scrT\) its subset of non-planar rooted trees.
  We will denote elements of \(\scrH\) by \(\scrf,\scrg,\ldots\) and elements of \(\scrT\) by \(\mathscr s, \mathscr t,\ldots\), and for \(\mathscr f \in \scrH\) we let \([\mathscr f] \in \scrT\) denote the tree obtained by joining every root in \(\mathscr f\) to a new vertex, the new root.   We denote by \(\{\scrf\}\) the underlying set of vertices of \(\scrf \in \scrH\).  This set inherits a natural partial order from its forest structure.  Given a finite set \(A\), we denote by \(\mathbb S_A\) the permutation group of \(A\), \(\mathbb S_k \coloneqq \mathbb S_{\{1,\ldots,k\}}\).
  We define \(\mathbb S_\scrf\) to be the group of order-preserving permutations of \(\{\scrf\}\) (the order coming from the forest structure), which we may view as a subgroup of the ordinary permutation group \(\mathbb S_{\{\scrf\}}\).
  The set \(\mathbb S_\scrf\) has the following inductive description: \(\mathbb S_{[\scrf]} = \mathbb S_\scrf\), and for \(\mathscr f = \mathscr t_1^{k_1} \cdots \mathscr t_n^{k_n}\) with \(\mathscr t_i \neq \mathscr t_j\) for \(i \neq j\), \(\sigma \in \mathbb S_{\{\scrf\}}\) belongs to \(\mathbb S_\scrf\) if we can write
	\begin{equation}\label{eq:indS}
		\begin{split}
	\sigma = ((\sigma_1^1 \sqcup \ldots \sqcup \sigma^1_{k_1}) \circ \tau_1) \sqcup \ldots \sqcup ((\sigma_1^n \sqcup \ldots \sqcup \sigma^n_{k_n}) \circ \tau_n) \\
	\text{for some } \sigma_i^j \in \mathbb S_{\scrt_i} \text{ acting on the $j^\text{th}$ copy of $\scrt_i$}, \quad \tau_i \in \mathbb S_{k_i}
	\end{split}
	\end{equation}
	i.e. \(\sigma\) acts disjointly on the vertices of each individual tree \(\scrt_i\) by an element of \(\mathbb S_{\scrt_i}\) and then disjointly permutes the sets of trees that are all equal. We denote the cardinality of $S_{\mathscr f}$ by $\varsigma(f)$.
	
	Given a (real) finite-dimensional vector space \(U\), we will now define a way of \say{labelling} the vertices of a forest without having to choose a basis of \(U\).  This and similar constructions have already appeared in the literature: \cite{CW17} (motivated by the need to consider \(U\) infinite-dimensional), and \cite[Ch.\ 5]{Lod12} (using \(\mathbb S\)-modules); we choose to adopt this notation in order to develop the theory in a coordinate-free manner.
    For \(\mathscr f = \mathscr t_1^{k_1} \cdots \mathscr t_n^{k_n}\), we define recursively
    \begin{equation}\label{eq:Vf}
      U^{\boxtimes\mathbf{1}} \coloneqq \mathbb R, \qquad U^{\boxtimes\mathscr f} \coloneqq (U^{\boxtimes\mathscr t_1})^{\boxdot k_1} \boxtimes \cdots \boxtimes (U^{\boxtimes\mathscr t_n})^{\boxdot k_n}, \qquad U^{\boxtimes[\mathscr f]} \coloneqq U^{\boxtimes\scrf} \boxtimes U\,,
    \end{equation}
    where \(\boxtimes\) is the usual tensor product among vector spaces and \(\boxdot\) is the symmetric tensor product (we prefer these symbols to \(\otimes\), \(\odot\), since the latter will be used more extensively later for external tensor products).
    
  The presence of \(\boxdot\) ensures that \(U^{\boxtimes \mathscr f}\) encodes all the symmetries of the forest \(\mathscr f\) not already present in the non-planarity of the undecorated forest. From \eqref{eq:indS} it follows that there is a canonical isomorphism between this definition that of \cite{CW17} as a quotient by the action of \(\mathbb S_\scrf\)
  \begin{equation}\label{eq:quotient}
    U^{\boxtimes \scrf} \cong  U^{\boxtimes \{\scrf\}}\Big/\mathbb{ S}_{\scrf}\,;
  \end{equation}
  both descriptions will be useful. At the two extremes, denoting \(\mathscr l_n\) the totally ordered tree with \(n\) vertices (called the \emph{\(n\)-ladder}) and \(\mathscr r_n\) the product of \(n\) disjoint roots, we have \(U^{\boxtimes \mathscr l_n} = U^{\boxtimes n}\) and \(U^{\boxtimes \mathscr r_n} = U^{\boxdot n}\). A more generic example would be
  \[
    U^{\boxtimes\, \Forest{[[][[]][[]]]}} = U \boxtimes (U \boxtimes U)^{\boxdot 2} \boxtimes U
  \,.\]

	Elementary tensors in \(U^{\boxtimes\scrf}\) should be thought of as decorations of \(\scrf\) with elements of \(U\); we denote \(\scrH(U)\) such generators and continue to use cursive letters to denote such decorated forests and trees, while we use Greek letters to denote elements of \(U\). For example
	\[
	\Forest{[\alpha[\beta][\gamma[\delta]][\epsilon[\zeta]]]} = \beta \boxtimes ((\delta \boxtimes \gamma) \boxdot (\zeta \boxtimes \epsilon)) \boxtimes \alpha \ \in U^{\boxtimes\, \Forest{[[][[]][[]]]}}\,.
	\]

	When considering decorated forests, we will use the notation \([\scrf]_\alpha\) to refer to the same operation as \([\scrf]\) in which we decorate the root with \(\alpha\).
  \begin{definition}
   For any finite-dimensional vector space \(U\) we define the vector spaces, graded by \emph{weight} of the forest (i.e.\ number of vertices, denoted \(|\scrf|\))
    \begin{equation}\label{eq:FV}
      \cH(U) \coloneqq \mathrm{span}_\mathbb{R} \scrH(U) = \bigoplus_{\mathscr f \in \scrH} U^{\boxtimes \scrf}, \qquad \cT(U) \coloneqq \mathrm{span}_\mathbb{R} \scrT(U) = \bigoplus_{\mathscr t \in \scrT} U^{\boxtimes \mathscr t}
    \end{equation}
    and note by $\cH$ the vector space freely generated by unlabelled forests (i.e. when  \(U = \bbR\)). Similarly, we use the notations \(\cH_+(U) \coloneqq \mathrm{span}_\mathbb{R} \scrH_+(U)\) and $\cH_+$.
  \end{definition}

  It is also helpful to introduce \emph{proper forests} \(\scrF \coloneqq \scrH \setminus \scrT\), i.e.\ forests that are not trees (this includes the empty forest), and their span \(\cF\), so that \(\cH = \cT \oplus \cF \). Let \(\underline{\mathrm{Vec}}\) denote the category of finite-dimensional \(\mathbb R\)-vector spaces. Thanks to functoriality of \(\boxtimes\) and \(\boxdot\), \(\cH\) can be viewed as a covariant functor \(\underline{\mathrm{Vec}} \to \underline{\mathrm{Vec}}\). Let \(\bigtimes(B) \coloneq \bigsqcup_{n \in \mathbb N} B^n\) denote the (disjoint) union of \(\mathbb N \ni n\)-fold Cartesian products (including the empty one) of a set \(B\), and \(\bigotimes ( \, \cdot \,) = \bigoplus_{n \in \mathbb N} (\, \cdot \,)^{\otimes n}\) the tensor algebra functor \(\underline{\mathrm{Vec}} \to \underline{\mathrm{Vec}}\) (note how we are now using the symbol \(\otimes\) for external tensor product). The next proposition guarantees that certain operations performed on individual vertices of forests are promoted to natural transformations.
  \begin{proposition}\label{prop:functor}
    Let $\Phi$ be a collection of real numbers indexed by bijections of vertex sets
    \[ \Phi = \big\{ \phi_b \in \mathbb R \colon \{b \colon \{\scrf_1\} \sqcup \cdots \sqcup \{\scrf_m\} \leftrightarrow \{\mathscr g_1\} \sqcup \ldots  \sqcup \{\mathscr g_n\} ,\ (\scrf_1,\ldots,\scrf_m), (\scrg_1,\ldots,\scrg_n) \in \bigtimes(\scrH) \}\big\}	\]
    with the property that $\phi_b = \phi_{b \circ \sigma}$ for any $\sigma \in \mathbb S_{\scrf_1} \times \cdots \times \mathbb S_{\scrf_m}$ (viewed as a subgroup of $\mathbb S_{\{\scrf_1\} \sqcup \ldots \sqcup \{\scrf_m\}}$). For any finite-dimensional vector space $U$ we set
    \begin{equation}\label{eq:phiU}
      \begin{split}
        \Phi_U(\scrf_1 \otimes \cdots \otimes \scrf_m)= \sum_{b \colon \{\scrf_1\} \sqcup \cdots \sqcup \{\scrf_m\} \leftrightarrow \{\mathscr g_1\} \sqcup \ldots  \sqcup \{\mathscr g_n\}} \phi_b \mathscr g_1 \otimes \cdots \otimes \mathscr g_n\,,
      \end{split}
    \end{equation}
    where each $\scrf_i$ is a decorated forest (an elementary tensors in $U^{\boxtimes \scrf_i}$) and the undecorated forests $\scrg_j$ are given the decorations corresponding to the images of their vertices through the bijections $b$. Then $\Phi$ defines a natural endomorphism of the functor $\bigotimes \circ \cH \colon \cat{Vec} \to \cat{Vec}$.
  \end{proposition}

    We comment on the definition of $\Phi$.
    The fact that $\Phi$ is indexed by bijections is to be interpreted by saying that it corresponds to performing certain weight-preserving operations on an ordered collection of forests, e.g.\ cutting and grafting of edges, and that $\phi_b$ is the coefficient of $(\scrg_1,\ldots,\scrg_n)$ for this operation applied to $(\scrf_1,\ldots,\scrf_m)$, where the vertices of the forests $\scrg_j$ correspond to specific vertices of the forests $\scrf_i$ through $b$.
    It would not be possible to lift an element of $\mathrm{End}(\bigotimes (\cH(\bbR)))$ to one of $\mathrm{End}(\bigotimes ( \cH(U)))$ without knowing the precise correspondence between vertices: for example, knowing that the former maps $\Forest{[[]]} \mapsto \Forest{[[]]}$ does not, on its own, determine whether the induced one should map $\Forest{[\beta[\alpha]]} = \alpha \boxtimes \beta$ to itself or to $\Forest{[\alpha[\beta]]} = \beta \boxtimes \alpha$.
    We are considering endomorphisms of the whole tensor algebra over $\cH$ for maximum flexibility, but often we will only need to consider $\Phi$ to be $\cH^{\otimes 2} \to \cH$ or similar, which can be obtained by the proposition by setting many of the coefficients to zero.
    The map $b$ is required to be a bijection, since injections may leave certain labels undefined and surjections may lead to ill-defined maps (for example assigning the sole non-zero coefficient to the map $\Forest{[[]]} \to \Forest{[]}$ that outputs the root does not lead to a well-defined map $U \boxtimes U \ni \alpha \boxtimes \beta \mapsto \beta \in U$); note also that allowing for multi-valued functions would violate linearity (e.g.\ duplicating a vertex would lead to $\Phi_U(\lambda \Forest{[\alpha]}) = \lambda^2 \Forest{[\alpha]} \Forest{[\alpha]}$).
    The invariance requirement rules out maps that would not be well-defined on non-planar forests, e.g.\
    $\Forest{[\alpha]} \Forest{[\beta]} \mapsto \Forest{[\alpha]} \otimes \Forest{[\beta]}$ (while $\Forest{[\alpha]}
    \Forest{[\beta]} \mapsto \Forest{[\alpha]} \otimes \Forest{[\beta]} + \Forest{[\beta]} \otimes \Forest{[\alpha]}$ is
    acceptable).

  \begin{proof}[Proof of \Cref{prop:functor}]
    First of all, \eqref{eq:phiU} is a well defined linear map $\bigotimes (\mathcal H(U)) \to \bigotimes (\mathcal H(U))$, as a direct sum of its restrictions to $U^{\boxtimes \scrf_1} \otimes \cdots \otimes U^{\boxtimes \scrf_m}$. Indeed, viewing each of these as a map from $U^{\boxtimes \{\scrf_1\}} \otimes \cdots \otimes U^{\boxtimes \{\scrf_m\}}$, which it is equal to a sum of maps induced by permutations and thus linear, it passes to the quotient (where each $U^{\boxtimes \scrf_i}$ is given the description of \eqref{eq:quotient}) thanks to $\mathbb S_{\scrf_1} \times \cdots \times \mathbb S_{\scrf_m}$-invariance.  Next, we must show that for any $\mathbb R$-linear map $\theta \colon V \to W$ the square
    \[\begin{tikzcd}
      \textstyle\bigotimes (\cH(V)) \arrow[r,"\Phi_V"] \arrow[d,"\bigotimes(\cH (\theta))",swap] & \bigotimes (\cH(V)) \arrow[d,"\bigotimes (\cH (\theta))"] \\
      \textstyle\bigotimes (\cH(W) ) \arrow[r,"\Phi_W"]  & \bigotimes (\cH(W))
    \end{tikzcd}
  \]
  commutes. We have
  \begin{align*}
  &\tens(\mathcal H(\theta)) \circ \Phi_V (\scrf_1 \otimes \cdots \otimes \scrf_m) \\
    ={} &\sum_{b \colon \{\scrf_1\} \sqcup \cdots \sqcup \{\scrf_m\} \leftrightarrow \{\mathscr g_1\} \sqcup \ldots  \sqcup \{\mathscr g_n\}} \phi_b \cH(\theta)(\mathscr g_1) \otimes \cdots \otimes \cH(\theta)(\mathscr g_n) \\
    ={} &\Phi_W \circ \tens(\mathcal H(\theta))  (\scrf_1 \otimes \cdots \otimes \scrf_m)
  \end{align*}
  where the second identity follows from the fact that $\mathcal H(\theta)$ and $b$ are interchangeable, by functoriality of (symmetric) tensor products. Since the identity holds on a generating set of $\bigotimes (\mathcal H(V))$ we conclude by linearity.
\end{proof}
This result makes it possible to regard \(\cH\) and spaces defined in terms of it as undecorated, so that operations on
it given of the form above, indexed by precise bijections between the vertex sets, automatically induce maps on the decorated spaces, natural in the vector space of decorations. This will be made precise once and for all in \Cref{cor:natural} for the main structures of interest to us.
	Recall that the \emph{Connes-Kreimer Hopf algebra} \(\mathcal H_\mathrm{CK} \coloneqq (\cH, \cdot, \dck)\) is the free abelian algebra over \(\scrT\), with coproduct \(\dck\) defined in terms of admissible cuts.
  By applying \Cref{prop:functor} with \(\Phi = \cdot, \Delta\), we obtain a graded connected Hopf algebra \(\mathcal H_\mathrm{CK}(U) \coloneqq (\cH(U), \cdot, \dck)\).
  Let \(\underline{\mathrm{gcHopf}}\) denote the category of graded, connected Hopf algebras, and recall that the forgetful functor from this category to that of graded, connected bialgebras is an equivalence, i.e.\ each graded connected bialgebra has a unique antipode, and these are automatically preserved by maps of graded, connected bialgebras, see \cite{tak71}.
  This is what makes it possible to treat the category \(\cat{gcHopf}\) without ever mentioning antipodes.

	We now introduce the secondary structure on \(\mathcal H_\mathrm{CK}(U)\) which will prove fundamental to our goals: cofree coalgebras.
  We refer the \cite[Section 1.2.]{Lod12} for a complete reference on the topic.
  In what follows, we will always consider graded, connected coassociative coalgebras \((C, \Delta, \varepsilon)\).
  In particular, as a vector space \(C\cong\bigoplus_{n\ge 0}C_n\), with \(C_0=\mathbb{R}\mathbf{1}_C\) of dimension one.
  The reduced coproduct will be denoted by
  \[
    \widetilde \Delta h\coloneqq\Delta h- \mathbf{1}_C \otimes h - h \otimes \mathbf{1}_C,
  \]
  and the space of primitive elements
  \[
    \mathrm{Prim}(C)\coloneqq\ker\widetilde{\Delta}=\{ p \in C : \Delta p = \mathbf{1}_C \otimes p + p \otimes \mathbf{1}_C \}.
  \]
  Both \(\Delta\) and \(\widetilde \Delta\) can be iterated, thanks to coassociativity, to yield linear maps \(\Delta^{n}, \widetilde \Delta^{n} \colon C \to C^{\otimes (n+1)}\) for \(n\geq 2\) defined recursively by the condition 
	\[\Delta^{n}= (\mathrm{id} \otimes\Delta^{n-1}) \circ \Delta, \quad \widetilde\Delta^{n}= (\mathrm{id} \otimes
  \widetilde\Delta^n)\circ \widetilde\Delta. \]
  We extend these operators by setting \(\Delta^{1} \coloneqq \Delta\), \(\Delta^{0} \coloneqq \text{id}\) and \(\Delta^{-1}(\mathscr{f}) \coloneqq
  \varepsilon(\mathscr{f})\mathbf{1}_C\), where \(\varepsilon\) is the counit of \(C\). We will frequently use (possibly sum-free) Sweedler notation
	\begin{align*}
	\Delta h &\eqqcolon \sum_{(h)} h_{(1)} \otimes h_{(2)} \eqqcolon h_{(1)} \otimes h_{(2)} \\
	\widetilde \Delta h &\eqqcolon \sum_{(h)} h^{(1)} \otimes h^{(2)} \eqqcolon h^{(1)} \otimes h^{(2)} \eqqcolon h' \otimes h''
	\end{align*}
and similar for higher order (reduced) coproducts.

\begin{definition}\label{def:cofree}
 \(C\) is cofree over a vector space \(P\) with projection \(\pi\colon C\to P\), sending \(\mathbf{1}_C\) to
  \(0\), if for all other graded and connected coalgebras \(D\) and linear maps \(\phi \colon D \to P\) sending
  \(\mathbf{1}_D\) to \(0\) there exists a unique graded coalgebra map \(\Phi \colon D \to C\) making the diagram
  \begin{equation}\label{eq:cofree}
    \begin{tikzcd}
      & C\arrow[d,"\pi"] \\
      D \arrow[r,"\phi",swap] \arrow[ur,"\Phi",dashed] & P
    \end{tikzcd}
  \end{equation}
  commute. We call \(\pi\) a \emph{cofreeness projection}.
\end{definition}

The canonical way of constructing a model of cofree coalgebra by considering the tensor algebra \(\bigotimes(P)\) with the deconcatenation coproduct 
\begin{equation}
  \Delta_\otimes (p_1 \cdots p_n) = \sum_{k = 0}^n (p_1  \cdots  p_k) \otimes (p_{k +1}  \cdots  p_n)\,.
\end{equation}
and \(\pi\) is the canonical projection \(\pi_1 \colon \bigotimes (P) \to P\), see \cite[Proposition 1.2.7]{Lod12}. It follows from the universal property \eqref{eq:cofree} that if \(C\) is cofree over \(P\) with projection \(\pi \colon C \to P\), \(\pi\) always restricts to a linear isomorphism \(\mathrm{Prim}(C) \to P\), i.e.\ the cofreeness map essentially consists of a projection of \(C\) onto \(\mathrm{Prim}(C)\). We emphasise that the structure of cofreeness depends on this projection: different maps \(\pi\) will yield different coalgebra isomorphisms \(C \cong \tens(P)\). Coalgebra morphisms from, and between, cofree coalgebras can be explicitly described in terms of their composition with the cofreeness projection.

\begin{proposition}\label{coalgebra_prop}
  If $C = \tens(P)$ and $\pi = \pi_1$  the unique coalgebra map $\Phi$ in \Cref{def:cofree} is given explicitly by \(\Phi(\mathbf{1_D})=\mathbf{1}_C\) and
  \begin{equation}
    \Phi:=  \sum_{n=1}^{\infty}\phi^{\otimes n} \circ \widetilde \Delta^{n-1}.
  \end{equation}
  In particular, if $D = (\tens(R), \Delta_\otimes)$ for some vector space $R$, we have 
  \begin{equation}
    \Phi(r_1 \otimes \cdots \otimes r_n) = \sum_{\substack{k = 1,\ldots, n \\ n_1,\ldots,n_k \geq 1 \\ n_1 + \ldots + n_k = n}} \phi(r_1 \otimes \cdots \otimes r_{n_1}) \otimes \cdots \otimes \phi(r_{n^{k-1} + 1}\otimes \cdots \otimes r_n) \,,
  \end{equation}
  where we let $n^h \coloneqq n_1 + \ldots + n_h$ and $r_1\,, \ldots r_n\in R$. Moreover, $\Phi$ is invertible if and only if $\phi$ is invertible.
\end{proposition}
\begin{proof}
  The first identity follows from \cite[Proposition 1.2.7]{Lod12}. The second follows from the explicit computation of $\Delta^{(k-1)}_\otimes$, see also \cite[Theorem 11.2]{Foi2002}, and \cite[Theorem 11.3]{Foi2002} for the claim of bijectivity.
\end{proof}

\begin{remark}
  The characterisation of cofreeness here presented works because we consider only connected graded coalgebras, whose coproduct is always conilpotent. Further generalisation to the general category of coalgebra would lead to a much more complicate universal object, see e.g. \cite[Theorem 6.4.1]{Sweedler1969}.
\end{remark}

	We now proceed to view \(\hck\) as cofree algebra. This is made possible thanks to the following operation:
	\begin{definition}[Natural growth, \cite{BK2000}]\label{def:T}
		Define \(\tee\colon\cH \otimes \cH \to\cH\), called \emph{natural growth} as the unique linear map satisfying 
		\(\scrf\tee \mathbf{1} = \scrf = \mathbf{1}\tee \scrf
		\) 
		on any forest $\scrf$ and for $\scrf \neq \mathbf{1} \neq \scrg$, we set $\scrf \tee \scrg$ as the sum over all the forests  obtained by grafting \(\mathscr{f}\) onto every vertex of \(\mathscr{g}\) and normalizing \footnote{Here we follow
    \cite{Foi02c} by including the normalizing factor, not present in the original definition.} by $|\mathscr{g}|$. In  each term, all roots of $\scrf$ must be grafted onto the same vertex of $\scrg$. 
\end{definition}
For instance, we have 
\[\Forest{[]}\Forest{[]}\tee\Forest{[[]]}=\frac12\left(
	\Forest{[[[][]]]}+\Forest{[[][][]]} \right)\,.	\]
 Note that \(\tee\) is not an associative operation, and when considering \(n\)-fold iterations of it, we will always compose it from the left, i.e.\
	\[
	h_1 \tee \ldots \tee h_n \coloneqq (h_1 \tee \ldots \tee h_{n-1}) \tee h_n.
	\]
The key property of this operation is obtained by relating \(\tee\) with the Connes-Kreimer primitive elements 
\[
\cP \coloneqq \mathrm{Prim}(\mathcal H_\mathrm{CK})\,.
\]
Namely, the map \(\tee\) restricted to \(\cH \otimes \pck\) satisfies the property
	\begin{equation}\label{eq:cocycle}
		\Delta_\mathrm{CK} (h \tee p) = (h \tee p) \otimes 1 + (\mathrm{id} \otimes (\cdot\, \tee p))\circ \Delta_\mathrm{CK}(h)\,.
	\end{equation}
	which is expressed by saying that for each \(p \in \cP\) the map \((\cdot\, \tee p) \colon \hck \to \hck\) is a 1-cocycle for the coproduct \cite[p.230]{CK98}; we will give this property an interpretation in terms of integration in \Cref{sec:brps} below. The following result uses the 1-cocycle property to establish the fundamental fact that \(\hck\) is cofree.
	\begin{theorem}[{\cite[Th\'eorème 48]{Foi02c}}]\label{thm:T}
		The following map (denoted by the same symbol)
		\begin{equation}
			\tee \colon \textstyle\bigotimes(\pck) \to \hck, \quad p_1 \otimes \cdots \otimes p_n \mapsto p_1 \tee \ldots \tee p_n
		\end{equation}
		is a coalgebra isomorphism.
	\end{theorem}
	The projection associated to this isomorphism is given by the following
	\begin{proposition}[{\cite[Theorem 9.6]{Foi2002}}]
	The projection $\pi\colon \hck \to\cP$ is given by the recursive formula
	\begin{equation}\label{recursive_pi}
    \pi =  \mathrm{id}-\tee\circ(\mathrm{id}\otimes\pi)\circ\widetilde{\Delta}_{\mathrm{CK}}
	\end{equation}
	\end{proposition}
	From now on, \(\pi\) will refer to the specific projection of the proposition above. We thus have two gradings on \(\hck\): 
	\begin{equation}\label{both-gradings}
		\hck = \bigoplus_{n \in \mathbb N} \cH^{(n)}, \qquad \hck = \bigoplus_{m \in \mathbb N}\cP^{\top m}
	\end{equation}
	where \(\cH^{(n)}\) denotes the space generated by forests of weight \(n\) and \(\cP^{\top m} \coloneqq \top(\cP^{\otimes m})\). Note that although the second grading is not locally finitely-dimensional (since \(\cP\) is infinite-dimensional), \(\cP\) can itself be graded into finitely-dimensional spaces \(\cP = \bigoplus_n \cP^{(n)}\), \(\cP^{(n)} \coloneqq \pi(\cH^{(n)})\). One may view elements of \(\cP^{\top m}\) corresponding to elementary tensor as \say{\(m\)-ladders decorated by primitives}, since \(\top(\Forest{[]}^{\otimes m}) = \mathscr l_m\). We will refer to such \(m\) as the \emph{degree of primitiveness}. Sometimes we will consider the associated filtrations, which we denote \(\cH^{n} \coloneqq \bigoplus_{k = 0}^n \cH^{(k)}\) in the first case. On the other hand, we do not have a particular notation for the filtration associated to the second grading, which coincides with the \emph{coradical filtration} of the coalgebra \((\hck, \dck)\), see \cite[section 11]{Foi2002}; note that this does not require cofreeness to be defined.
		Since \(\tee\) preserves weight, degree of primitiveness is bounded by weight, these two gradings can be combined as
	\begin{equation}\label{eq:twoGradings}
	\hck = \bigoplus_{0 \leq m \leq n} \pi_m(\cH^{(n)})
	\end{equation}
	where we denote \(\pi_m\) the canonical projection \(\hck \twoheadrightarrow \cP^{\top m}\).

We can use \Cref{coalgebra_prop} to express the inverse isomorphism \(\top^{-1}\colon  \hck\to \textstyle\bigotimes(\pck)\) in terms of \(\pi\):
\begin{equation}\label{eq:Tinverse}
\top^{-1} = \sum_{n=0}^{\infty}\pi^{\otimes n} \circ \widetilde \Delta_\mathrm{CK}^{n-1}
\end{equation}
from which, writing \(\top \circ \top^{-1} = \mathrm{id}\), it follows that
\begin{equation}\label{eq:pim}
\pi_m = \tee \circ \pi^{\otimes m} \circ \wdck^{m-1}
\end{equation}
It is also possible to obtain a non-recursive formula for \(\pi_m\), \(m \geq 1\) using the operator \(R_n\colon\cH^{\otimes n}\to\cH\), defined by
	\begin{align*}
		R_{n+1}(\scrf_0\otimes\dotsm\otimes\scrf_n) &\coloneqq \scrf_0\tee R_n(\scrf_1\otimes\dotsm\otimes\scrf_n)\,;
	\end{align*}
notice that  \(R_n|_{\cP^{\otimes n}}\neq \top|_{\cP^{\otimes n}}\), due to lack of associativity of the binary operation \(\tee\).
\begin{proposition}\label{prop:pirec}
For any finite-dimensional vector space $U$ the following non-recursive formulae for $\pi\colon \hck \to \pck$ and $\tee^{-1} \colon\hck  \to \textstyle\bigotimes(\pck)$ hold
		\begin{align*}
			\pi &= \sum_{k\ge 1}(-1)^{k+1}R_k\circ\wdck^{(k-1)}, \\
			\tee^{-1} &=\sum_{n\geq 0} \sum_{k\geq n}\sum_{k_1+\dotsb+k_n=k}(-1)^{k+n}(R_{k_1}\otimes\dotsm\otimes R_{k_n})\circ\wdck^{(k-1)}\,.
		\end{align*}		
\end{proposition}
\begin{proof}
	The first identity is obtained by iterating \eqref{recursive_pi}; the second by substituting it into \eqref{eq:Tinverse}.
\end{proof}

\begin{remark}
  Identity \eqref{eq:cocycle} implies that the pair of operations \(\top\) and \(\dck\)
  turn \(\hck\) into an infinitesimal bialgebra \cite[Definition 2.1.]{LR06}.
  The structure theorems and some of the formulae we obtain are a reflection of this fact.
\end{remark}

		\begin{example}[Change of coordinates for \(\hck^3\)]\label{expl:P3}
		 All elements of order \(1\), i.e.\ single vertices, are primitive. At level \(2\), we have
		 \[
		 \pi(\Forest{[\alpha]}\Forest{[\beta]}) = \Forest{[\alpha]}\Forest{[\beta]} - \Forest{[\beta[\alpha]]} - \Forest{[\alpha[\beta]]}
		 \]
		 so that
		 \[
		  \Forest{[\alpha]}\Forest{[\beta]} = \pi(\Forest{[\alpha]}\Forest{[\beta]}) + \Forest{[\beta[\alpha]]} + \Forest{[\alpha[\beta]]}, \qquad \Forest{[\beta[\alpha]]} = \Forest{[\alpha]} \tee \Forest{[\beta]}
		 \]
		 At level 3 we have
		\begin{align*}
			\pi(\Forest{[\alpha]}\Forest{[\beta]}\Forest{[\gamma]})&=\Forest{[\alpha]}\Forest{[\beta]}\Forest{[\gamma]}-\frac12\sum_{\sigma\in\mathbb{S}_3}\Forest{[\sigma(\alpha)]}\Forest{[\sigma(\beta)[\sigma(\gamma)]]}
			+\frac12\sum_{\sigma\in\mathbb{S}_3}\Forest{[\sigma(\alpha)[\sigma(\beta)[\sigma(\gamma)]]]}\\
			\pi\left(\Forest{[\alpha]}\Forest{[\c[\b]]}\right)&=\frac12\left( \Forest{[\alpha]}\Forest{[\c[\b]]}- \Forest{[\c]}\Forest{[\alpha[\b]]}\right)
			+\frac12\left( \Forest{[\alpha[\c][\b]]}-\Forest{[\c[\alpha][\b]]} \right) + \frac12\left( \Forest{[\alpha[\c[\b]]]}-\Forest{[\c[\alpha[\b]]]} \right)\,.
		\end{align*}
	In particular, the last term is antisymmetric in the lower two indices; notice how this is not visible in the 1-dimensional undecorated case: this is another example of how, using the notation of \Cref{prop:functor}, \(\Phi_U\) can be zero for a particular choice of \(U\) even though \(\Phi\) itself is not. More complicated relations can be expected at higher orders. From these computations we obtain the change of coordinates at level 3 (\(\mathscr l_3\) is already expressed in the \(\cP^{\top}\) grading)
	\begin{align*}
		\Forest{[\gamma[\alpha][\beta]]} &= \pi(\Forest{[\alpha]} \Forest{[\beta]}) \tee \Forest{[\gamma]} + \Forest{[\gamma[\beta[\alpha]]]} + \Forest{[\gamma[\alpha[\beta]]]} \\
	\Forest{[\a]}\Forest{[\b]}\Forest{[\c]} &= \pi(\Forest{[\a]}\Forest{[\b]}\Forest{[\c]}) + \sum_{\sigma \in \mathbb S_3} \Big( \Forest{[\sigma(\a)]} \tee \pi(\Forest{[\sigma(\b)]}\Forest{[\sigma(\c)]}) + \pi(\Forest{[\sigma(\b)]}\Forest{[\sigma(\c)]}) \tee \Forest{[\sigma(\a)]} \Big) + \sum_{\sigma\in\mathbb{S}_3}\Forest{[\sigma(\alpha)[\sigma(\beta)[\sigma(\gamma)]]]} \\
	\Forest{[\alpha]}\Forest{[\c[\b]]} &= \pi\left(\Forest{[\alpha]}\Forest{[\c[\b]]}\right) + \Forest{[\b]} \tee \pi(\Forest{[\a]}\Forest{[\c]}) + \pi(\Forest{[\a]}\Forest{[\b]}) \tee \Forest{[\c]} +  \Forest{[\c[\b[\a]]]} + \Forest{[\c[\a[\b]]]} + \Forest{[\a[\c[\b]]]}
	\end{align*}\qedhere
	\end{example}

	\begin{remark}[Bases of primitives]\label{rem:bases}
		Let \(B\) be a basis of \(U\), and denote \(\scrH(B)\) to be elementary tensors in \(\scrH(U)\) which are decorated with elements of \(B\), a basis of \(\cH(U)\). The previous example shows that \(\{ \pi_1(\scrf) \colon \scrf \in \scrH(B) \}\), which generates \(\cP(U)\), is not an independent set. Extracting a basis from this set would involve understanding the precise structure of \(\cP(U)\), which seems like a complex task (and perhaps similar, in flavour, to how Hall bases are defined \cite{Reu93}). Since in this article we are interested in describing natural operations between algebraic structures, and since such a basis would necessarily \say{break the symmetry}, we will avoid it and rather carry out all computations using forests natural maps. We will sometimes find it convenient to fix the basis \(B\) and perform computations in the induced forest basis \(\scrH(B)\).

    One can exactly compute the dimension of the span of primitive elements with exactly \(n\) nodes as a function of
    \(d\coloneqq\dim U\).
    Denoting \(h_n(d)\coloneqq\dim\cP^{(n)}(U)\), following the same argument as in \cite[Proposition 7.2]{Foi2002} one can show
    that the generating function
    \[
      H_d(t)\coloneqq\sum_{n\geq 1}h_n(d)t^n
    \]
    satisfies
    \[
      1-H_d(t) = \frac{1}{R_d(t)}
    \]
    where \(R_d\) is the generating function of the sequence \(\dim\hck^{(n)}(U)\).
    Using the coefficients found in \cite[Example 1.2]{Foi2021} we obtain, for the first few terms:
    \begin{align*}
      h_1(d)&=d,\\
      h_2(d)&=\frac{d(d+1)}{2},\\
      h_3(d)&=\frac{d(2d^2+1)}{3},\\
      h_4(d)&=\frac{d(9d^3+2d^2+3d+2)}{8},\\
      h_5(d)&=\frac{d(64d^4+20d^3-5d^2+6)}{30}.\qedhere
    \end{align*}
	\end{remark}

Cofree (graded, connected) Hopf algebras are intimately connected to the following notion for which we refer to
\cite[Ch.\ 2]{F2017}.
	\begin{definition}[commutative \(\mathbf{B}_\infty\)-algebra]
		A vector space \(P\) equipped with a map \(\lAngle  \argument, \argument\rAngle\colon \bigotimes(P)\odot
    \bigotimes(P)\to P\), for which we will use the shorthand $\lAngle  \argument, \argument\rAngle_{i,j} \coloneqq
    \left.\lAngle  \argument, \argument\rAngle\right|_{P^{\otimes i}\odot P^{\otimes j}}$, is said to be a commutative \(\mathbf{B}_\infty\)-algebra if the following conditions are satisfied:
		\begin{enumerate}
			\item for any \(k\ge 0\),
			\[
			\lAngle  \argument, \argument\rAngle_{k,0} =\begin{dcases}\mathrm{id}_P&k=1\\0&\text{otherwise}\end{dcases}
			\]
			\item for any tensors \(u,v,w\in \bigotimes(P)\) (and using Sweedler notation w.r.t.\ the unreduced $\Delta_\otimes^k$)
			\[
			\sum_{k\geq 1} \sum_{(u)^k,(v)^k}\lAngle \lAngle u_{(1)},v_{(1)}\rAngle \otimes \dotsm\otimes\lAngle u_{(k)},v_{(k)}\rAngle,w\rAngle=\sum_{k\ge 1}\sum_{(v)^k,(w)^k}\lAngle u,\lAngle v_{(1)},w_{(1)}\rAngle\otimes \dotsm\otimes\lAngle v_{(k)},w_{(k)}\rAngle\rAngle\,.
			\]
		\end{enumerate}
	\end{definition}
Commutative \(\mathbf{B}_\infty\)-algebra are intimately related to study of all possible commutative products on the tensor coalgebra.
	
\begin{theorem}[{\cite[Theorem 2]{F2017}}]
	Let \(\operatorname{Bialg}(P)\) be the set of commutative products on \(\bigotimes(P)\) compatible with
	\(\Delta_\otimes\) and \(\mathbf{B}_\infty(P)\) be the set of \(\mathbf{B}_\infty\) structures on \(P\).
	These two sets are in bijection via the maps:
	\begin{align*}
		\mathbf{B}_\infty(P)&\to\operatorname{Bialg}(P)\\
		\lAngle \argument, \argument\rAngle&\mapsto u*v\coloneqq\sum_{k\ge 1}\sum_{(u)^k,(v)^k}\lAngle u_{(1)},v_{(1)}\rAngle\otimes \dotsm\otimes\lAngle
		u_{(k)},v_{(k)}\rAngle\\
		\operatorname{Bialg}(P)&\to\mathbf{B}_\infty(P)\\
		*&\mapsto \lAngle u,v\rAngle\coloneqq\pi_1(u*v)
	\end{align*}
	\end{theorem}
We note that symmetry of \(\lAngle  \argument, \argument\rAngle\) corresponds to commutativity of the product; the
result holds more generally for non-symmetric brackets and non-commutative bialgebra products. Applied to our setting,
we have that the free commutative product on forests, written in \(\tens(\cP) \cong_\top \hck\), can be recovered from the \(\mathbf{B}_\infty\) Bracket operation associated to it. From this we deduce that under the isomorphism $\tee$ the product $\cdot$ is not free any more.

\begin{corollary}\label{cor:fprod}
  Let \(p_1,\dotsc,p_m,q_1,\dotsc,q_n\in\cP\). Then
  \begin{equation*}\label{eq:bInfProd}
    \begin{split}
      &(p_1 \tee \cdots \tee p_m)\cdot (q_1 \tee \cdots \tee q_n) \\
      ={} &\sum_{k \geq 1} \sum_{\substack{m^k = m \\ n^k = n}} \pi((p_1 \tee \ldots  p_{m_1}) \cdot (q_1 \tee \ldots  q_{n_1})) \tee \cdots \tee \pi((p_{m^{k-1}+1} \tee \ldots  p_m) \cdot (q_{n^{k-1}+1} \tee \ldots  q_n))\,,
\end{split}
  \end{equation*}
 where we set \(m^k= m_1 + \ldots + m_k\) and \(n^k= n_1 + \ldots + n_k\).
\end{corollary}

We regard the category of commutative \(\binf\)-algebras, \(\bnohat\), as the full subcategory of Hopf algebras, consisting of graded, connected, commutative and cofree Hopf algebras with specified cofreeness projection, i.e.\ isomorphism to \(\bigotimes (\cP)\) endowed with a compatible product. We do not require the cofreeness projection to be preserved: this would amount to only considering morphisms that preserve the primitiveness grading, whereas general Hopf maps only preserve the associated filtration, and would be too restrictive (with a notable exception following shortly). We end this subsection with the promised functoriality statement:
\begin{corollary}[\(\binf\)-Connes-Kreimer functor]\label{cor:natural}
  \((\hck, \cdot, \dck, \pi)\) defines a functor \(\cat{Vec} \to \bnohat\), naturally isomorphic through \(\top\) to \(\tens(\cP)\) (where \(\cP\) is viewed as a functor \(\cat{Vec} \to \cat{Vec}\) and \(\tens\) as one \(\cat{Vec} \to \cat{Coalg}\)) endowed with the product $\cdot$. Moreover, for a linear map \(\theta \colon V \to W\), the induced map preserves the cofreeness projection: \(\pi_W \circ \cH(\theta) = \cH(\theta) \circ \pi_V\).
  \begin{proof}
    We must show that induced maps $\hck(\theta)$ are bialgebra maps that commute with the projections $\pi$. We perform this check for the coproduct, the other two checks are analogous: this amounts to the statement that, for $\theta \colon V \to W$, $\cH(\theta)^{\otimes 2} \circ \Delta_{\mathrm{CK},V} = \Delta_{\mathrm{CK},W} \circ \cH(\theta)$. Define $\Phi$ as in \Cref{prop:functor} by $\phi_b = 1$ if $m = 1$, $n = 2$ and $b$ is precisely the bijection between the top part $\scrg_1$ and the bottom part $\scrg_2$ of an admissible cut of $\scrf_1$, and $\phi_b = 0$ otherwise. Then for every finite-dimensional vector space $U$, $\Phi_U = \Delta_{\mathrm{CK},U}$, the Connes-Kreimer coproduct on $\hck(U)$, and the required statement follows from the proposition. Note how, in particular, this implies that $\hck(\theta)$ restricts to a map $\cP(\theta) \colon \cP(V) \to \cP(W)$, which is required for functoriality of $\cP$.
  \end{proof}
\end{corollary}

\subsection{Duality with Grossman-Larson}\label{subsec:GL}
The \emph{Grossman-Larson Hopf algebra}, which we denote \(\hgl\), was defined independently of (and prior to) the Connes-Kreimer Hopf algebra, in \cite{GL89}.
While \(\hck\) is useful to encode iterated integration, \(\hgl\) was introduced to represent composition of differential operators.
The fact that these two Hopf algebras could be related by duality was noticed only a decade later \cite{Pan00,Hof03}.
The main purpose of this section is to revisit this dual pairing in the setting of functors to (co)free Hopf algebras.

While we follow \cite{Hof03} as our main reference on the Grossman-Larson Hopf algebra, we proceed with one technical difference. We take the underlying vector space of \(\hgl\) to be \(\cH = \mathrm{span}_\mathbb{R}(\scrH)\), not \(\mathrm{span}_\mathbb{R}(\scrT)\) (this was also done, for example, in \cite{HK15}). In this notation, the Grossman-Larson product of two forests is defined by 
\begin{equation}\label{eq:glProd}
\mathscr{f}\star\mathscr{t}_1\dotsm\mathscr{t}_n=(\mathscr{f}_{(1)}\curvearrowleft\mathscr{t}_1)\dotsm(\mathscr{f}_{(n)}\curvearrowleft\mathscr{t}_n)\mathscr{f}_{(n+1)}\, \quad \scrt_k \in \scrT, \ \scrf \in \scrH\,,
\end{equation}
where the term \(\scrf \curvearrowleft \scrg\) is defined by summing over all possible graftings of roots in \(\scrf\) to vertices in \(\scrg\), and the terms $\mathscr{f}_{(i)}$ are the Sweedler notation for cofree cocommutative coproduct over \(\cT\) (which partitions the multiset of trees that makes up \(\scrf\) in all possible ways), denoted  \(\dgl\); in particular \(\mathscr s \star \scrg = \mathscr s \scrg + \mathscr s \curvearrowleft \scrg\) for any $\mathscr s, \mathscr g\in  \scrT $. Note that unlike in the definition of \(\scrf \tee \scrg\), different roots in \(\scrf\) can be grafted to different vertices (also, there is no normalising factor in the definition of \(\star\)). The triple \(\hgl \coloneqq (\cH, \star, \Delta_\mathrm{GL})\) is shown to be a graded connected bialgebra, which we view again as a functor \(\cat{Vec} \to \cat{gcHopf}\).

We now wish to define a graded dual pairing \(\cH(V^*) \otimes \cH(V) \to \bbR\), or equivalently (under the identification \(\cH(V^*)^* = \cH(V)\), where the external duality is intended as graded duality) a graded map \(\cH(V) \to \cH(V)\) which is full-rank at each level. In the notation of \Cref{prop:functor}, set \(\Phi_b = 1\) when \(m = 1 = n\), \(\scrf_1 = \scrg_1\) and \(b \in \mathbb S_{\scrf_1}\) and \(\Phi_b = 0\) otherwise: this induces a natural map
\begin{equation}\label{eq:dual}
	\langle \argument, \argument \rangle \colon \hck(V^*) \otimes \hgl(V) \to \mathbb R
\end{equation}
given by a direct sum of pairings \(\langle\argument, \argument \rangle_{\scrf} \colon V^{*\boxtimes \scrf} \otimes V^{\boxtimes \scrf} \to \bbR\) which counts the order automorphisms of \(\scrf\). For example
\begin{align*}
	\big\langle \Forest{[\gamma[\alpha][\beta]]}, \Forest{[\zeta[\xi][\eta]]} \big\rangle_{\Forest{[[][]]}} &= \langle \alpha \boxdot \beta, \xi \boxdot \eta \rangle \langle \gamma, \zeta \rangle  \\ &= \langle \alpha, \xi\rangle\langle \beta,\eta \rangle\langle \gamma, \zeta\rangle + \langle \alpha, \eta\rangle\langle \beta, \xi\rangle\langle \gamma, \zeta\rangle
\end{align*}

Define \(\mathcal N\) to be the orthogonal complement to \(\cP\) defined by \(\top\), i.e.\ \(\mathcal N \coloneqq \bigoplus_{n = 2}^\infty \cP^{\top n}\), so that \(\hck = \bbR \oplus \cP \oplus \mathcal N\), and define \(\cQ\) to be the annihilator of \(\bbR \oplus \cN\) in \(\hgl\) w.r.t.\ the pairing \eqref{eq:dual}, i.e.\ \say{dual primitives}
\begin{equation}
	\cQ \coloneqq \{q \in \hgl \colon \forall h \in \bbR \oplus \cN \ \langle h, q \rangle = 0\}
\end{equation}
The isomorphism of \Cref{thm:T} dualises to one \(\top^* \colon (\hgl(V), \star) \to (\tens (\cQ(V)),\otimes )\) which identifies \(\hgl(V)\) with the free algebra over \(\cQ(V)\), naturally in \(V\). Putting everything together, we have the following result, original to \cite[Proposition 4.4]{Hof03} (see \cite[Theorem 4.8]{CW17} for the decorated case). We recall that a pairing of bialgebras induces an obvious duality between products and coproducts.

\begin{theorem}[Duality between Connes-Kreimer and Grossman-Larson]\label{thm:dual}
	The pairing \eqref{eq:dual} is a non-degenerate pairing of graded Hopf bialgebras, and the following diagram
	\begin{equation}
	\begin{tikzcd}[column sep = large]
	\hck(V^*) \otimes \hgl(V) \arrow[rd,"\pairing"]  \arrow[dd,"\top^{-1} \otimes \top^*"', "\cong"] \\ & \bbR \\
	\tens\cP(V^*) \otimes \tens\cQ(V) \arrow[ru,"\pairing_{\cP,\cQ}",swap]
	\end{tikzcd}
	\end{equation}
	in which \(\langle \argument, \argument \rangle_{\cP, \cQ}\) denotes the canonical extension to the tensor algebras
  of the restriction of \eqref{eq:dual} to \(\cP(V^*) \otimes \cQ(V)\), commutes.
\end{theorem}
Here non-degeneracy refers to non-degeneracy at each level of the grading by weight. Note that an application of \Cref{prop:functor} shows how this diagram is functorial in pairs of maps \(\psi \colon W \to V\), \(\varphi \colon V \to W\) (the maps going in opposite directions because of contravariance of the dualisation functor); from now on we omit \(V\) and \(V^*\) from the notation whenever possible, with the understanding that \(\cP\) and \(\cQ\) are to be evaluated on dual vector spaces, with functoriality intended in the above sense.

We now consider computations in coordinates. Let \(B\) be a basis of \(V\) and \(B^*\) be its dual basis, and we identify \(\scrf \in \scrH(B)\) with its dual basis element in \(\scrH(B^*)\).
For \(x \in \hgl\), \(y \in \hck\) and \(\scrf \in \scrH(B)\), \(\scrg \in \scrH(B^*)\) we denote \(x^{\scrg} \coloneqq \langle \scrg, x \rangle\) and \(y_\scrf \coloneqq \langle y, \scrf \rangle\).
Given an undecorated forest \(\scrg \in \scrH\) we denote \(\ell_B(\scrf)\) the set of labellings of \(\scrf\) with elements of \(B\), i.e.\ its elements are maps \(l \colon \{\scrf\} \to B\) and \(\scrf_l \in \scrH(B)\) is the forest so labelled.
In the unlabelled case, it follows directly from the definition of the pairing that 
\begin{equation}
	\langle \scrf, \scrg \rangle = \deltaup_{\scrf\scrg}\varsigma(\scrf),\qquad \scrf,\scrg \in \scrH.
\end{equation}
In the labelled case, if \(\scrf \in \scrH(B)\), \(\scrg \in \scrH(B^*)\), \(\langle \scrf, \scrf \rangle\) is the cardinality of the \emph{labelled} symmetry group of \(\scrf\), i.e.\ the group of order automorphisms of \(\scrf\) that preserve the labelling, which may not coincide with \(\varsigma(\scrf)\).

\begin{lemma}[Coordinates and pairing in the forest basis]\ \\[-2ex] \label{lem:labelPairing}
	\begin{enumerate}
		\item For \(\scrf \in \scrH(B)\), and \(\widetilde\scrf \in \scrH\) the same forest stripped of its labels, there are \(\frac{\varsigma(\widetilde\scrf)}{\langle \scrf, \scrf \rangle}\) labellings \(l \in \ell(\widetilde\scrf)\) s.t.\ \(\widetilde\scrf_l = \scrf\);
		\item For  \(x \in \hgl\),
		\[
		x = \sum_{\scrf \in \scrH(B)} \langle \scrf, \scrf \rangle^{-1} x^\scrf \scrf = \sum_{\scrh \in \scrH} \varsigma(\scrh)^{-1} \sum_{l \in \ell(\scrh)} x^{\scrh_l} \scrh_l\,;
		\]
		and an analogous statement holds for \(y \in \hck\);
		\item For \(x \in \hgl\) and \(y \in \hck\)
		\[
		\langle y, x \rangle = \sum_{\scrf \in \scrH(B)} \langle \scrf, \scrf \rangle^{-1} x^{\scrf} y_\scrf = \sum_{\scrh \in \scrH} \varsigma(\scrh)^{-1} \sum_{l \in \ell(\scrh)}  x^{\scrh_l} y_{\scrh_l} \,.
		\]
	\end{enumerate}
	\begin{proof}
		The set of labellings of \(\widetilde \scrf\) that result in \(\widetilde\scrf_l = \scrf\) carries a transitive \(\mathbb S_{\widetilde \scrf}\)-action whose stabiliser is the subgroup of label-preserving order automorphisms of \(\scrf\). The first assertion then follows from the interpretation of \(\langle \scrf, \scrf \rangle\) as the cardinality of this subgroup and the orbit-stabiliser theorem. For all \(\scrg \in \scrH(B^*)\) 
		\[
		\Big\langle \scrg, \sum_{\scrf \in \scrH(B)} \langle \scrf, \scrf \rangle^{-1} x^\scrf \scrf \Big\rangle =  \sum_{\scrf \in \scrH(B)} x^\scrf \deltaup_{\scrf \scrg} = \langle \scrg, x \rangle \,,
		\]
		implying the first identity in (2) by non-degeneracy of the pairing; the second identity follows by decomposing the sum over \(\scrf \in \scrH(B)\) into one over \(\scrh \in \scrH\) and one over \(l \in \ell(\scrh)\) and applying (1) to count the number of collisions by which to divide. (3) is now a straightforward application of (2). 
	\end{proof}
\end{lemma}

We now describe how to compute the isomorphism \(\top^* \colon (\hgl(\argument), \star) \to (\tens \cQ(\argument),\otimes )\) and its inverse.
This broadly follows the same principle identified by \cite{Pan00,Hof00}, according to which grafting is the dual operation to cutting.
Since \(\top\) involves a grafting operation on \(\hck\), its dual will involve a cutting operation on \(\hgl\).
\begin{definition}\label{def:snips}
  We define a \emph{snip} of \(\scrh \in \scrH\) to be a non-empty collection of edges which share their lower endpoint; removing such edges separates \(\scrh\) into two forests, call the one that contains the roots \(\scrg\) and the one consisting of everything above the snipped edges \(\scrf\): we use the notation \((\scrf, \scrg) \in \mathrm{Snip}(\scrh)\) to denote this correspondence (note that the forests \((\scrf, \scrg)\) themselves do not identify the snip, and the same term in the sum may appear more than once). We then define
  \[
    \widetilde\Delta_\top \colon \hgl \to \hgl \otimes \hgl,\qquad \scrh \mapsto \sum_{(\scrf, \scrg) \in
    \mathrm{Snip}(\scrh)} \frac{1}{|\scrg|} \scrf \otimes \scrg.\qedhere
  \]
\end{definition}
Note that this operation does not define a coproduct, as it is not coassociative (we have only used the tilde as a reminder that both \(\scrf\) and \(\scrg\) are never empty). We define its iterations by
\begin{equation}
	\widetilde\Delta_\top^{(n)} \coloneqq ( \widetilde\Delta_\top^{(n-1)} \otimes \mathrm{id} ) \circ \widetilde\Delta_\top \colon \hgl \to \hgl^{\otimes (n+1)}
\end{equation}

\begin{proposition}
The adjoint of the binary operation \(\top\) \Cref{def:T} restricted to \(\cH_+^{\otimes 2}\) is \(\widetilde\Delta_\top\).
\begin{proof}
	We follow the same ideas of \cite[Proposition 4.4]{Hof03} with the main difference that snips replace admissible cuts (in fact, the proof is simplified since snips are easier to handle than admissible cuts). In the notation of \Cref{prop:functor}, we consider two collections \(\Phi = \{\phi_b\}_b, \Psi = \{\psi_b\}_b\), both set to vanish unless \(m = 2\) (call the two forests \(\scrf\) and \(\scrg\)) and \(n = 1\) (call it \(\scrh\)). Set \(\phi_b = |\scrg|^{-1}\) if \(b\) carries \(\scrf\) and \(\scrg\) to two subforests of \(\scrh\), in a way that is an order isomorphism for both, and so that all the roots of the image of \(\scrf\) are directly connected to the same vertex in the image of \(\scrg\); set \(\phi_b = 0\) otherwise. Set \(\psi_b = |\scrg|^{-1}\) for each \(b\) obtained by first grafting \(\scrf\) onto \(\scrg\), with all roots grafted onto the same vertex, and mapping the resulting forest to \(\scrh\) order-isomorphically; set \(\psi_b = 0\) otherwise. Then (under the identification \(\cH(V^*)^* = \cH(V)\)) the induced maps \(\Phi_{V^*}, \Psi_{V^*} \colon \cH(V^*)^{\otimes 2} \otimes \cH(V) \to \bbR\) are given by
	\begin{align*}
	x \otimes y \otimes z &\mapsto \langle x \otimes y, \widetilde\Delta_\top z \rangle \\
	x \otimes y \otimes z &\mapsto \langle x \tee y, z \rangle
	\end{align*}
	To show  these two maps are equal for all \(V\) we prove that \(\phi_b = \psi_b\) for all \(\scrf,\scrg,\scrh\) and \(b\) as above. Assume \(\phi_b = |\scrg|^{-1}\): then, since all roots of \(b(\scrf)\) are grafted onto the same vertex \(v\) in \(b(\scrg)\), the order isomorphisms \(b|_\scrf\) and \(b|_\scrg\) glue to an order isomorphism going from the forest obtained by grafting \(\scrf\) onto \(b^{-1}(v)\). This implies \(\phi_b = |\scrg|^{-1}\). Conversely, assume the latter: \(b\) must restrict to an order-isomorphism on both \(\scrf\) and \(\scrg\), and calling \(u\) the vertex of \(\scrg\) onto which the roots of \(\scrf\) are grafted in the construction of \(b\), \(b(\scrf)\) will be grafted onto \(b(u)\). Therefore \(\phi_b = |\scrg|^{-1}\) and since \(\phi_b,\psi_b\) can only take two values the proof is complete.
\end{proof}
\end{proposition}
The following explicit description of \(\top^*\) immediately follows by dualising \Cref{thm:T}, \eqref{recursive_pi} and \eqref{eq:pim}.
\begin{corollary}[Description of \(\hgl \cong \bigotimes ( \cQ)\)]\label{prop:dual}
	The map
	\[
	\tens ( \cQ ) \to \hgl, \qquad q_1 \otimes \cdots \otimes q_n \mapsto q_1 \star \cdots \star q_n
	\]
	is an algebra isomorphism with inverse \(\top^*\), and, denoting \(\pi_n^* \colon \hgl \to \cQ^{\star n}\) the canonical projections, \(\pi^* \coloneqq \pi^*_1\), we have
	\begin{align*}
		\pi^* &= \mathrm{id} - \star \circ (\mathrm{id} \otimes \pi^*) \circ \widetilde{\Delta}_\top \\
		\pi_n^* &= \star^{(n-1)} \circ \pi^{*\otimes n} \circ \widetilde\Delta_\top^{(n-1)}
	\end{align*}
	where \(\star^{(n-1)}\) denotes the Grossman-Larson product of \(n\) arguments.
\end{corollary}
Note the slight abuse of notation regarding the projections: strictly speaking, \(\pi_n^*\) is not the dual of \(\pi_n \colon \hck \twoheadrightarrow \cP^{\top n}\), rather of \(\top^n \circ \pi_n \colon \hck \to \hck\).

\begin{example}[Dual change of coordinates up to level \(3\)]\label{ex:dual}
	The following is a system of generators of \(\cQ^3\) obtained from \(\pi^*(\scrf)\), \(\scrf \in \hgl^3\); to make it a basis it suffices pick a system of representatives for the symmetry and antisymmetry relations and take the labels to range in a basis of \(V\).
	\begin{align*}
	\big\{\Forest{[\c]}, \ \Forest{[\a]} \Forest{[\b]}, \ \Forest{[\a]} \Forest{[\b]} \Forest{[\c]}, \ \tfrac 12 \big(\Forest{[\a]} \Forest{[\c[\b]]} - \Forest{[\c]} \Forest{[\a[\b]]}\big) \ \colon \ \alpha,\beta,\gamma \in V \big\}
	\end{align*}
	Note how, while ladders are well-adapted to the primitiveness grading of \(\hck\) (\(\mathscr l_n \in \cP^{\top n}\)), products of nodes \(\mathscr r_n \in \cQ\), as can be seen from \Cref{prop:dual} and the fact that \(\widetilde\Delta_\top \mathscr r_n = 0\) since there are no edges. We proceed to write all forest of weight \(\leq 3\) (not already contained in \(\cQ\)) in the \(\cQ^{\star\bullet}\) grading:
	\begin{align*}
	\Forest{[\b[\a]]} &= \Forest{[\a]} \star \Forest{[\b]} - \Forest{[\a]} \Forest{[\b]} \\
	\Forest{[\c[\b[\a]]]} &= \frac 12 \big(\Forest{[\b]} \Forest{[\c[\a]]} - \Forest{[\c]} \Forest{[\b[\a]]} \big)+\frac 12 \Forest{[\a]} \Forest{[\b]} \Forest{[\c]}  - \Forest{[\a]} \Forest{[\b]} \star \Forest{[\c]} - \frac 12 \Forest{[\a]} \star \Forest{[\b]} \Forest{[\c]} + \Forest{[\a]} \star \Forest{[\b]} \star \Forest{[\c]} \\
	\Forest{[\c[\a][\b]]} &= \frac 12 \big( \Forest{[\c]} \Forest{[\b[\a]]} - \Forest{[\b]} \Forest{[\c[\a]]}\big) +\frac 12 \big( \Forest{[\c]} \Forest{[\a[\b]]} - \Forest{[\a]} \Forest{[\c[\b]]} \big) - \frac 12 \Forest{[\a]} \star \Forest{[\b]}\Forest{[\c]} - \frac 12 \Forest{[\b]} \star \Forest{[\a]}\Forest{[\c]}  + \Forest{[\a]}\Forest{[\b]} \star \Forest{[\c]} \\
	\Forest{[\a]}\Forest{[\c[\b]]} &= \frac 12 \big( \Forest{[\a]}\Forest{[\c[\b]]} - \Forest{[\c]}\Forest{[\a[\b]]} \big)- \frac 12 \Forest{[\a]} \Forest{[\b]} \Forest{[\c]}   + \frac 12 \Forest{[\b]} \star \Forest{[\a]} \Forest{[\c]} 
	\end{align*}
	Note that the first term in each expansion is the image through \(\pi^*\) of the left hand side.
\end{example}

    \subsection{Pre-Lie algebras and higher-order differential operators}\label{subsec:preLie}
We start by recalling the notion of (left) pre-Lie algebra.
\begin{definition}
  A pre-Lie algebra is a vector space \(L\) endowed with a bilinear operation \(\rhd\colon L\otimes L\to L\) such that
  the associator
  \[
    \operatorname{a}_\rhd(x,y,z)\coloneqq x\rhd(y\rhd z)-(x\rhd y)\rhd z
  \]
  is symmetric in \(x\) and \(y\).
\end{definition}

It follows immediately from the definition that the Pre-Lie commutator of a Pre-Lie product  \([-,-]_\rhd\colon L\otimes L\to L\) given by $[x,y]_\rhd\coloneqq  x\rhd y- y\rhd x $ defines a Lie bracket on \(L\). We are led thus to consider the universal enveloping algebra \(\mathcal U(L)\), defined as the quotient of \(\tens(L)\) by the Hopf
ideal generated by the relation \(x\otimes y-y\otimes x-[x,y]_\rhd\).
Describing the product of \(\mathcal U(L)\) is difficult for a general Lie algebra, but in our present setting, one can
characterize it as a modification of the symmetric algebra \(\bigodot(L)\).
We now briefly recall the construction, due to J.-M. Oudom and D. Guin \cite{OG08}.

The first step is to extend the pre-Lie product to \(L\otimes\bigodot(L)\to\bigodot(L)\).
This is done by simply imposing that $\rhd$ is a derivation, that is \(x\rhd 1\coloneqq 0\)
\begin{equation*}
  x\rhd y_1\odot\dotsm\odot y_n \coloneqq \sum_{j=1}^n y_1\odot\dotsm\odot y_{j-1}\odot(x\rhd y_j)\odot y_{j+1}\odot\dotsm\odot y_n\,.
\end{equation*}
The extension to \(\bigodot(L)\otimes\bigodot(L)\to\bigodot(L)\) is more involved.
Consider \(\bigodot(L)\) as a bialgebra with its cofree cocommutative coproduct (\(\Delta_\odot x = x\otimes 1+1\otimes x\) for
\(x\in L\)).
Define recursively, for \(x,y\in L\) and \(A, B, C\in\bigodot(L)\),
\begin{align*}
      1\rhd x &\coloneqq x\\
      (x\odot A)\rhd y &\coloneqq x\rhd(A\rhd y)-(x\rhd A)\rhd y\\
      A\rhd (B\odot C) &\coloneqq (A_{(1)}\rhd B)\odot(A_{(2)}\rhd C)\,.
\end{align*}
Given this extension, one introduces a non-commutative product on \(\bigodot(L)\) by setting
\begin{equation}\label{eq:GO}
  A\circledast B\coloneqq (A_{(1)}\rhd B)\odot A_{(2)}\,.
\end{equation}
Note that from the recursive definition it immediately follows that,
\begin{equation}\label{eq:og.poly}
  A\circledast(y_1\odot\dotsm\odot y_n)=(A_{(1)}\rhd y_1)\odot\dotsm\odot(A_{(n)}\rhd y_n)\odot A_{(n+1)}\,.
\end{equation}
Using these definitions we can indeed summarize the main result of \cite{OG08} into the following:
\begin{theorem}\label{thm:GO}
  The map \(\mathcal U(L) \to (\bigodot(L),\circledast,\Delta)\) induced by $\mathrm{id}_L$ is a Hopf isomorphism. Moreover, for any \(A,B,C\in\bigodot(L)\) one has 
  \begin{equation}\label{lem:brace}
    A\rhd(B\rhd C)=(A\circledast B)\rhd C\,.
  \end{equation}  
\end{theorem}
Iterating the pre-Lie operation, we define the maps \(r^{(n)}_\rhd\colon L^{\otimes n}\to L\) recursively by
  \begin{align*}
    r^{(1)}(x)&\coloneqq x\\
    r^{(n+1)}(x\otimes x_1\otimes\dotsm\otimes x_n)&\coloneqq x\rhd r^{(n)}_{\rhd}(x_1\otimes\dotsm\otimes x_n)\,.
  \end{align*}
More generally, given \(x_1,\dotsc,x_n\in L\) and \(I=\{i_1<\dotsb<i_k\}\subset[n]\) a generic partition of \([n]\coloneqq\{1,\dotsc,n\}\), we set 
\[
  r_\rhd(x_I)\coloneqq r^{(k)}_\rhd(x_{i_1}\otimes\dotsm\otimes x_{i_k})\,.
\]
Thanks to these operators we have the following formula.
\begin{proposition}\label{prop:starSym}
  Let \(x_1,\dotsc,x_n\in L\), then we have the  equality  in \(\mathcal{U}(L)\).
  \[
    x_1\circledast\dotsm\circledast x_n = \sum_{B\in\operatorname{P}(n)}\bigodot_{I\in B}r_\rhd(x_I)\,,
  \]
  where we denoted by $\operatorname{P}(n) $ the set of all partitions of $[n]$.
\end{proposition}
\begin{proof}
  The proof is by induction, the case \(n=2\) being straightforward.
  Note that, by definition
  \[
    x_1\circledast\dotsm\circledast x_n\circledast x_{n+1} = \left( (x_{1,(1)}\circledast\dotsm\circledast
    x_{n,(1)})\rhd x_{n+1} \right)\odot(x_{1,(2)}\circledast\dotsm\circledast x_{n,(2)}).
  \]
  By \eqref{lem:brace} the right-hand side equals
  \[
    r^{(n+1)}_\rhd(x_{1,(1)}\otimes\dotsm\otimes x_{n,(1)}\otimes x_{n+1})\odot(x_{1,(2)}\circledast\dotsm\circledast x_{n,(2)}).
  \]
  However, since the \(x_i\) are primitive for the cofree cocommutative coproduct, this equals
  \[
    \sum_{I\subset[n]}r_\rhd(x_I\otimes x_{n+1})\odot x^{\circledast}_{[n]\setminus I}
  \]
  where \(x^{\circledast}_J\coloneqq x_{j_1}\circledast\dotsm\circledast x_{j_k}\). By the induction hypothesis we obtain
  \begin{align*}
    x_{[n]}^{\circledast}\circledast x_{n+1} &= \sum_{I\subset[n]}\sum_{J\in\mathcal{P}([n]\setminus
    I)}r_\rhd(x_I\otimes x_{n+1})\odot r_\rhd(x_{J_1})\odot\dotsm\odot r_\rhd(x_{J_k})\\
                                             &= \sum_{B\in\mathcal{P}(n+1)}\bigodot_{I\in B}r_\rhd(x_I).\qedhere
  \end{align*}
\end{proof}
Concisely, this formula relates the two exponentials in \(\bigodot(L)\), namely
\[
  \exp_\circledast\left(\sum_{i=1}^k x_i\right) = \exp_\odot\left(\sum_{n=1}^\infty r_\rhd(x_{[n]})\right).
\]

The following two examples of pre-Lie algebras are relevant to this article.
Given two (undecorated) trees \(\mathscr{s},\mathscr{t}\in\mathscr{T}\) and a vertex \(v\in \{\mathscr{t}\}\), the tree
\(\mathscr{s}\curvearrowleft_v\mathscr{t}\) is obtained by joining the root of \(\mathscr{s}\) onto \(v\) by means of
a new edge.
We then extend this to an operator \(\curvearrowleft\colon\mathscr{T}\otimes\mathscr{T}\to\mathscr{T}\) by setting
\[
  \mathscr{s}\curvearrowleft\mathscr{t}\coloneqq\sum_{v\in \{\mathscr{t}\}}\mathscr{s}\curvearrowleft_v\mathscr{t}.
\]
It is not difficult to show that this turns \(\mathscr{T}\) into a pre-Lie algebra.
Moreover, we have the following result.
\begin{theorem}[\cite{CL01}]\label{thm:CL}
  The pair \((\mathscr{T},\curvearrowleft)\) is the free pre-Lie algebra on one generator.
\end{theorem}

In this particular case, one can obtain an explicit form of the Oudom--Guin extension of grafting.
\begin{proposition}[\cite{OG08}]
  Let \(\mathscr{t}_1,\dotsc,\mathscr{t}_n\in\mathscr{T}\) and \(\mathscr{f}\in\mathscr{H}\). Then
  \[
    \mathscr{t}_1\dotsm\mathscr{t}_n\curvearrowleft\mathscr{f}=\sum_{(v_1,\dotsc,v_n)\in
  \{\mathscr{f}\}^n}\mathscr{t}_1\curvearrowleft_{v_1}(\dotsb(\mathscr{t}_n\curvearrowleft_{v_n}\mathscr{f})).
  \]
\end{proposition}

By the Cartier--Milnor--Moore theorem, the Oudom--Guin product obtained in the previous section coincides with the Grossman-Larson product \eqref{eq:glProd}.

Now consider the space of smooth vector fields on a finite-dimensional vector space \(W\), denoted here by \(\mathfrak{X} = \mathfrak{X}(W) \coloneqq C^\infty(W,W)\).
Recall that for \(F \in C^\infty(U,V)\), its \(n^\text{th}\) derivative can be seen as a map \(\D^n F \in C^\infty (U,
\mathcal L(U^{\otimes n},V))\). Here and for the rest of the article we denote \(\cL(W,W')\) the space of linear maps
between two vector spaces \(W\) and \(W'\). We also write \(\mathrm{End}(W)\coloneqq\cL(W,W)\).
\begin{definition}
  Let \(F,G\in\mathfrak{X}\). Define \(F\rhd G\in\mathfrak{X}\) by
  \[
    (F\rhd G)(x)\coloneqq \D G(x)F(x).\qedhere
  \]
\end{definition}
\begin{proposition}
  The pair \((\mathfrak{X},\rhd)\) is a left pre-Lie algebra.
\end{proposition}
\begin{proof}
  Follows from an easy computation: from the chain rule we see that
  \begin{align*}
    (F\rhd(G\rhd H))(x) &= \D^2H(x)(F(x),G(x))+\D H(x)\D G(x)F(x)\\
    ((F\rhd G)\rhd H)(x) &= \D H(x) \D G(x)F(x).
  \end{align*}
  Subtracting both expressions yields
  \[
    \operatorname{a}_\rhd(F,G,H)(x)= \D^2H(x)(F(x),G(x))
  \]
  which is symmetric in \(F\) and \(G\).
\end{proof}

Recall that \(\mathfrak{X}\) is a Lie algebra with bracket
\[
  [F, G](x)=\D G(x)F(x)-\D F(x)G(x) = (F\rhd G)(x)-(G\rhd F)(x),
\]
that is, the Lie bracket is induced by the pre-Lie product. One obtains a representation \(\varrho\) of \(\mathfrak{X}(W)\) on \(C^\infty(W)\) by setting
\[
  \varrho(F)\varphi(x)\coloneqq \D \varphi(x)F(x).
\]
Indeed, note that if \(F,G\in\mathfrak{X}\)
\begin{align*}
\varrho(F)\varrho(G)\varphi(x)-\varrho(G)\varrho(F)\varphi(x) &= \begin{multlined}[t]\D\varphi(x)\D G(x)F(x)+ \D^2\varphi(x)(F(x),G(x))\\ -
\D\varphi(x)\D F(x)G(x)- \D^2\varphi(x)(F(x),G(x))\end{multlined} \\
                    &= \D \varphi(x) [F, G](x)\\
                    &= \varrho([F,G])\varphi(x).
\end{align*}

The universal property of the universal envelope ensures that \(\varrho\) lifts uniquely to a representation of
\(\mathcal{U}(\mathfrak{X})\) on \(\operatorname{End}(C^\infty(W))\).
We let \(\mathfrak D = \mathfrak D(W)\coloneqq \operatorname{im}(\varrho)\subset \operatorname{End}(C^\infty(W))\) be the algebra of differential operators on \(W\).
In particular, for any \(F_1,\dotsc,F_k\in\mathfrak{X}(W)\) we have that
\[
 \varrho(F_1)\circ\dotsm\circ\varrho(F_k)\varphi = \varrho(F_1\circledast\dotsm\circledast F_k)\varphi.
\]
The following result makes the computation of \(\varrho\) fully explicit.
\begin{lemma}\label{lem:symprod}
  Let \(F_1,\dotsc,F_n\in\mathfrak{X}(W)\) and \(\varphi\in C^\infty(W)\).
  Then
  \[
    \varrho(F_1\odot\dotsm\odot F_n)\varphi(x) = \D^n\varphi(x)(F_1(x),\dotsc,F_n(x)).
  \]
\end{lemma}
\begin{proof}
  The proof is by induction, the base case being true by definition.
  Let \(F_0,\dotsc,F_n\in\mathfrak{X}\).
  From \eqref{eq:og.poly} we see that
  \[
    F_0\circledast(F_1\odot\dotsm\odot F_n)=F_0\odot F_1\odot\dotsm\odot F_n + \sum_{j=1}^nF_1\odot\dotsm\odot(F_0\rhd
    F_j)\odot\dotsm\odot F_n.
  \]
  By the induction hypothesis and the isomorphism \(\mathcal{U}(\mathfrak{X})\cong\mathfrak{D}(W)\) the action of the left-hand side on \(\varphi\) 
  equals
  \begin{align*}
    \D^{n+1}\varphi(x)(F_0(x),F_1(x),\dotsc,F_n(x))+\sum_{j=1}^n\D^n\varphi(x)(F_1(x),\dotsc,\D F_j(x)F_0(x),\dotsc,F_n(x)).
  \end{align*}
  Likewise, the action of the second term on the right-hand side equals
  \[
    \sum_{j=1}\D^n\varphi(x)(F_1(x),\dotsc,(F_0\rhd F_j)(x),\dotsc,F_n(x)),
  \]
  and since by definition we have that \((F_0\rhd F_j)(x)=\D F_j(x)F_0(x)\) the result follows.
\end{proof}

With this proposition in place we may compute composition of vector fields algebraically. As an example, consider
\(\varphi\in C^\infty(W,V)\) and \(F,G,H\in\mathfrak{X}\). Then
\begin{align*}
  \varrho(F)\varphi &= \D\varphi F \\
  \varrho(F)\varrho(G)\varphi &= \varrho(F\circledast G)\varphi \\
  &= \varrho(F\rhd G + F\odot G)\varphi \\
  &= \D \varphi(\D G F) +  \D^2\varphi (F, G) \\
  \varrho(F)\varrho(G)\varrho(H)\varphi &= (F\circledast G\circledast H)\varphi \\
  &= \bigl(F\rhd (G\rhd H) + F\odot(G\rhd H)+G\odot(F\rhd H)+H\odot(F\rhd G) + F\odot G\odot H\bigr)\varphi \\
  &= \D\varphi \D^2H(F,G)+ \D\varphi \D H\D G F +  \D^2\varphi (F, \D HG) + \D^2 \varphi (G, \D HF)\\
    &\phantom{={}}+ \D^2\varphi (H ,\D GF) + \D^3 \varphi (F, G, H)\,.
\end{align*}

Combining \Cref{thm:CL} with the universal property of the universal enveloping algebra and \Cref{lem:symprod} we have the following connection between $\hgl$ and vector fields:
\begin{proposition}\label{prp:prelie.vf}
 For any map \(F\in C^\infty(W,\cL(V, W))\cong\cL(V,\mathfrak{X}(W))\) there is a unique pre-Lie morphism \(\widehat{F}\colon\mathscr{T}(V)\to\mathfrak{X}(W)\) and a Hopf morphism \(\mathbf{F}\colon\hgl(V)\to\mathfrak{D}(W)\) such that  \(\left.\mathbf{F}\right\rvert_{\mathscr{T}(V)}=\widehat{F}\) and  \(\widehat{F}(\Forest{[\alpha]})= F_\alpha\). In particular,  for any \(\mathscr{f}=\mathscr{s}_1\dotsm\mathscr{s}_n\in\mathscr{H}(V)\), \(\mathscr{t}\in\mathscr{T}(V)\) one has 
  \[
    \widehat{F}_{\mathscr{f}\curvearrowleft\mathscr{t}}(x)=\D^n\widehat{F}_{\mathscr{t}}(x)(\widehat{F}_{\mathscr{s}_1}(x),\dotsc,\widehat{F}_{\mathscr{s}_n}(x))\,.\]
\end{proposition}
	
\section{The It\^o formula}\label{sec:ito}
	\subsection{Branched rough paths and the change-of-variables formula}\label{sec:brps}
	In this section we relate the algebraic structure developed before with the notion of branched rough paths \cite{Gub10}. We will follow \cite{HK15}, with the caveat that the dual pairing introduced in the previous section is not used therein, which will make our formulae look a little different. The theory of integration against branched rough paths has been since reworked more than once; we additionally refer the reader to \cite{FZ18} for the extension to rough paths with jumps. For \(T \geq 0\) let \(\Delta_T \coloneqq \{(s,t) \in [0,T]^2 \colon s \leq t\}\).
  A \emph{control} on \([0,T]\) is a continuous function \(\omega \colon \Delta_T \to [0,+\infty)\) s.t.\ \(\omega(t,t)
  = 0\) for \(0 \leq t \leq T\) and is super-additive, i.e.\ \(\omega(s,u) + \omega(u,t) \leq \omega(s,t)\) for \(0 \leq
  s \leq u \leq t \leq T\); the main example of a control is \(\omega(s,t) \coloneqq t-s\), but allowing for more
  general ones make it possible to consider paths of bounded \(\rho\)-variation that may only be \(\rho^{-1}\)-Hölder up to reparametrisation.
  Throughout this article, \(\rho\) will denote a real number \(\in [1,+\infty)\).
  Recall that \(\cH^n(V) \coloneqq \bigoplus_{|\scrf| \leq n} V^{\otimes \scrf}\) (and similar) where \(|\argument|\) refers to forest weight, i.e.\ number of vertices.
  (Co)algebra operations on such a space are understood to be automatically truncated.
  To distinguish forests of different weight we will also adopt the notation \(\cH^{(n)}(V) \coloneqq \bigoplus_{|\scrf| = n} V^{\otimes \scrf}\).
	\begin{definition}[Branched rough path]\label{def:brp}
		A \(V\)-valued \(\rho\)-\emph{branched rough path} on \([0,T]\) controlled by \(\omega\) is a map
		\begin{equation}
			\bX \colon \Delta_T \to \hgl^\p(V), \quad (s,t) \mapsto \bX_{s,t}
		\end{equation}
		satisfying the following three properties:
		\begin{itemize}
			\item \textbf{Chen property}: \(\bX_{s,t} = \bX_{s,u} \star \bX_{u,t}\), i.e. in coordinates
			\[\langle\scrf,\bX_{s,t}\rangle  = \langle{\scrf_{(1)}},\bX_{s,u} \rangle\langle{\scrf_{(2)}}, \bX_{u,t}\rangle\]
			for \(\scrf \in \scrH^\p(V^*)\), and \(0 \leq s \leq u \leq t \leq T\);
			\item \textbf{Character property}: \(\Delta_\mathrm{GL} \bX_{s,t} = \bX_{s,t} \otimes \bX_{s,t}\), i.e. in coordinates
			\[\langle\scrf\scrg,\bX_{s,t}\rangle = \langle\scrf, \bX_{s,t}\rangle \langle\scrg, \bX_{s,t}\rangle\]
			for \(\scrf,\scrg \in \scrH(V^*)\) with \(|\scrf| + |\scrg| \leq \p\), and \(0 \leq s \leq t \leq T\);
      \item\textbf{Regularity}: for \(\scrf \in  \scrH^\p(V^*)\)

			\[\displaystyle\sup_{0 \leq s < t \leq T} \frac{|\langle\scrf, \bX_{s,t}\rangle|}{\omega(s,t)^{|\scrf|/\rho}} < \infty\,.\]

		\end{itemize}
		We denote the set of these \(\mathscr C^\rho_\omega([0,T],V)\).
	\end{definition}
  The intuitive meaning of a branched rough path is given by \(\langle\mathbf{1},\bX_{s,t}\rangle = 1\), \(\langle \scrt_1 \cdots \scrt_n,\bX_{s,t}\rangle = \langle\mathscr t_1\bX_{s,t} \rangle\cdots \langle\mathscr t_n,\bX_{s,t}\rangle\) for \(\scrt_k \in \scrT(V^*)\) and
	\begin{equation}\label{eq:brPostulate}
		\langle [\mathscr f]_\alpha,\bX_{s,t}\rangle \mathbin{\text{\say{\(=\)}}} \int_s^t\langle {\mathscr f} ,\bX_{s,u}
    \rangle\, \dif X_u^\alpha.
	\end{equation}

	The meaning of this last identity is only heuristic: the integral is not well-defined by Stieltjes/Young integration
	when \(\rho \geq 2\). When equipped with an initial value \(X_0\), the components of \(\bX\) indexed by single
  labelled vertices are the increments of components \(X^\alpha\) of a continuous function \(X \colon [0,T] \to V\)
  called the \emph{trace}; \(X\) belongs to \(\mathcal C^\rho_\omega([0,T],V)\), the set of functions \(Z \colon [0,T] \to V\) with the property that, assuming \(V\) from now on to be finite-dimensional
	\begin{equation}
		\sup_{0 \leq s < t \leq T} \frac{|Z_{s,t}|}{\omega(s,t)^{1/\rho}} < \infty
	\end{equation}
	where \(Z_{s,t} \coloneqq Z_t - Z_s\) is the increment and \(|\cdot|\) is a generic norm on \(V\).

	In what follows we will write \(\approx_m\) between two real-valued quantities dependent on \(0 \leq s \leq t \leq T\) to mean that their difference is of order \(O(\omega(s,t)^{m/\rho})\) on \([0,T]\) and simply \(\approx\) (\emph{almost} equal) to mean \(\approx_{\p + 1}\). Note that, since the vector spaces in which these quantities take values will always be finite-dimensional, the meaning of \(\approx_m\) is independent of the choice of a norm of \(W\). When considering expansions, quantities that are \(\approx 0\) will become negligible, which is why will will always be able to truncate sums at order \(\p\). We now give the definition of path controlled by a branched rough path, which will be used to define integration. We use sub/superscripts to elements of \(\cH\) to refer to (forest weight) grading and \(\mathcal L\) refers to the space of linear maps. The following definition is original to \cite{Gub04,Gub10}, see also \cite{HK15} for a formulation more similar to the one below.
	\begin{definition}\label{def:contr}
		Let \(\bX \in \scrC_\omega^\rho([0,T],V)\). A \(U\)-valued \(\bX\)-controlled path of regularity $m\leq \p+1$  is an element \(\bH\) belonging to  \(\mathcal C^\rho_\omega([0,T],\mathcal L(\hgl^{m- 1}(V),U))\) such that for \(n = 0,\ldots, m\)
		\begin{equation}\label{eq:brContr}
			\langle \bH_{t},h \rangle \approx_{m - n} \langle  \bH_{s}, \bX_{s,t} \star h \rangle, \quad h \in
      \hgl^{(n)}(V),
		\end{equation}
		where to apply the pairing we are identifying \(\mathcal L(\hgl^{m - 1}(V),U) = \hck^{m - 1}(V^*) \otimes U\). Call the set of these \(\mathscr D_{\bX}^{m}(U)\).
	\end{definition}

	The intuition behind \Cref{def:contr} is that the \emph{trace} \(H \coloneqq \langle\bH,1\rangle\in \mathcal
  C_{\omega}^\rho([0,T],U)\) is a path defined in terms of \(\bX\) (such as the solution to an \(\bX\)-driven
  differential equation); higher terms in \(\bH\), the so-called \emph{Gubinelli derivatives}, represent formal derivatives of \(H\) w.r.t. each component\ \(\bX\). We denote the portion of \(\bH\) that constitutes the Gubinelli derivatives \cite{Gub04}, i.e.\ the restriction of \(\bH\) to forests of positive degree, \(\bH_\bullet\).

  An essential property of the space of \(\bX\)-controlled paths is that it is stable under images through smooth
  functions. Given two finite dimensional vector spaces $W, Z$ we denote by \(\in C^m(W,Z)\) the space of functions of class $C^m$ between $W$
 and $Z$, see \cite[Lemma 8.4]{Gub10} and \cite[Proposition 2.8]{Bel20}.
 \begin{lemma}[Smooth functions of controlled paths are controlled paths]\label{lem:smoothContr}
		Let \(\bH \in \mathscr D_{\bX}^{m}(W)\) and \(\varphi \in C^{m}(W,Z)\). Then the definition
		\begin{equation}\label{eq:fContr}
			\varphi(\bH)_\bullet \coloneqq \sum_{n = 1}^{m - 1} \frac{1}{n!} \mathrm{D}^n \varphi(H) \circ \bH^{\otimes n}_\bullet \widetilde\Delta^{(n-1)}_\mathrm{GL}
		\end{equation}
		defines an element \(\varphi(\bH) \in \mathscr D_{\bX}^{m}(Z)\) with trace \(\varphi(H)\).
	\end{lemma}

	Recall that for any $\varphi\in C^m(W,Z) $ and $ n\leq m $ one has \(\D^n\varphi \in C^{m-n}(W,\mathcal{L}(W^{\otimes n}, Z))\); in coordinates this reads
	\[
	\langle
	\varphi(\bH),\scrf\rangle=\sum_{n=1}^{m-1}\frac{1}{n!}\D^n
  \varphi(H)(\bH_{\scrf^{(1)}}\otimes\dotsm\otimes\bH_{\scrf^{(n)}}),
	\]
	where recall that the Sweedler notation indicates summing over all ways of partitioning up the forest \(\scrf\) into \(n\) non-empty subforests.

	Integration against controlled paths is usually introduced by considering \(\bH \in \mathscr D_{\bX}^{\p}(\mathcal
  L(V,W))\) and defining \(\int_0^T\bH \dif \bX \in W\) as the limit
	\begin{equation}\label{eq:intLim}
		\int_0^T\bH\dif \bX =\lim_{|\pi_n|\to 0}\sum_{[s,t] \in \pi_n} \langle \bH_s ,  \bX_{s,t}\rangle
	\end{equation}
  where \(|\pi_n|\) is  the mesh-size of the partition \(\pi_n\) of the interval \([0,T]\).
	The intuition is the same as in the better-known case of \(\rho \in [2,3)\) \cite{FH20}: compensating the ordinary Riemann sums with the higher terms of \(\bX\) contracted against the Gubinelli derivatives of \(\bH\) result in almost additivity of the two parameter function \(\langle \bH_s ,  \bX_{s,t}\rangle\).
  Here we are interested in a generalisation of this operation of integration, which will be sufficiently expressive to write general controlled paths (with an extra degree of controlledness --- see \Cref{rem:itoContr} below) as integrals.

  We begin by observing that primitive elements have a special significance in the context of rough paths. Indeed, \(X^p_{s,t} \coloneqq \langle p,\bX_{s,t} \rangle\) for \(p \in \cP\) are increments of paths: for \(s \leq u \leq t\)
	\[
	X^p_{s,u} + X^p_{u,t} = \langle p \otimes 1 + 1 \otimes p, \bX_{s,u} \otimes \bX_{u,t} \rangle = \langle \Delta_\mathrm{CK} p, \bX_{s,u} \otimes \bX_{u,t} \rangle = \langle p, \bX_{s,t} \rangle= X^p_{s,t}
	\]
	(note the use of normal font for \(X\) to emphasise evaluation on primitives, i.e.\ this is viewed as an extension of the trace of \(\bX\)). It is therefore natural to wish for integration against \(X^{\alpha \in V}\) to extend to \(X^{p \in \cP}\). The algebraic counterpart to this operation is provided by \(\top\): the following heuristic identity
	\begin{equation}\label{eq:Tint}
		\langle h \tee p, \bX_{s,t} \rangle  \text{\say{\(=\)}}  \int_s^t \langle h,\bX_{s,u} \rangle\,\dif X^p_{u}
	\end{equation}
	is supported by the fact that it satisfies the correct rule for breaking up the integral on \([s,t]\) into one on
  \([s,u]\) and one on \([u,t]\) (the so-called Chasles relation):
	\begin{align*}
		\langle h \tee p, \bX_{s,t} \rangle  &\text{\say{\(=\)}} \int_s^t \langle h,\bX_{sr} \rangle\,\dif X^p_r = \int_s^u
    \langle h,\bX_{sr}\rangle\,\dif \bX^p_r + \int_u^t \langle h,\bX_{sr}\rangle\,\dif X^p_r \\
		&=\int_s^u \langle h,\bX_{sr}\rangle\,\dif X^p_r +  \langle h_{(1)},\bX_{s,u}\rangle  \int_u^t \langle
    h_{(2)},\bX_{ur}\rangle\,\dif \bX^p_r \\
		&\text{\say{\(=\)}} \langle h \tee p , \bX_{s,u} \rangle + \langle h_{(1)} \otimes (h_{(2)} \tee p), \bX_{s,u} \otimes \bX_{u,t} \rangle\,.
	\end{align*}
	Even though the first and last identities above are heuristic, the first and last items are actually equal by the
  Chen identity and the 1-cocycle property of \(\top\), that is, \eqref{eq:cocycle}.
  Note that \eqref{eq:Tint} extends \eqref{eq:brPostulate} thanks to the fact that \(h \tee\Forest{[\alpha]} = \EuScript{B}_+^\alpha(h)\) where \(\EuScript{B}_+\) is defined by
	\[
	\EuScript{B}_+ \colon \hck \otimes V^* \to \hck, \qquad \EuScript{B}_+^\alpha(\scrf) \coloneqq [\scrf]_\alpha
	\]
	still satisfies \eqref{eq:cocycle} but is defined on a restricted space. The presence of this operator is what made it possible, in the literature published so far, to define integration against \(X^{\alpha \in V}\). Below we will see examples of Hopf algebras that do not have such an operator, and that are thus not good candidates to describe an integration theory.

	In fact, something more can be said in the way of making \eqref{eq:Tint} rigorous. Recall that for \(\bX \in \mathscr
  C_\omega^\rho([0,T],V)\) we can extend it to arbitrarily high degree by defining its \emph{branched signature} as the convergent limit of Grossman-Larson products
	\begin{equation}\label{eq:signature}
		\EuScript S (\bX)_{s,t} \coloneqq \lim_{n \to \infty}\underset{[u,v] \in \Pi_n}{\bigstar} \bX_{u,v} \quad\in \hgl \doubleBr{V}
	\end{equation}
	independently of the partition \(\Pi_n\) with vanishing mesh size as \(n \to \infty\) \cite[Theorem 7.3]{Gub10}. The
  signature extends the rough path, in the sense that the projection of \(\EuScript{S}(\bX)\) onto \(\hgl^\p\) equals
  \(\bX\), but \(\EuScript{S}(\bX)\) belongs to the algebraic (ungraded) dual of \(\hck(V^*)\), defined by taking a
  direct product in \eqref{eq:FV} instead of the direct sum. From this it follows that \(\top^*
  (\EuScript{S}(\bX)_{s,t}) \in \tens \doubleBr{\cQ \doubleBr{V}}\) satisfies a similar expression to the above, with a
  tensor product over intervals in \(\Pi_n\) in place of the Grossman-Larson product, the same expression valid for the
  signature of a geometric rough path \cite[Theorem 2.2.1]{Lyo98}. In particular, if \(p \in \cP\) and \(h \in \hck \tee
  (\bigoplus_{n > \p - |p|} \cP^{(n)})\), i.e.\ \(h\) is a linear combination of elements of the form \(k \tee q\) with
  \(|p| + |q| \geq \p + 1\), the path \(u \mapsto \langle h,\EuScript{S}(\bX)_{s,u}\rangle\) is Young integrable against the path \(u \mapsto X^p_u \coloneqq \langle p, \EuScript{S}(X)_{s,u}\rangle\), and \eqref{eq:Tint} (with \(\langle h, \bX\rangle\) replaced with \(\langle h,\EuScript{S}(\bX)_{s,u}\rangle\)) holds literally, with the integral taken in the sense of Young.

	The next lemma extends this notion of integration to general \(\bX\)-controlled paths.

	\begin{proposition}[Rough integration of controlled paths against \(X^{p \in \cP}\)]\label{lem:int}
		Let \(\bH \in \mathscr D_{\bX}^{\p}(\mathcal L(\cQ, W))\).
    Then
		\[
		\langle \bH_s, \bX_{s,t} \rangle - \langle \bH_s, \bX_{s,u} \rangle - \langle \bH_u, \bX_{u,t} \rangle \approx 0
		\]
		where we are identifying
		\begin{equation}\label{eq:ident}
			\cL(\hgl , \cL(\cQ, W)) = \cL(\hgl \otimes \cQ, W) = (\hck \otimes \cP) \otimes W \xrightarrow{\top \otimes \mathrm{id}} \hck \otimes W
		\end{equation}
		to apply the pairing. This implies that the limit
		\[
		 \int_s^t \bH\,\dif \bX \coloneqq \lim_{|\pi_n| \to 0} \sum_{[u,v] \in \pi_n} \langle \bH_u, \bX_{uv} \rangle
		\]
		exists independently of the sequence of partitions \((\pi_n)_n\) on \([s,t]\) with vanishing mesh size as \(n \to
    \infty\). The resulting limit is additive. i.e. \(\int_s^u \bH\,\dif
    \bX + \int_u^t \bH \,\dif \bX= \int_s^t \bH\,\dif
    \bX \) for \(s \leq u \leq t\) and it satisfies \ \(\int_s^t \bH\,\dif \bX \approx \langle \bH_s,
    \bX_{s,t} \rangle\). Moreover, setting \(\langle \int \! \bH\,\dif \bX, h \star q \rangle \coloneqq \langle \bH_q, h \rangle \in W\) for \(h \in \hgl\) and \(q \in \cQ\) (where now \(\bH_q \in \cL(\hgl,W) = \hck \otimes W\)) defines an element of \(\mathscr D_{\bX}^{\p+1}(W)\).
	\end{proposition}
    \begin{remark}
      The contractions in the statement of \Cref{lem:int} contain summations over forests up to a
      certain degree depending on the regularity of \(\bX\). We chose not to include this in the notation at the expense
      of possibly including terms of high regularity.
      In practice this is not an issue since condition \eqref{eq:brContr} and its version in \Cref{lem:int} only hold up to
      a remainder, so when doing computations one can just discard these extra terms.
    \end{remark}

    \begin{proof}[Proof of \Cref{lem:int}]
			The proof is standard and easily adapted from the case of geometric rough paths. Denote \(\bX^{m_1\ldots m_k}\) the restriction of \(\bX\) to \(\cP^{m_1 \tee \ldots \tee m_k} \coloneqq \top(\cP^{(m_1)} \otimes \cdots  \otimes \cP^{(m_k)})\) (recall that \(\cP^{(m)}\) are Connes-Kreimer primitives of homogeneous weight equal to \(m\)), and similarly for \(\bH_{m_1\ldots m_k}\) the projection of \(\bH\) onto \(\cP^{m_1 \tee \ldots \tee m_k} \otimes W\). Then, implying below a sum on \(k\) and on \(m_i \geq 1\) subject to \({m_1 + \ldots + m_k \leq \p}\)
			\begin{align*}
				&\langle \bH_s, \bX_{s,t} \rangle - \langle \bH_s, \bX_{s,u} \rangle - \langle \bH_u, \bX_{u,t} \rangle \\
				={}& \langle \bH_{m_1\ldots m_k;s} , \bX^{m_1 \ldots m_k}_{s,t} \rangle - \langle \bH_{m_1\ldots m_k;s}, \bX^{m_1 \ldots m_k}_{s,u}\rangle - \langle \bH_{m_1\ldots m_k;u}, \bX^{m_1 \ldots m_k}_{u,t} \rangle \\
				={}& \langle \bH_{m_1\ldots m_k;s} , \bX^{m_1 \ldots m_k}_{s,t} - \bX^{m_1 \ldots m_k}_{s,u} - \bX^{m_1 \ldots m_k}_{u,t}\rangle - \langle\bH_{m_1\ldots m_k;s,u} ,\bX^{m_1 \ldots m_k}_{u,t}\rangle \\
				={}&\sum_{i = 1}^{k-1}\langle \bH_{m_1\ldots m_k;s},  \bX^{m_1\ldots m_i}_{s,u} \otimes \bX^{m_{i+1}\ldots m_k}_{u,t}\rangle - \langle\bH_{m_1\ldots m_k;s,u} ,\bX^{m_1 \ldots m_k}_{u,t}\rangle \\
				\approx{}& \sum_{i = 1}^{k-1}\langle \bH_{m_1\ldots m_k;s},  \bX^{m_1\ldots m_i}_{s,u} \otimes \bX^{m_{i+1}\ldots m_k}_{u,t}\rangle - \langle \bH_s, \bX_{s,u} \otimes \bX_{u,t}^{m_1 \ldots m_k} \rangle \\
				={}& 0
			\end{align*}
			where we have used \eqref{eq:brContr} combined with the fact that \(y = \bX_{u,t}^{m_1 \ldots m_k}\) contributes an additional \(O(\omega(s,t)^{m_1+\ldots+m_k})\). The classical statement on almost-additivity \cite[Theorem 3.3.1]{Lyo98} implies existence of the limit, uniqueness of the additive functional, and controlledness at level \(0\). Since the statement extends linearly, we may assume \(y\) to be of the form \(y = h \star q\), and we have
			\[
			\textstyle	\langle \int_0^t \bH\,\dif \bX,y \rangle = \langle\bH_{q;t}, h \rangle \approx_{\p - |h|} \langle
      \bH_{q;s}, \bX_{s,t} \star h \rangle = \langle \int_0^s \bH\,\dif \bX,\bX_{s,t} \star y \rangle
			\]
			and notice that integration has introduced one extra (inhomogeneous) degree of controlledness (corresponding to \(|q|\) above); stated otherwise, some terms in this controlled path could be discarded, but this is of little importance. This concludes the proof.
		\end{proof}

	We now consider computations in coordinates. When considering, say, the integral against a \(V\)-valued bounded
  variation path \(\int H\,\dif X\) with \(H\) is a \(\cL(V,W)\)-valued path, given a basis of \(V\) it is often helpful
  to decompose the integral as a sum \(\int H_\gamma\,\dif X^\gamma\) over the chosen coordinates. This carries over straightforwardly to the case of rough integrals in the case in which \(\bH \in \mathscr D_\bX(\mathcal L(V,W))\), but if \(\bH \in \mathscr D_\bX(\mathcal L(\cQ,W))\) a basis of \(\cQ\) would be needed. As stated in \Cref{rem:bases} such bases are difficult to construct, and we will not do so here; the question then becomes whether it is possible to express rough integrals as linear combinations of ones against one-dimensional components of the rough path (in which the integrand however is still controlled by the whole of \(\bX\)), using only a basis of \(V\) and the projection \(\pi \colon \hck \twoheadrightarrow \cP\) to non-injectively obtain the components over which to sum. To do this, having picked a basis \(B\) of \(V\) (with dual basis \(B^*\) of \(V^*\)), we introduce the family of linear maps indexed by \(\scrf \in \scrH(B)\)
	\begin{equation}\label{eq:bigPi}
		\Pi_\scrf \colon \cH_{\mathrm{CK},+} \to \cH_{\mathrm{CK},+}, \qquad h \tee p \mapsto \frac{\langle p, \scrf \rangle}{\langle \pi(\scrf), \scrf \rangle} h \tee \pi(\scrf)
	\end{equation}
	where the argument is uniquely decomposed in the form considered with \(p \in \cP\) and \(h \in \hck\). Note that if \(p\) is in the linear span of \(\pi(\scrf)\), say for \(p = \lambda\pi(\scrf)\)
	\[
	\Pi_\scrf(h \tee p) = \frac{\langle \lambda \pi(\scrf), \scrf \rangle}{\langle \pi(\scrf), \scrf \rangle} h \tee \pi(\scrf) = h \tee p
	\]
	implying that \(\Pi_\scrf\) is a projection onto \(\{\scrg \tee \pi(\scrf) \mid \scrg \in \scrH(B^*) \}\). In the next lemma we use the same notation as in \Cref{lem:labelPairing}.
	\begin{lemma}[Rough integral in coordinates]\label{lem:intCoords}\ \\[-2ex]
		\begin{enumerate}
			\item Any element \(y \in \cH_{\mathrm{CK},+}\) can be written as
			\[
			y = \sum_{\scrf \in \scrH(B)} \frac{\langle \pi(\scrf), \scrf \rangle}{\langle \scrf,\scrf \rangle} \Pi_\scrf(y) \,;
			\]
			\item Let \(\bH \in \mathscr D_\bX^{\p}(\cL(\cQ,W))\). For any \(\scrf \in \scrH(B)\), identifying \(\bH\) with an element of \(\hck \otimes W\) thanks to \eqref{eq:ident} one has  \(\Pi_\scrf(\bH) \in \mathscr D_\bX^{\p}(\cL(\cQ,W))\) and we obtain

\[\langle \Pi_\scrf\bH_s, \bX_{s,t} \rangle= \sum_{\scrg \in \scrH} \varsigma(\scrg)^{-1} \sum_{l \in \ell(\scrg)} \bH_{\scrg_l \star \pi^*(\scrf);s} \bX^{\scrg_l \tee \pi(\scrf)}_{s,t}\]
			where again we use \eqref{eq:ident} to identify \(\bH\) with an element of \(\hck \otimes W\) which we evaluate on \(\scrg_l \star \scrf\). We will use the notation
			\[
			\int_s^t \bH_{\pi^*(\scrf)}\,\dif \bX^{\pi(\scrf)} \coloneqq \int_s^t
      \Pi_\scrf(\bH)\,\dif \bX\,.
			\]

			\item We have the identities
			\[
			\int_s^t\bH\,\dif \bX = \sum_{\scrf \in \scrH} \varsigma(\scrf)^{-1} \sum_{l \in \ell(\scrf)} \int_s^t
      \bH_{\pi^*(\scrf_l)}\,\dif \bX^{\pi(\scrf_l)} \]
      \[\langle \bH_s, \bX_{s,t} \rangle =\sum_{\scrf,\scrg \in \scrH} \varsigma(\scrf)^{-1}\varsigma(\scrg)^{-1} \sum_{\substack{l \in \ell(\scrf) \\ m \in \ell(\scrg)}} \bH_{\scrg_m \star \pi^*(\scrf_l);s} \bX^{\scrg_m \tee \pi(\scrf_l)}_{s,t}  \,.
			\]
		\end{enumerate}
			\begin{proof}
			For \(h \in \hck\) and \(p \in \cP\) we have
			\begin{align*}
				p &= \sum_{\scrf \in \scrH(B)} \frac{\langle p, \scrf \rangle}{\langle \scrf,\scrf \rangle} \pi(\scrf) \\
        \intertext{so that}
        h \tee p &= \sum_{\scrf \in \scrH(B)} \frac{\langle p, \scrf \rangle}{\langle \scrf,\scrf \rangle} h \tee \pi(\scrf) = \sum_{\scrf \in \scrH(B)} \frac{\langle \pi(\scrf), \scrf \rangle}{\langle \scrf,\scrf \rangle} \Pi_\scrf(h \tee p),
			\end{align*}
			and (1) follows from \eqref{lem:labelPairing}. Now, for \(k \in \hgl\) and \(q \in \cQ\)
			\begin{align*}
        \langle \Pi_\scrf(h \tee p), k \star q \rangle &= \frac{\langle p, \scrf \rangle}{\langle \pi(\scrf), \scrf
      \rangle} \langle h , k \rangle \langle \pi(\scrf), q \rangle = \langle h \tee p, \frac{\langle \scrf, q
      \rangle}{\langle f , \pi^*(\scrf) \rangle} k \star \pi^*(\scrf) \rangle,\\
      \intertext{whence}
      \Pi_\scrf(k \star q) &= \frac{\langle \scrf, q \rangle}{\langle f , \pi^*(\scrf) \rangle} k \star \pi^*(\scrf).
			\end{align*}
			This allows us to check the controlledness condition (making the necessary adjustments to take into account that we are operating under the identification of \eqref{eq:ident}):
			\begin{align*}
				\langle \Pi_\scrf(\bH_t), k \star q \rangle &= \langle \bH_t, \Pi_\scrf^*(k \star q) \rangle \\
				&= \frac{\langle \scrf, q \rangle}{\langle f , \pi^*(\scrf) \rangle}  \langle \bH_t, k \star \pi^*(\scrf) \rangle \\
				&\approx_{\p - |k|} \frac{\langle \scrf, q \rangle}{\langle f , \pi^*(\scrf) \rangle} \langle \bH_s,  \bX_{s,t} \star k \star \pi^*(\scrf) \rangle \\
				&= \langle \bH_t, \Pi^*_\scrf(\bX_{s,t} \star k \star q) \rangle \\
				&= \langle \Pi_\scrf(\bH_t),\bX_{s,t} \star k \star q \rangle
			\end{align*}
			Using \Cref{lem:labelPairing}, we compute
			\begin{align*}
				\langle \Pi_\scrf(h \tee p), k \star q \rangle &= \frac{\langle p , \scrf \rangle}{\langle \pi(\scrf) , \scrf \rangle} \langle h , k \rangle \langle \pi(\scrf) , q \rangle \\
				&= \frac{\langle p , \pi^*(\scrf) \rangle}{\langle \pi(\scrf) , \scrf \rangle} \sum_{\scrg \in \scrH} \varsigma(\scrg)^{-1} \sum_{l \in \ell(\scrg)} \langle h , \scrg_l \rangle\langle \scrg_l , k \rangle\langle \pi(\scrf) , q \rangle \\
				&= \langle \pi(\scrf) , \scrf \rangle^{-1} \sum_{\scrg \in \scrH} \varsigma(\scrg)^{-1} \sum_{l \in \ell(\scrg)} \langle h \tee p , \scrg_l \star \pi^*(\scrf) \rangle \langle \scrg_l \tee \pi(\scrf), k \star q \rangle
			\end{align*}
			and (2) now follows \(\int_s^t \Pi_\scrf(\bH)\,\dif \bX \approx \langle \bH_s , \bX_{s,t} \rangle\) and bilinearity of the pairing.
      Assertion (3) is a straightforward consequence of (1) and (2).
		\end{proof}
		\end{lemma}
		\begin{remark}
		Note that the coefficients \(\langle \pi(\scrf), \scrf \rangle\) (which depend in a complicated way on the labelling) disappear in the final version of the integral in coordinates. The definition
    \[
      \langle \pi(\scrf), \scrf \rangle^{-1}\int_0^\cdot\bH_{\pi^*(\scrf)}\,\dif \bX^{\pi(\scrf)} \coloneqq \int_0^\cdot\Pi_\scrf(\bH)\,\dif \bX
    \]
    is what makes this happen, and is justified by the presence of this coefficient as a factor to the Davie expansion that follows, while the use of integral notation is justified that this is an integral of the controlled path \(\Pi_\scrf(\bH)\). If we were provided with a basis \(\mathscr P(B)\) of \(\cP\) (say, extracted from the generating set \(\{\pi(\scrf) \colon \scrf \in \scrH(B) \}\)) with dual basis \(\mathscr Q(B)\) (in the sense that \(\langle p_i, q_j \rangle = 0\) if \(i \neq j\)), in which case we could rewrite (3) as
		\begin{equation}\label{eq:intBasis}
			\begin{split}
				\int_s^t\bH\,\dif \bX &= \sum_{\langle p, q \rangle \neq 0} \langle p, q \rangle^{-1} \int_s^t \bH_{p}\,\dif \bX^{q}
      \approx \sum_{\substack{n = 1,\ldots,\p \\ \langle q_k, p_k \rangle \neq 0}} \left(\prod_{k = 1}^n \langle p_k,
        q_k \rangle\right)^{-1} \bH_{q_1 \star \cdots \star q_n;s} \bX^{p_1 \tee \ldots \tee p_n}_{s,t}
			\end{split}
		\end{equation}
		We prefer, however, to provide formulae in the canonical forest coordinates, since this does not mean stating results in terms of a basis which is complicated to extract. Finally, we mention there is a third way of writing the expansion, which combines the \((\star,\tee)\)-powers of \eqref{eq:intBasis} and the basis-free formulation of \Cref{lem:intCoords}:
		\[
		\int_s^t\bH\,\dif \bX \approx \sum_{\scrf_1,\ldots, \scrf_n} \left( \prod_{k = 1}^n \varsigma(\scrf_n) \right)^{-1} \bH_{\pi^*(\scrf_1) \star \cdots \star \pi^*(\scrf_n);s} \bX^{\pi(\scrf_1) \tee \cdots \tee \pi(\scrf_n)}_{s,t}
		\]
	\end{remark}
	Using the extended notion of integration we expand the original definition of branched rough differential equation, to the case where the underlying vector field \(F\) is also allowed to take values in \(\mathcal L(\cQ(V),W)\), thanks to the natural inclusion \(V \hookrightarrow \cQ(V)\).

	\begin{definition}[RDEs driven by \(X^{p \in \cP}\)]\label{def:RDE}
		Let \(\bX \in \scrC_\omega^\rho([0,T],V)\), \(F \in C^\infty(W, \mathcal L(\cQ(V),W))\) and \(y_0 \in W\). We say that \(Y \colon [0,T] \to W\) solves the \emph{rough differential equation} (RDE)
		\begin{equation}\label{eq:rde}
			\dif Y_t = F(Y_t)\,\dif \bX_t, \quad Y_0 = y_0
		\end{equation}
		if there exists \(\bY \in \mathscr D^{\p}_\bX(W)\) with trace \(Y\) such that
		\[
		\bY_t = y_0 + \int_0^t F(\bY)\,\dif \bX,
		\]
		where \(F(\bY)\) is defined according to \eqref{eq:fContr} and the integral is defined thanks to \Cref{lem:int}.
	\end{definition}
\begin{remark}\label{rk:RDE} Since the focus of this work is on the algebraic relations underlying the equation \eqref{eq:rde}, we have chosen not to emphasise issues of existence and uniqueness. A unique solution exists as long as the vector fields are $C^{\p + 1}_\mathrm{b}$ \cite[Theorem 8.8]{Gub10}, cf.\ \Cref{rem:kelly} below.
\end{remark}

Unlike the case where the equation is driven by the paths \(X^\alpha\), the Davie expansion of the solution to
  \eqref{eq:rde} requires a bit more care since the map \(h\mapsto\langle \bY,h\rangle\) will not be induced by a pre-Lie
  morphism (see \Cref{prp:prelie.vf}).
  Nonetheless, we can make use of the pre-Lie structure of vector fields to obtain a recursive formula for the coefficients. In order to keep a  compact notation, we set
  \[
    \widetilde{\Delta}_\top\scrf=\mathscr{f}^1\otimes\mathscr{f}^0,
  \]
 and we refer to \Cref{def:snips} for the definition of the map \(\widetilde{\Delta}_\top\).
  The reason for the 0 index in \(\widetilde{\Delta}_\top\) is to reflect the fact that this operator is not coassociative and so can only be iterated on the left.

  \begin{definition}\label{def:Fhat}
    Let \(F\in\mathcal{L}(\mathcal{Q}(V),\mathfrak{X}(W))\). We recursively define a map
    \(\widehat{F}\colon\hgl(V)\to\mathfrak{X}(W)\) by
    \begin{equation}\label{eq:Fhat}
      \widehat{F}_h\coloneqq
      F_{\pi^*(h)}+\sum_{n\ge1}\frac{1}{n!}\widehat{F}_{(h^1)^{(1)}}\odot\dotsm\odot\widehat{F}_{(h^1)^{(n)}}\rhd
      F_{\pi^*(h^0)}.
    \end{equation}
    (with the symmetric product taking precedence on the pre-Lie product) and introduce the map \(\mathbf{F}\colon\hgl(V)\to\mathfrak{D}(W)\) by \(\mathbf{F}_1=\mathrm{id}\) and
  \begin{equation}\label{eq:Fbold}
    \mathbf{F}_h\coloneqq\sum_{n\ge1}\frac{1}{n!}\widehat{F}_{h^{(1)}}\odot\dotsm\odot\widehat{F}_{h^{(n)}}\,,
  \end{equation}
  for \(h\in(\hgl)_+(V)\), where the $n$-fold coproducts are taken according to $\Delta_\mathrm{GL}$.
  \end{definition}
  We note that in the case \(F\in\mathcal{L}(V,\mathfrak{X}(W))\), \(\widehat{F}\) coincides with the map $\mathbf{F}$ in \Cref{prp:prelie.vf} and it is also a pre-Lie morphism. This fact is well-known in the literature, and has been extended to the case of vector fields on a manifold in \cite{kern2023flow}.
  However, in the general setting the map \(\widehat{F}\) is \emph{not} a pre-Lie morphism, as e.g.
  \[
      \widehat{F}_{\Forest{[\alpha[\beta]]}}(y)=\D
      F_{\Forest{[\alpha]}}(y)F_{\Forest{[\beta]}}(y)-F_{\Forest{[\alpha]}\Forest{[\beta]}}(y)\neq(\widehat{F}_{\Forest{[\alpha]}}\rhd\widehat{F}_{\Forest{[\beta]}})(y).
    \]
    Nonetheless, we still have the recursion of \Cref{def:Fhat}, \(\mathbf{F}_{\mathbf{1}}=\widehat{F}_{\mathbf{1}}=\mathrm{id}\) and
    \begin{align*}
    \mathbf{F}_h\varphi(y) &= \sum_{n \geq 1} \frac{1}{n!} \D^n \varphi(y)(\widehat{F}_{h^{(1)}}(y), \ldots, \widehat{F}_{h^{(n)}}(y)) \\
    \widehat{F}_{h \star q}(y) &= \mathbf{F}_h F_q(y) = \sum_{n \geq 1} \frac{1}{n!} \D^n F_q(y)(\widehat{F}_{h^{(1)}}(y), \ldots, \widehat{F}_{h^{(n)}}(y)), \qquad \widehat{F}|_{\cQ} = F|_{\cQ}
    \end{align*}
    for $q \in \cQ$, $h \in(\hgl)_+$, $\varphi \in C^\infty(W)$.

  \begin{proposition}
    The map \(\mathbf{F}\) defined in \eqref{eq:Fbold} is a Hopf morphism.
  \end{proposition}
  \begin{proof}
    We begin by noting that by \Cref{thm:GO}, the coalgebra \((\mathfrak{D}(W),\Delta)\) is cofree cocommutative over
    \(\mathfrak{X}(W)\). In particular, given any other cocommutative coalgebra \((C,\Delta_C)\) every map \(\phi\colon
    C\to\mathfrak{X}(W)\) extends uniquely to a coalgebra morphism \(\Phi\colon C\to\mathfrak{D}(W)\) such that
    \(\pi_{\mathbf{X}(W)}\Phi=\phi\).
    Moreover,
    \[
      \Phi(c) = \sum_{n\ge 1}\frac{1}{n!}\phi(c_{(1)})\odot\dotsm\odot\phi(c_{(n)}).
    \]
    Applying this to the map \(\widehat{F}\colon\hgl\to\mathfrak{X}(W)\) we immediately see that \(\mathbf{F}\) is a
    coalgebra morphism.
    We note that in particular \(\pi_{\mathfrak{X}}\mathbf{F}=\widehat{F}\).
    A cofreeness argument, together with the fact that in \(\mathfrak{D}(W)\) we have \(\pi_{\mathfrak{X}(W)}(G\circledast H)=G\rhd\pi_{\mathfrak{X}(W)}H\) for
    any \(F,G\in\mathfrak{D}(W)\), yields that in order for \(\mathbf{F}\) to be an algebra morphism it is enough that
    \begin{equation}\label{eq:Fhatbar}
      \widehat{F}_{h\star h'}=\mathbf{F}_h\rhd\widehat{F}_{h'}.
    \end{equation}

    The proof proceeds in two steps.
    First, we show that \eqref{eq:Fhatbar} holds for \(q\in\mathcal{Q}\) and \(h\in\hgl^+\).
    Indeed, noting that for any \(\scrf\in\hck^+\) and \(p\in\mathcal{P}\) the formula
    \[
      \widetilde{\Delta}_{\mathrm{CK}}(\scrf\tee p)=\scrf'\otimes\scrf''\tee p+f\otimes p
    \]
    holds, we see immediately that
    \begin{align*}
      \langle \widetilde{\Delta}_\top(h\star q),\mathscr{f}\otimes p\rangle &= \langle h\star q,\scrf\tee p\rangle \\
      &= \langle h\otimes q,\scrf'\otimes\scrf''\tee p+\scrf\otimes p\rangle \\
      &= \langle h,\scrf\rangle\langle q,p\rangle,
    \end{align*}
    that is, \(\widetilde{\Delta}_\top(h\star q)=h\otimes q\). It follows immediately that
    \begin{align*}
      \widehat{F}_{h\star q}=\mathbf{F}_h\rhd\widehat{F}_q.
    \end{align*}
    Second, any \(h\in\hgl^+\) may be written as \(h=\hat h\star q\) for some \(\hat h\in\hgl^+\) with \(|\hat h|< |h|\) and
    \(q\in\mathcal{Q}\). Taking \(h,h'\in\hgl^+\) and letting \(h'=\hat h'\star q\) inductively we see that
    \begin{align*}
      \widehat{F}_{h\star h'} &= \mathbf{F}_{h\star\hat h'}\rhd\widehat{F}_q\\
                              &= (\mathbf{F}_h\circledast\mathbf{F}_{\hat h'})\rhd\widehat{F}_q \\
                              &= \mathbf{F}_h\rhd(\mathbf{F}_{\hat h'}\rhd\widehat{F}_q) \\
                              &= \mathbf{F}_h\rhd\widehat{F}_{h'}\qedhere.
    \end{align*}
  \end{proof}
  The following theorem extends Davie expansion to RDEs driven by the collection \(X^{p\in\cP}\).
  \begin{theorem}[Davie expansion with drifts]\label{thm:davie.drift}
    Let \(F\in\mathcal{L}(\mathcal{Q},\mathfrak{X}(W))\). A path \(Y\colon[0,T]\to W\) solves the RDE \eqref{eq:rde} if and only if $Y_0 = y_0$ and
    \begin{equation}\label{eq:davie.drift}
      Y_{s,t}\approx\langle \widehat F(Y_s),\bX_{s,t}\rangle \approx\sum_{\mathscr{f}\in\mathscr{H}}\varsigma(\mathscr{f})^{-1}\widehat{F}_{\mathscr{f}_l}(Y_s)\bX^{\mathscr{f}_l}_{s,t}
    \end{equation}
    and the map defined by \(\langle \bY_t,h\rangle\coloneqq\widehat{F}_h(Y_t)\) belongs to
    \(\mathscr{D}^{\p+1}_\bX(W)\).
  \end{theorem}
  \begin{proof}
    Let us first suppose that \(Y\) solves \eqref{eq:rde}. By evaluating both sides on the empty forest we immediately
    obtain from \Cref{lem:int} that
    \[
      Y_{s,t}\approx\langle F(\bY_s),\bX_{s,t}\rangle=\sum_{\scrf\in\mathscr{F}_+}\varsigma(\scrf)^{-1}\langle
      F(\bY_s),\scrf\rangle\bX^\scrf_{s,t}.
    \]
    We now show the expression for \(\langle \bY_t,\scrf\rangle\) by induction. In the case of \(q\in\mathcal{Q}\) we
    immediately see that
    \[
      \langle \bY_t,q\rangle=\langle F_q(\bY_t),\mathbf{1}\rangle = F_q(Y_t).
    \]
    If \(h\in\ker\pi^*\), assume the identity is true for all \(h'\in\hgl\) with \(|h'|< |h|\).
    Noting that \(h=\pi^*(h)+h^1\star\pi^*(h^0)\) we obtain
    \begin{align*}
      \langle \bY_t,h\rangle &= F_{\pi^*(h)}(Y_t) + \langle \bY_t,h^1\star\pi^*(h^0)\rangle\\
                           &= F_{\pi^*(h)}(Y_t)+\left\langle \int_0^t F(\bY_t)\,\dif\bX,h^1\star\pi^*(h^0)\right\rangle \\
                           &= F_{\pi^*(h)}(Y_t)+\langle F_{\pi^*(h^0)}(\bY_t),h^1\rangle\\
                           &= F_{\pi^*(h)}(Y)+\sum_{n\ge
                           1}^{\p-1}\frac{1}{n!}\D^nF_{\pi^*(h^0)}(Y_t)\left(\bY_{(h^1)^{(1)};t},\dotsc,\bY_{(h^1)^{(n)};t}\right)\\
                           &= F_{\pi^*(h)}(Y_t)+\sum_{n\ge
                           1}^{\p-1}\frac{1}{n!}\D^nF_{\pi^*(h^0)}(Y_t)\left(\widehat{F}_{(h^1)^{(1)}}(Y_t),\dotsc,\widehat{F}_{(h^1)^{(n)}}(Y_t)\right)\\
                           &= F_{\pi^*(h)}(Y_t)+\sum_{n\ge 1}^{\p-1}\frac{1}{n!}\left(
                           \widehat{F}_{(h^1)^{(1)}}\odot\dotsm\odot\widehat{F}_{(h^1)^{(n)}}\rhd F_{\pi^*(h^0)}
                         \right)(Y_t).
    \end{align*}

    To show the converse, assume that \(Y\) has the local expansion in \eqref{eq:davie.drift} and that \(\bY\) is a
    controlled path above it with components given by \(\langle \bY_t,h\rangle=\widehat{F}_h(Y_t)\).
    Fix \(n\ge 1\) and \(h\in\hgl^{(n)}\).
    By Taylor expansion and \eqref{eq:davie.drift} we can easily see that
    \begin{align*}
      \widehat{F}_h(Y_t)-\widehat{F}_h(Y_s)&\approx_{\p-n}\sum_{k=1}^{\p-n-1}\frac{1}{k!}\D^k\widehat{F}_h(Y_s)Y_{s,t}^{\otimes
      k}\\
                                           &\approx_{\p-n}\sum_{k=1}^{\p-n-1}\frac{1}{k!}\sum_{\scrf_1,\dotsc,\scrf_k\in\mathcal{F}_+}\left(\prod_{j=1}^k\varsigma(\scrf_k)\right)^{-1}\D^k\widehat{F}_h(Y_s)(\widehat{F}_{\scrf_1}(Y_s)\otimes\dotsm\otimes\widehat{F}_{\scrf_k}(Y_s))\bX_{s,t}^{\scrf_1\dotsm\scrf_k}\\
                                           &\approx_{\p-n}\sum_{\scrf\in\cF_+}\varsigma(\scrf)^{-1}(\mathbf{F}_\scrf\rhd\widehat{F}_h)(Y_s)\bX^{\scrf}_{s,t}\\
                                           &\approx_{\p-n}\sum_{\scrf\in\cF_+}\varsigma(\scrf)^{-1}\widehat{F}_{\scrf\star
                                           h}(Y_s)\bX_{s,t}^\scrf,
    \end{align*}
    that is, \(\bY\) satisfies \Cref{def:contr}.
  \end{proof}
  \begin{example}[Davie expansion for \(\p=3\)]\label{ex:Davie}
  Let us consider the RDE with drifts
  \def\qq{\frac12\left( \Forest{[\alpha]}\Forest{[\gamma[\beta]]}-\Forest{[\gamma]}\Forest{[\alpha[\beta]]} \right) }
  \[
    \dif
    Y_t=\begin{multlined}[t]F_{\Forest{[\alpha]}}(Y)\,\dif\bX_t^{\Forest{[\alpha]}}+ 
      \frac12 F_{\Forest{[\alpha]}\Forest{[\beta]}}(Y)\,\dif\bX_t^{\pi(\Forest{[\alpha])}\Forest{[\beta]})}+ \frac16
      F_{\Forest{[\alpha]}\Forest{[\beta]}\Forest{[\gamma]}}(Y)\,\dif\bX_t^{\pi(\Forest{[\alpha]}\Forest{[\beta]}\Forest{[\gamma]})}+F_{\qq}(Y)\,\dif\bX_t^{\pi\left(
    \Forest{[\alpha]}\Forest{[\gamma[\beta]]}\right) }.\end{multlined}
  \]
  Its solution has the following Davie expansion:
\begin{equation*}
	\begin{split}
		Y_{s,t} \approx{}&F_{\Forest{[\gamma]}} \bX_{s,t}^{\Forest{[\gamma]}} + \frac 12 F_{\Forest{[\alpha]}\Forest{[\beta]}} \bX_{s,t}^{\Forest{[\alpha]}\Forest{[\beta]}} + \big( \mathrm{D}F_{\Forest{[\beta]}}(F_{\Forest{[\alpha]}}) - F_{\Forest{[\alpha]}\Forest{[\beta]}} \big) \bX_{s,t}^{\Forest{[\beta[\alpha]]}} + \frac 16 F_{\Forest{[\alpha]}\Forest{[\beta]}\Forest{[\gamma]}} \bX_{s,t}^{\Forest{[\alpha]}\Forest{[\beta]}\Forest{[\gamma]}} \\
		&+ \frac 12 \bigg( F_{\frac 12 \big(\Forest{[\gamma]} \Forest{[\beta[\alpha]]} - \Forest{[\beta]} \Forest{[\gamma[\alpha]]} \big) } + F_{\frac 12 \big(\Forest{[\gamma]} \Forest{[\alpha[\beta]]} - \Forest{[\alpha]} \Forest{[\gamma[\beta]]} \big) } - \frac 12 \D F_{\Forest{[\beta]} \Forest{[\gamma]}} (F_{\Forest{[\alpha]}}) - \frac 12 \D F_{\Forest{[\alpha]} \Forest{[\gamma]}} (F_{\Forest{[\beta]}}) + \D F_{\Forest{[\gamma]}} (F_{\Forest{[\alpha]} \Forest{[\beta]}})
		\\
		&+ \D^2 F_{\Forest{[\gamma]}}(F_{\Forest{[\alpha]}}, F_{\Forest{[\beta]}})
		\bigg)
		\bX_{s,t}^{\Forest{[\gamma[\alpha][\beta]]}}+\bigg( F_{\frac 12 \big( \Forest{[\alpha]}\Forest{[\gamma[\beta]]} - \Forest{[\gamma]}\Forest{[\alpha[\beta]]}\big)} - \frac 12 F_{\Forest{[\alpha]}\Forest{[\beta]}\Forest{[\gamma]}} + \frac 12 \D F_{\Forest{[\alpha]}\Forest{[\gamma]}}(F_{\Forest{[\beta]}}) \bigg)\bX_{s,t}^{\Forest{[\alpha]}\Forest{[\gamma[\beta]]}}\\
		&+ \bigg( F_{\frac 12 \big( \Forest{[\beta]} \Forest{[\gamma[\alpha]]} - \Forest{[\gamma]} \Forest{[\beta[\alpha]]}\big)} + \frac 12 F_{\Forest{[\alpha]}\Forest{[\beta]}\Forest{[\gamma]}} - \frac 12 \D F_{\Forest{[\beta]} \Forest{[\gamma]}}(F_{\Forest{[\alpha]}}) - \D F_{\Forest{[\gamma]}}(F_{\Forest{[\alpha]} \Forest{[\beta]} }) + \D F_{\Forest{[\gamma]}}(\D F_{\Forest{[\beta]}} (F_{\Forest{[\alpha]}}) ) \bigg) \bX_{s,t}^{\Forest{[\gamma[\beta[\alpha]]]}}
	\end{split}
\end{equation*}
\end{example}

	We come to one of our main results. The purpose of an It\^o formula is to express functionals of a path, in our case smooth functions of RDE solutions, as integrals; see \Cref{expl:quasi3} below for the link with It\^o's celebrated formula for continuous semimartingales. A type of It\^o formula for branched rough paths that depends on additional structure is to be found in \cite[Ch.\ 5]{Kelly}, see \Cref{rem:kelly} below. We continue to implicitly use the identification \eqref{eq:ident}.

\begin{theorem}[Itô formula for RDEs with drifts]\label{thm:ito}
  Let \(\varphi \in C^{\p+1}(W,U)\). Then the
   path \(\varphi(Y)\) satisfies
  \[
    \varphi(Y_t) = \varphi(Y_s) + \int_s^t\mathbf{F}\varphi(Y_u)\,\mathrm{d}\bX_u.
  \]
\end{theorem}
\begin{proof}
  It is enough to check that
  \[
    \varphi(Y_t)-\varphi(Y_s)\approx\langle \mathbf{F}\varphi(Y_s),\bX_{s,t}\rangle.
  \]
  By Taylor's formula, we have
  \[
    \varphi(Y_t)-\varphi(Y_s)\approx\sum_{n=1}^{\p}\frac{1}{n!}\D^n\varphi(Y_s)Y_{s,t}^{\otimes n}.
  \]
  Since \(Y\) solves \eqref{eq:rde}, by \Cref{thm:davie.drift} we may substitute the expansion in \eqref{eq:davie.drift}
  for \(Y_{t}-Y_s\) to obtain
  \begin{align*}
    \varphi(Y_t)-\varphi(Y_s)&\approx\sum_{n=1}^\p\frac{1}{n!}\sum_{\mathscr{f}_1,\dotsc\mathscr{f}_n\in\mathcal{F}_+}\D^n\varphi(Y_s)\left(\widehat{F}_{\mathscr{f}_1}(Y_s)\otimes\dotsm\otimes\widehat{F}_{\mathscr{f}_n}(Y_s)\right)\bX_{s,t}^{\mathscr{f}_1\dotsm\mathscr{f}_n}\\
                             &\approx\sum_{\substack{\mathscr{f}\in\mathcal{F}_+\\|\mathscr{f}|\le\p}}\left(\sum_{n=1}^{|\mathscr{f}|}\frac{1}{n!}\D^n\varphi(Y_s)\left(\widehat{F}_{\mathscr{f}^{(1)}}\otimes\dotsm\otimes\widehat{F}_{\scrf^{(n)}}\right)\right)\bX_{s,t}^\scrf\\
                             &=\sum_{\substack{\mathscr{f}\in\mathcal{F}_+\\|\mathscr{f}|\le\p}}\mathbf{F}_\scrf\varphi(Y_s)\bX_{s,t}^\scrf\\
                             &=\left\langle\mathbf{F}\varphi(Y_s),\bX_{s,t}\right\rangle.\qedhere
  \end{align*}
\end{proof}

	\begin{remark}[It\^o formula for controlled paths]\label{rem:itoContr}
		In the extended setting of \(\mathcal(\cQ,W)\)-valued controlled integrands, it could be more generally stated that, given any \(\bH \in \mathscr D^{\p+1}_\bX(W)\) and \(\varphi\),
		\[
		\varphi(H_t) - \varphi(H_s) = \sum_{\scrf \in \scrH} \varsigma(\scrf)^{-1} \sum_{l \in \ell(\scrf)} \int_s^t \varphi(\bH)_{\pi^*(\scrf_l)} \dif \bX^{\pi(\scrf_l)},
		\]
		that is, one can get an Itô formula for an arbitrary controlled path at the cost of requiring \(\bH\) to have \say{an extra degree of controlledness}. This is not an
    issue with RDE solutions, which have enough Gubinelli derivatives, i.e., they always define elements of \(\mathscr{D}^{\p+1}_X(W)\).
	\end{remark}

	\begin{example}[It\^o formula for \(\p = 3\)]\label{ex:ito3}
    Continuing with \Cref{ex:Davie}, we see that for any smooth function \(\varphi\in C^\infty(W)\),
		\begin{align*}
      \varphi(Y_t)-\varphi(Y_s) = &\int_s^t \mathbf{F}_{\Forest{[\c]}}\varphi(Y_u) \dif \bX_u^{\Forest{[\c]}} \\
			&+ \frac 12 \int_s^t \mathbf{F}_{\Forest{[\a]}\Forest{[\b]}}\varphi(Y_u) \dif \bX_u^{\pi(\Forest{[\a]}\Forest{[\b]})} + \frac 16 \int_s^t \mathbf{F}_{\Forest{[\a]}\Forest{[\b]}\Forest{[\c]}}\varphi(Y_u) \dif \bX_u^{\pi(\Forest{[\a]}\Forest{[\b]}\Forest{[\c]})} \\
			&+ \int_s^t \mathbf{F}_{\frac 12 (\Forest{[\a]}\Forest{[\c[\b]]} - \Forest{[\c]}\Forest{[\a[\b]]})}\varphi(Y_u) \dif \bX_u^{\pi(\Forest{[\a]}\Forest{[\c[\b]]})}
		\end{align*}
		where in each integral, the sum is intended over labels, not labelled forests.
    Here, the differential operators are given explicitly by (using the notation introduced in \Cref{subsec:preLie})
    \begin{align*}
      \mathbf{F}_{\Forest{[\gamma]}} &= F_{\Forest{[\gamma]}} \\
      \mathbf{F}_{\Forest{[\alpha]}\Forest{[\beta]}} &=
      F_{\Forest{[\alpha]}\Forest{[\beta]}}+F_{\Forest{[\alpha]}}\odot F_{\Forest{[\beta]}} \\
      \mathbf{F}_{\Forest{[\alpha]}\Forest{[\beta]}\Forest{[\gamma]}} &=
      F_{\Forest{[\alpha]}\Forest{[\beta]}\Forest{[\gamma]}}+F_{\Forest{[\alpha]}}\odot
      F_{\Forest{[\beta]}\Forest{[\gamma]}}+F_{\Forest{[\beta]}}\odot
    F_{\Forest{[\alpha]}\Forest{[\gamma]}}+F_{\Forest{[\gamma]}}\odot
  F_{\Forest{[\alpha]}\Forest{[\beta]}}+F_{\Forest{[\alpha]}}\odot F_{\Forest{[\beta]}}\odot F_{\Forest{[\gamma]}}\\
      \mathbf{F}_{\frac12\left( \Forest{[\alpha]}\Forest{[\gamma[\beta]]}-\Forest{[\gamma]}\Forest{[\alpha[\beta]]}
  \right) } &= F_{\frac12\left(
    \Forest{[\alpha]}\Forest{[\gamma[\beta]]}-\Forest{[\gamma]}\Forest{[\alpha[\beta]]} \right)
  }+F_{\Forest{[\alpha]}}\odot\left( F_{\Forest{[\beta]}}\rhd F_{\Forest{[\gamma]}} \right)-F_{\Forest{[\gamma]}}\odot
  \left( F_{\Forest{[\beta]}}\rhd F_{\Forest{[\alpha]}} \right) +F_{\Forest{[\gamma]}}\odot
  F_{\Forest{[\alpha]}\Forest{[\beta]}}-F_{\Forest{[\alpha]}}\odot F_{\Forest{[\beta]}\Forest{[\gamma]}}.
    \end{align*}
    We further Davie-expand each integral:
		\begin{align*}
	\int_s^t \mathbf{F}_{\Forest{[\gamma]}}\varphi(Y) \dif \bX^{\Forest{[\gamma]}} \approx{}& \mathbf{F}_{\Forest{[\c]}}\varphi(Y_s) \bX_{s,t}^{\Forest{[\c]}} + \mathbf{F}_{\Forest{[\a]}} \mathbf{F}_{\Forest{[\b]}}\varphi(Y_s) \bX_{s,t}^{\Forest{[\a]} \tee \Forest{[\b]}} \\
	&+ \frac 12 \mathbf{F}_{\Forest{[\a]}\Forest{[\b]} }\mathbf{F}_{ \Forest{[\c]}}\varphi(Y_s) \bX_{s,t}^{\Forest{[\a]}\Forest{[\b]} \tee \Forest{[\c]}}+ \mathbf{F}_{\Forest{[\beta[\alpha]]}}\mathbf{F}_{ \Forest{[\gamma]}}\varphi(Y_s) \bX{}^{\Forest{[\beta[\alpha]]} \tee \Forest{[\gamma]}}_{s,t}  \\
	\int_s^t \mathbf{F}_{\Forest{[\a]}\Forest{[\b]}}\varphi(Y) \dif \bX^{\pi(\Forest{[\a]}\Forest{[\b]})} \approx{}& \mathbf{F}_{\Forest{[\a]}\Forest{[\b]}}\varphi(Y_s) \bX^{ \pi(\Forest{[\a]}\Forest{[\b]})}_{s,t} + \mathbf{F}_{\Forest{[\c]} }\mathbf{F}_{ \Forest{[\a]}\Forest{[\b]}}\varphi(Y_s) \bX^{\Forest{[\c]} \tee \pi(\Forest{[\a]}\Forest{[\b]})}_{s,t} \\
	\int_s^t \mathbf{F}_{\Forest{[\a]}\Forest{[\b]}\Forest{[\c]}}\varphi(Y) \dif \bX^{\pi(\Forest{[\a]}\Forest{[\b]}\Forest{[\c]})} \approx{}& \mathbf{F}_{\Forest{[\a]}\Forest{[\b]}\Forest{[\c]}}\varphi(Y_s) \bX^{ \pi(\Forest{[\a]}\Forest{[\b]}\Forest{[\c]})}_{s,t} \\
	\int_s^t \mathbf{F}_{\frac 12 (\Forest{[\a]}\Forest{[\c[\b]]} - \Forest{[\c]}\Forest{[\a[\b]]})}\varphi(Y) \dif \bX^{\pi(\Forest{[\a]}\Forest{[\c[\b]]})} ={}& \approx \mathbf{F}_{\frac 12 (\Forest{[\a]}\Forest{[\c[\b]]} - \Forest{[\c]}\Forest{[\a[\b]]})}\varphi(Y_s) \bX_{s,t}^{\pi(\Forest{[\a]}\Forest{[\c[\b]]})}
\end{align*}
		where we were able to simplify the antisymmetrisation in the last integral due to the fact that there are no higher derivatives (if \(\p\) were \(4\) or higher, this might not be possible, as is revealed by an inspection of the proof of \Cref{lem:intCoords}). Summing these four Davie expansions and writing the result in the forest basis, one obtains an expression that generalises that of \Cref{ex:Davie} to non-linear $\varphi$. The reader may check that extracting a dual pair of bases of \(\cP\) and \(\cQ\) (which is straightforward up to degree \(3\) by working with symmetric/antisymmetric tensors, see \Cref{expl:P3}, \Cref{ex:dual}) and writing the alternative expansion \eqref{eq:intBasis} one obtains the same result.
	\end{example}

	\begin{remark}[Kelly's approach via bracket extensions]\label{rem:kelly}
		As already mentioned, a version of \Cref{thm:ito} was already obtain by Kelly \cite[Theorem 5.3.11]{Kelly} in the case of \(F \in C^\infty(W,\cL(V,W))\). The main difference between the two approaches is that we leverage an algebraic property of \(\hgl\) --- its freeness over \(\cQ\) --- in order to write a change-of-variable formula that does not depend on anything that is not already contained in the original branched rough path \(\bX\). On the other hand, Kelly introduces the notion of \emph{bracket extension}, a branched rough path whose trace takes its values in an enlarged space, which in our notation would be \(\cH(V)\): that is, while the original \(\bX\) is indexed by forests labelled with \(\alpha \in V\), the extended \(\widehat\bX\) is indexed by forests \emph{indexed with by forests labelled with \(\alpha \in V\)}. These new labels \(\scrf\) actually represent the primitives \(\pi(\scrf)\), so that \(\widehat \bX\) is meant to encode branched integrals indexed by primitives (which in our approach is handled by \eqref{eq:Tint}). Beyond extending \(\bX\), \(\widehat\bX\) is required to satisfy certain \say{bracket polynomial} relations from which it follows that a result analogous to the theorem above with integration against \(\dif \widehat \bX^{(\scrf_1 \cdots \scrf_n)}\). The bracket extension \(\widehat\bX\), however, is highly non unique, and only shown to exist using non-constructive methods, which have been avoided here. Another point worth noting is that, in the bracket extension approach, \(\bX\)-controlled paths can be expanded as integrals against \(\widehat \bX\); however, due to redundancies introduced by the additional labels, \(\widehat \bX\)-controlled paths would need a further iteration of the bracket extension to be represented as integrals; this, too, is avoided in our approach, since we do not need to extend \(\bX\), only re-organise its internal structure. We remark that Kelly's change of variable formula takes a slightly simplified expression in coordinates compared to ours; in particular the \(\pi(\scrf)\) that indexes the integrator for us is a separate, newly added label for him, and \(\scrf\) replaces \(\pi^*(\scrf)\) as subscript of the integrand: this substitution, which appears not to be possible in our case (see a special case where it is, in \Cref{ex:ito3}), jointly with the fact that \(F\) is of simplified form, makes it possible for him to sum over trees.

		The two approaches can be reconciled by observing that setting, recursively, for \(\scrf \in \scrH(\scrH(V))\) and \(\scrg \in \scrH(V^*)\)
		\begin{equation}\label{eq:Tbracket}
			\langle [\scrf]_{(\scrg)}, \widehat \bX \rangle = \langle \scrf \tee \pi(\scrg), \widehat \bX \rangle
		\end{equation}
		defines a \emph{canonical} bracket extension: the bracket relations to be checked read
		\[
		\langle \pi(\scrg) , \widehat \bX \rangle = \langle \scrg - \scrg' \tee \pi(\scrg''), \widehat \bX \rangle,\qquad \scrg \in \scrH(V)
		\]
		which is the content of \Cref{prop:pirec}. One must also check that \eqref{eq:Tbracket} defines a branched rough path, which is implied by the fact that the algebra morphism \(\top \colon \hck(\cH(V^*)) \to \hck(V^*)\) defined recursively by
		\[
		[\scrf]_{(\scrg)} \mapsto \top(\scrf) \tee \pi(\scrg),\quad \scrf \in \scrH(\scrH(V^*)), \ \scrg \in \scrH(V^*)
		\]
		is a Hopf morphism. This is again checked thanks to the cocycle property \eqref{eq:cocycle}:
		\begin{align*}
			\widetilde\Delta_\mathrm{CK} \top ([\scrf]_{(\scrg)}) &= \top(\scrf) \otimes \pi(\scrg) + \top(\scrf)' \otimes (\top(\scrf)'' \top \pi(\scrg)) \\
			&= \top(\scrf) \otimes \pi(\scrg) + \top(\scrf') \otimes (\top(\scrf'') \tee \pi(\scrg)) \\
			&= \top(\scrf) \otimes \top(\bullet^{(\pi(\scrg))}) + \top(\scrf') \otimes \top([\scrf'']_{(\pi(\scrg))}) \\
			&= (\top \otimes \top) \widetilde \Delta_\mathrm{CK} ([\scrf]_{(\scrg)})
		\end{align*}
		where in the second step we have argued inductively since \(\scrf\) is of lower degree in the coradical filtration, jointly with the observation that \(\dck\) is an algebra morphism, so that (taking \(\scrf = \scrt_1 \cdots \scrt_n\))
		\[
		\widetilde \Delta_\mathrm{CK} \top (\scrt_1 \cdots \scrt_n) = \widetilde \Delta_\mathrm{CK} ( \top (\scrt_1) \cdots \top(\scrt_n)) = \top^{\otimes 2} \widetilde \Delta_\mathrm{CK} (\scrt_1) \cdots \widetilde \Delta_\mathrm{CK} (\scrt_n)
		\]
		The fact that our approach and Kelly's can be reconciled in this way can be helpful when bootstrapping statements in the extended setting that are already known in the classical branched setting: for example, local existence and uniqueness of the RDEs driven by \(X^{p \in \cP}\) can be inferred from the fact that these are RDEs (in the ordinary branched sense \cite[Theorem 8.8]{Gub10}) driven by this canonical bracket extension.
	\end{remark}
	
	\subsection{A criterion for quasi-geometricity and the simple change-of-variables formula}\label{subsec:quasi}
We begin by motivating this section by considering the restriction of \Cref{thm:ito} to the \say{trivial RDE}, from
which we obtain the It\^o formula for functions of \(X\). The proof is a straightforward application of the recursive
definition of \(\mathbf{F}\) in \Cref{thm:ito}. We recall the notation \(\mathscr r(\gamma_1,\ldots,\gamma_n) \coloneqq \Forest{[\gamma_1]} \cdots \Forest{[\gamma_n]}\).
\begin{corollary}[Simple change of variable formula]\label{cor:simple}
	Let \(\bX \in \scrC^\rho([0,T],V)\), \(\cQ(V)=W=V\) and \(F(y) \equiv \mathrm{id}_\cQ\) in \eqref{eq:rde}. Then, for
  \(\varphi \in C^{\p+1}(V)\) we have
	\[
    \mathbf{F}_h \varphi (Y) = \sum_{n \geq 1} \frac{1}{n!} \partial_{\pi^*(h^{(1)}) \ldots \pi^*(h^{(n)})}\varphi(X)
	\]
  for any \(h\in\hgl^{\p}\).
  If, furthermore, \(W = V\), \(\mathbf{F}_{\scrr(\gamma_1,\ldots,\gamma_n)} \varphi(Y) = \partial_{\gamma_1 \ldots
  \gamma_n}\varphi(X)\) and \(\Fbar_h \varphi = 0\) otherwise, so that
	\begin{equation}\label{eq:itoSimple}
	\varphi(X_t)-\varphi(X_s) = \sum_{n = 1}^\p \frac{1}{n!}\int_s^t \partial_{\gamma_1\ldots \gamma_n}\varphi(X_u)\,\dif
  \bX_u^{\pi(\scrr(\gamma_1,\ldots,\gamma_n))}
	\end{equation}
\end{corollary}
The use of the term \emph{simple} is due to Kelly \cite{Kelly}, who first pointed out that the terms of a bracket
extension needed to make sense of \eqref{eq:itoSimple} were fewer than those needed for the general RDE case; the same
formula is available in the quasi-geometric case \cite[Theorem 4.20]{Bel20}. There is a setting in which the general
change-of-variables formula of \Cref{thm:ito} takes on the simpler form \eqref{eq:itoSimple}. Hoffman's quasi-shuffle
algebra, introduced in \cite{Hof00} to study renormalisation of multiple zeta values, has been connected with rough path
theory \cite{Bel20} and stochastic differential equations \cite{EMPW15}. We refer to \cite{Reu93} for an introduction to the shuffle algebra, and to \cite{Hof00, HI17} for an introduction to the quasi-shuffle algebra, and assume familiarity with these notions. In the spirit of the approach taken here we view these as covariant functors (see \cite{NJPT13,FPT16} for a similar perspective on quasi-shuffle algebras and their deformations)	
\begin{equation}
	\begin{split}	
		\bigshuffle(U) = (\tens(U), \shuffle,  \Delta_\otimes), \qquad &\widetilde\bigshuffle(U) = (\tens({\textstyle \bigodot}(U)), \widetilde\shuffle,  \Delta_\otimes)
	\end{split}
\end{equation}
We are considering a free commutative bracket for our quasi-shuffle algebra, i.e.\ \([x,y]\) in the notation of \cite{Hof00} is given by \(x \odot y\) here (\(x,y \in \bigodot (U)\)).

A first direct link with $\hck$ is provided by the arborifications maps, see \cite[Section 2.3]{BCE20} and \cite[Section 4.1]{HK15}.
\begin{proposition}[(Contracting-)arborification maps]\label{prop:arb}
There exist unique surjective Hopf algebra morphisms \(\mathfrak a \colon \hck(U) \to \bigshuffle(U)\),  \(\widetilde{\mathfrak a}_\odot \colon \hck(\bigodot ( U)) \to \widetilde\bigshuffle(U)\) defined recursively by the condition
	\[
\mathfrak a ([\scrf]_r) \coloneqq \mathfrak a (\scrf) \otimes r\,. 
	\]
with $r \in U$ and $r \in \bigodot (U)$ respectively. Following the literature, we call them the \textnormal{arborification} and \textnormal{contracting-arborification} maps respectively. These maps are natural and are respectively left inverses to the natural injective coalgebra maps \(\iota \colon \bigshuffle (U) \to \hck(U)\), \(\widetilde \iota_\odot \colon \bigqshuffle (U) \hookrightarrow \hck ( \bigodot(U))\) mapping words to decorated ladders, with the last letter corresponding to the root.
\end{proposition}
This result shows how \(\bigshuffle (U)\) and \(\bigqshuffle (U)\) can be obtained as natural quotients of the Connes-Kreimer Hopf algebra over \(U\) and \(\bigodot(U)\) respectively.
However, in order to relate quasi-geometric and branched rough paths, we need to represent \(\bigqshuffle (U)\) as a quotient of \(\hck(U)\), not \(\hck(\bigodot(U))\): for this reason, from now on and unless otherwise mentioned, we revert to working over \(U\).

We denote by \(\cQ_\scrr\) the subspace of \(\cF\) spanned by elements of the form \(\scrr(\gamma_1,\ldots,\gamma_n)\), and as remarked in \Cref{ex:dual}, \(\cQ_r \subset \cQ\); similarly, denote \(\mathcal{P}_\scrr \coloneqq \pi(\cQ_\scrr)\). Note that there is an obvious, natural isomorphism \(\cQ_r(U) \cong \bigodot(U)\). Using the definition of \(\pi\), we have the following recursive formula:
\begin{equation}\label{form_basis_prim}
\pi (\scrr(\c_1,\ldots,\c_n))= \scrr(\c_1,\ldots,\c_n) - \sum_{I\sqcup \,J= \{\c_1\,, \ldots \,, \c_n\}} \scrr(\c_I)\tee \pi(\scrr(\c_J))\,,
\end{equation}
where \(I\) and \(J\) are two non-empty subsets partitioning \(\{\c_1, \ldots , \c_n\}\) and \(\scrr(\gamma_I)\), \(\scrr(\gamma_J)\) are the associated forests obtained decorating the roots with elements indexed by \(I\) and \(J\) respectively.
\begin{lemma}\label{lemma_independence}
\(\pi \colon \cQ_\scrr \to \cP_\scrr\) is an isomorphism.
\end{lemma}
\begin{proof}
\eqref{form_basis_prim} implies that, given a basis \(B\) of \(U\), \(\{\pi(\scrr(\c_1,\ldots,\c_n)) \colon \c_1,\ldots, \c_n \in B\}\) is a basis of \(\cP_\scrr\). Indeed, \(\{\scrr(\c_1,\ldots,\c_n) \colon \c_1,\ldots, \c_n \in B\}\) is a basis of \(\cQ_\scrr\), and no non-trivial linear combination of the terms of the form \(\scrr(\c_I)\tee \pi(\scrr(\c_J))\) can belong to \(\cQ_\scrr\) since all forests involved have at least one edge.
\end{proof}

We are now in a position to define a version of the quasi-arborification on \(\hck(U)\) instead of on \(\hck(\bigodot (U))\) and to construct a natural right inverse to it.

\begin{theorem}[Intrinsic quasi-arborification map]\label{prop:arb.intr}
There exists a surjective Hopf algebra morphism, the \emph{intrinsic quasi-arborification} map \(\widetilde{\mathfrak a} \colon \hck(U) \to \widetilde\bigshuffle(U)\), defined recursively by the condition
	\[
\widetilde{\mathfrak a} ([\scrf]_\gamma) \coloneqq \widetilde{\mathfrak a} (\scrf)\otimes \gamma\,. 
	\]
\(\widetilde{\mathfrak a}\) is natural and is left inverse to the natural coalgebra map \(\widetilde \iota \colon \bigqshuffle (U) \hookrightarrow \hck ( U)\) defined by \(\c_1 \odot \cdots \odot \c_n \mapsto \pi(\scrr(\c_1,\ldots,\c_n))\) and extended by \(s_1 \otimes \cdots \otimes s_n \mapsto s_1 \tee \ldots \tee s_n\) for \(s_1,\ldots,s_n \in \bigodot(U)\).
\end{theorem}
\begin{proof}
Existence and uniqueness of an algebra morphism is clear from the defining property and the fact that algebra morphism property forces \(\widetilde{\mathfrak a}(\scrt_1 \cdots \scrt_n) = \widetilde{\mathfrak a}(\scrt_1) \qshuffle \ldots \qshuffle \widetilde{\mathfrak a}(\scrt_n)\). \(\widetilde{\mathfrak a}\) is a coalgebra morphism by the universal property of the Connes-Kreimer Hopf algebra (see for example \cite[Theorem 3]{Foi13}), since \(\widetilde{\mathfrak a} \circ \EuScript{B}_+^\gamma = (- \otimes \gamma) \circ \widetilde{\mathfrak a}\).
We now make use of the following property of quasi-shuffles, which can be verified directly from the definition of \(\qshuffle\): if \(w_1,\ldots,w_n \in \bigqshuffle(U)\) are words (i.e.\ elementary tensors) at least one of which has length strictly greater than \(1\) (where length refers to length in the tensor algebra, disregarding weight of the letters in \(\bigodot(U)\)), then \(w_1 \qshuffle \ldots \qshuffle w_n\) is a linear combination of words of length strictly greater than \(1\). This implies that 
\begin{equation}\label{eq:edgeToZero}
\pi_1 \circ \widetilde{\mathfrak a}(\scrf) = 0 \ \text{if} \ \scrf \ \text{is a forest containing at least one edge}
\end{equation}
We then have, by \Cref{coalgebra_prop} (and using notation therein), for \(p_1,\ldots,p_n \in \cP\)
\begin{align*}
	\widetilde{\mathfrak a}(p_1 \tee \ldots \tee p_n) &= \sum_{\substack{n_1 + \ldots + n_k = n \\ n_1,\ldots,n_k \geq 1}} \pi_1 \circ\widetilde{\mathfrak a}(p_1 \tee \ldots \tee p_{n_1}) \otimes \cdots \otimes \pi_1 \circ\widetilde{\mathfrak a}(p_{n^{k-1}+1} \tee \ldots \tee p_n) \\
	&= \widetilde{\mathfrak a}(p_1) \otimes \cdots \otimes \widetilde{\mathfrak a}(p_n)
\end{align*}
once again because \(\tee\) always introduces edges, and thus length strictly greater than \(1\) in the corresponding images through \(\widetilde{\mathfrak a}\). Finally, using \eqref{form_basis_prim}
\begin{align*}
\widetilde{\mathfrak a} \circ \widetilde \iota (\c_1 \odot \cdots \odot \c_n) &=  \widetilde{\mathfrak a} \circ \pi(\scrr(\gamma_1,\ldots,\gamma_n)) \\
 &= \pi_1 \circ \widetilde{\mathfrak a} \circ \pi(\scrr(\gamma_1,\ldots,\gamma_n)) \\
&= \pi_1 \circ \widetilde{\mathfrak a} \circ \scrr(\gamma_1,\ldots,\gamma_n) \\
&= \pi_1( \c_1 \qshuffle \cdots \qshuffle \c_n ) \\
&= \c_1 \odot \cdots \odot \c_n
\end{align*}
and since \(\widetilde{\mathfrak a} \circ \widetilde \iota (s_1 \otimes \cdots \otimes s_n) = \widetilde{\mathfrak a}(\widetilde\iota(s_1) \tee \ldots \tee \widetilde\iota(s_n)) = s_1 \otimes \cdots \otimes s_n\) (as \(\widetilde\iota\) maps \(\bigodot(U)\) to \(\cP_\scrr\)), we conclude that \(\widetilde{\mathfrak a} \circ \widetilde\iota = \mathrm{id}\).
\end{proof}

We note that the intrinsic arborification map, despite admitting a right inverse that is defined in terms of \(\top\), is not itself defined in terms of the \(\binf\)-structure of \(\hck\). The fact that \(\bigshuffle\) and \(\bigqshuffle\) are quotients of \(\hck\) immediately implies that geometric and quasi-geometric rough paths --- rough paths defined on these Hopf algebras in total analogy to \Cref{def:brp} --- can be viewed as branched rough paths.

 Conversely, one may ask whether a branched rough path comes from a (quasi-)geometric one in this manner. Before we
 answer this question, we provide some intuition as to what the answer entails from the point of view of integration.
 Just in \eqref{eq:Tint} and the discussion below, we showed that the \(\tee\) operation could be thought of as the algebraic counterpart to integration, we now show how \(\pi\) can be thought of the algebraic counterpart to a generalised integration-by-parts formula. From \Cref{prop:pirec} it follows that for \(\scrf = [\scrg_1]_{\gamma_1} \cdots [\scrg_n]_{\gamma_n}\), the product of the \say{integrals} (in the interpretation \eqref{eq:brPostulate}) \(\langle {[\scrg]_{\gamma_k}},\bX^{[\scrg]_{\gamma_k}}\rangle\) can be written as a sum of \say{integrals} as
\begin{equation}\label{eq:ibp}
	\begin{split}
	\langle\scrf, \bX_{s,t}\rangle ={} &\sum_{(\scrf)} \langle \scrf_{(1)} \tee \pi(\scrf_{(2)}), \bX_{s,t} \rangle \\
	={}&\sum_{k = 1}^n \big\langle \big[ [\scrg_1]_{\gamma_1} \cdots [\scrg_{k-1}]_{\gamma_{k-1}} \cdot \scrg_k \cdot [\scrg_{k+1}]_{\gamma_{k+1}} \cdots [\scrg_n]_{\gamma_n} \big]_{\gamma_k} , \bX_{s,t} \big\rangle \\
	&+\sum_{\substack{I \sqcup J = [n] \\ |I| \geq 2}} \big\langle \prod_{i \in I}\scrg_i  \cdot \prod_{j \in J} [\scrg_j]_{\gamma_j} \big) \tee \pi \big( \prod_{i \in I} \Forest{[\gamma_i]} \big),\bX_{s,t}  \big\rangle \\
	&+\text{terms involving cuts above the roots of the forests \(\scrg_k\)}
	\end{split}
\end{equation}
Truncating the last expression after its first line sign is the classical integration-by-parts formula that holds for smooth paths, which continues to hold for geometric rough paths (which, despite their roughness, obey the same laws as ordinary calculus). Stochastic examples of geometric rough paths are given by Stratonovich integration of a semimartingale. For quasi geometric rough paths, instead, one must stop after the second line: this means that \say{interactions} involving \(n\) of the integrators (\say{\(n\)-variations}) are non-negligible. All \(2 < \rho\)-branched rough paths are automatically quasi-geometric (this will be obvious from the next result), the It\^o integral being the prime stochastic example. The correct integration-by-parts for general branched rough path, however, involves further terms (third line), \say{interactions between integrals} not captured by those mentioned previously.  We note that a characterisation of \(\bigodot(V)\)-valued branched rough paths that are quasi-geometric was provided in \cite[Theorem 3.1]{F22}. This is different from a characterisation of \(V\)-valued branched rough paths that are quasi geometric which does not require the extension of the rough path to \(\bigodot(V)\) to be checked, and which is only possible to state thanks to our intrinsic arborification map above.

\begin{corollary}[Characterisation of (quasi-)geometric rough paths]\ \\[-2ex] \label{Thm_quasi_geo}
	\begin{itemize}
		\item A \(V\)-valued branched rough path \(\bX\) is geometric (i.e.\ lies in the image of \(\mathfrak a^*\)) if and only if, for \(p_1,\ldots,p_n \in \cP\), \(\langle p_1 \tee \ldots \tee p_n, \bX \rangle = 0\) if \(\exists k \ p_k \in \cP/\iota(V^*)\) (where the quotient is identified with the direct complement of \(\iota(V^*)\) in \(\cP\) generated by forests that have more than one vertex);
		\item A \(V\)-valued branched rough path \(\bX\) is quasi-geometric (i.e.\ lies in the image of \(\widetilde{\mathfrak a}{}^*\)) if and only if, for \(p_1,\ldots,p_n \in \cP\), \(\langle p_1 \tee \ldots \tee p_n, \bX \rangle = 0\) if \(\exists k \ p_k \in \cP/\cP_\scrr\) (where the quotient is identified with the direct complement of \(\cP_\scrr\) in \(\cP\) generated by forests that have an edge).
	\end{itemize}
\end{corollary}
\begin{proof}
We only prove the second assertion; the first is proved similarly (using \Cref{prop:arb}). Note that if \(\bX = \widetilde{\mathfrak a}{}^*(\bZ) = \widetilde{\mathfrak a}{}^* \circ \widetilde \iota^* \circ \widetilde{\mathfrak a}{}^* (\bZ)\), then \(\bX = \widetilde{\mathfrak a}^* \circ \widetilde \iota^*(\bX)\). By \Cref{thm:T}, the elements \(p_1 \tee \ldots \tee p_n\) with \(p_k \in \cP\) span \(\hck\), so by non-degeneracy of the pairing \Cref{thm:dual}, \(\bX = \widetilde{\mathfrak a}{}^* \circ \widetilde \iota^*(\bX)\) if and only if
\begin{align*}
\langle p_1 \tee \ldots \tee p_n , \bX \rangle &= \langle \widetilde\iota \circ \widetilde{\mathfrak a} (p_1 \tee \ldots \tee p_n) , \bX \rangle = \langle \widetilde\iota \circ \pi_1 \circ \widetilde{\mathfrak a}(p_1) \tee \ldots \tee \widetilde\iota \circ \pi_1 \circ \widetilde{\mathfrak a}(p_n) , \bX \rangle
\end{align*}
and the assertion now follows from \eqref{eq:edgeToZero}.
\end{proof}

We now state the restriction of \Cref{thm:ito} to the case of quasi-geometric rough paths (in the sense of \Cref{Thm_quasi_geo}). An important difference with the general branched case is that no derivatives of the vector fields appear in the integral form of the change of variable formula. We use bold letters to denote tuples of elements in a basis \(B\) of \(V\).

\begin{corollary}[Quasi geometric change-of-variable formula for RDEs]
	In the notation of \Cref{thm:ito}, if \(\bX\) is quasi-geometric the formula reduces to
	\[
	\varphi(Y_t)-\varphi(Y_s)= \sum_{\substack{|\boldsymbol \gamma^1|, \ldots, |\boldsymbol \gamma^n| \geq 1 \\ |\boldsymbol
  \gamma^1|, \ldots, |\boldsymbol \gamma^n| \leq \p}} \frac{1}{n!|\boldsymbol \gamma^1|! \cdots |\boldsymbol \gamma^n|!}
  \int_s^t \partial_{k_1 \ldots k_n} \varphi(\bY) F{}^{k_1}_{\mathscr r(\boldsymbol{\gamma}^1)}(\bY) \cdots
  F{}^{k_n}_{\mathscr r(\boldsymbol{\gamma}^n)}(\bY)\,\dif \bX^{\pi(\scrr(\boldsymbol{\gamma}^1\ldots\boldsymbol{\gamma}^n))} \,.
	\]
	We consider two special cases of this: first, if \(F \in C^\infty(W,\cL(V,W))\), 
	\[
		\varphi(Y_t)-\varphi(Y_s) = \sum_{\substack{\gamma_1,\ldots,\gamma_n \\ 1 \leq n \leq \p}} \frac{1}{n!} \int_s^t \partial_{k_1
  \ldots k_n}\varphi(\bY) F^{k_1}_{\gamma_1}(\bY) \cdots F^{k_1}_{\gamma_1}(\bY)\,\dif \bX^{\pi(\scrr(\c_1 \ldots \c_n))} \,.
	\]
	Secondly, if \(W,F\) are as in \Cref{cor:simple},
	\[
		\varphi(X_t)-\varphi(X_s)= \sum_{\substack{|\boldsymbol \gamma^1|, \ldots, |\boldsymbol \gamma^n| \geq 1 \\ |\boldsymbol
  \gamma^1|, \ldots, |\boldsymbol \gamma^n| \leq \p}} \frac{1}{n!|\boldsymbol \gamma^1|! \cdots |\boldsymbol \gamma^n|!}
  \int_s^t \partial_{\scrr(\boldsymbol \gamma^1) \ldots  \scrr(\boldsymbol\gamma^n)} \varphi(\bX)\,\dif \bX^{\pi(\mathscr r(\boldsymbol{\gamma}^1\ldots\boldsymbol{\gamma}^n))}\,
	\]
\end{corollary}
 The result follows straightforwardly from \Cref{thm:ito}, using the fact that \(\cQ_\scrr\) is a sub-coalgebra of
 \(\hgl\) to write the integrals explicitly in terms of the vector fields \(F\), without having to resort to the more
 complicated general expression of \(\mathbf{F}\).

\begin{remark}[Simple change of variable formula and quasi-geometric rough paths]
	We note that the intersection of the two special cases corresponds to \Cref{cor:simple} (restricted to the
  quasi-geometric setting, already shown in \cite[Theorem 4.20]{Bel20}). Since \(\widetilde \iota^* \colon \hgl \twoheadrightarrow \bigqshuffle\) is not a Hopf algebra morphism, \(\widetilde \iota^*(\bX)\) cannot be expected to define a rough path. Nevertheless, the integration of \eqref{eq:itoSimple} only depends on this projection: this is because, in the Davie expansion of that integral, \(\bX\) is only evaluated on linear combinations of terms of the form \((\Forest{[\alpha_1]} \cdots \Forest{[\alpha_m]}) \tee \pi(\Forest{[\beta_1]} \cdots \Forest{[\beta_n]})\); using \Cref{prop:pirec} it is then seen inductively that such elements belong to \(\top(\bigotimes \bigodot (V^*))\), i.e.\ the range of \(\widetilde\iota\).
\end{remark}

\begin{example}[It\^o/Stratonovich calculus and obstructions to quasi-geometricity]\label{expl:quasi3}
	The degree to which a branched rough path is free to be non (quasi-)geometric is constrained, to a certain extent, by
  its regularity. When \(\rho < 2\) there is only one rough path, and it is geometric. When \(2 \leq \rho < 3\), a rough
  path \(\bX\) is always quasi-geometric and geometric if and only if \(\langle \pi(\Forest{[\a]}\Forest{[\b]}), \bX
  \rangle = 0\): in the notation of \cite{FH20}, this corresponds to \([\bX] = 0\). Recall that if \(X\) is a continuous
  semimartingale, this condition is satisfied when \[\bX_{s,t}^{\Forest{[\alpha[\beta]]}} \coloneqq \int_s^t (X^\a_{u}- X^\a_{s})\circ \dif
  X^\b_u\]
   is defined by Stratonovich integration, but when it is defined by It\^o integration 
  \[\bX_{s,t}^{\Forest{[\alpha[\beta]]}} \coloneqq   \int_s^t(X^\a_{u}- X^\a_{s})\,\dif X^\b_u\]
  it holds that \(\langle \pi(\Forest{[\a]}\Forest{[\b]}), \bX \rangle = [X]\), the quadratic variation of \(X\). We take the opportunity to remark that the framework laid out in this and the previous section makes it possible to recover the classical It\^o SDE with drift \Cref{eq:SDE} from the equation \(\dif Y = F(Y)\,\dif \mathbf W\), with \(F_\gamma = \sigma_\gamma\) and \(F_{\scrr(\alpha,\beta)} = \frac{\deltaup_{\alpha\beta}}{2} \mu\). This is because \(\mathbf W^{\pi(\scrr(\alpha,\beta))}_t = [W]_t^{\alpha\beta} = \deltaup^{\alpha\beta}t\), and the dependence of the coefficients on \(t\) is easily handled by coupling the equation with the trivial equation \(\dif t = \dif \mathbf W^{\pi(\scrr(1,1))}_t\).
  Moving on to higher-order cases, the only obstructions for quasi-geometricity in the case \(3 \leq \rho < 4\) is
	\begin{equation}
		\langle \pi(\Forest{[\a]}\Forest{[\c[\b]]}) , \bX \rangle = 0
	\end{equation}
	by \Cref{expl:P3}. Expressions get increasingly complicated as \(\rho\) increases, but the quasi-geometricity condition is still tractable for \(4 \leq \rho < 5\) and \(\mathrm{dim}(V) = 1\), in which case (see \Cref{prof:P4Q3}) it reads 
	\begin{equation}\label{eq:fourQuasi}
		\langle \pi(\Forest{[[]]}\Forest{[[]]}) , \bX \rangle = 0
	\end{equation}
	since the previous order-\(3\) condition is automatically cleared.
\end{example}

\section{The It\^o-Stratonovich correction formula}\label{sec:itostrat}

	\subsection{Commutative \texorpdfstring{$\mathbf B_\infty$}{B∞}-algebras are shuffle algebras}\label{subsec:iso}
The main purpose of this section is to define an explicit, natural Hopf isomorphism between the Connes-Kreimer Hopf algebra and the shuffle algebra over its primitive elements, and to use this isomorphism to write an RDE driven by a branched rough path \(\bX\) as an equivalent RDE driven by a geometric rough path, defined algebraically in terms of \(\bX\), and with the new vector fields defined algebraically in terms of the old.
	
Recall that on any Hopf algebra \((H, \times , \Delta)\) the \emph{convolution product} on linear endomorphisms is defined as
\[
f * g \coloneqq (\argument \times \argument) \circ (f \otimes g) \circ \Delta, \qquad f,g \in \cL(H,H)
\]
making \((\cL(H,H),*)\) an associative algebra. For the following definition and lemma, we follow \cite{P93} (see also \cite{AL15}).

\begin{definition}[Eulerian idempotent]\label{def:euler}
Given a graded connected Hopf algebra \((H, \times , \Delta)\), which we assume to be commutative or cocommutative, we define its Eulerian idempotent \(\e_H \in \cL(H,H)\)
\begin{align*}
	\e_H &\coloneqq \log_*(\mathrm{id}) \\
	&\coloneqq \sum_{n \geq 1} \frac{(-1)^{n-1}}{n} (\mathrm{id} - \eta_H\circ\varepsilon_H)^{*n} \\
	&= \sum_{n \geq 1} \frac{(-1)^{n-1}}{n} \times^{(n-1)} \circ (\mathrm{id} - \eta_H\circ\varepsilon_H)^{\otimes n} \circ \Delta^{(n-1)} \\
	&= \sum_{n \geq 1} \frac{(-1)^{n-1}}{n} \times^{(n-1)} \circ \widetilde \Delta^{(n-1)}
\end{align*}
where \(\varepsilon_H\) denotes the counit \(H \twoheadrightarrow \bbR\) and \(\eta_H\) denotes the unit
\(\bbR\hookrightarrow H\), \(\times^{(n-1)}\) is the product taking \(n\) arguments, and \(\widetilde \Delta\) the reduced coproduct.
\end{definition}
 In Sweedler notation the map \(\e_H\) reads
\begin{equation}\label{eq:e}
\e_H(h) = \sum_{n \geq 1} \frac{(-1)^{n-1}}{n} h^{(1)} \times \cdots \times h^{(n)}.
\end{equation}
We summarise a few properties of \(\e_H\), which distinguish between the case of \(H\) commutative or cocommutative
below \cite{P93}.

\begin{proposition}[Properties of the Eulerian idempotent]\label{prop:e}
Let \(H\) be as above.
\begin{enumerate}
	\item Setting 
	\[
	\e_n \coloneqq \frac{1}{n!} \e_H^{*n}, \qquad \text{it holds that} \quad H = \bigoplus_{n = 0}^\infty \e_n(H)
  \quad\text{and}\quad \e_m \circ \e_n = \deltaup_{mn}\e_n. 
	\]
  In particular, \(\e_n\) is idempotent for every \(n\ge 1\).
	\item Let \(H^{\circ}\) be the graded dual Hopf algebra of \(H\). Then \(\e_H^{\circ} = \e_{H^{\circ}}\).
	\item \(\mathrm{Prim}(H)\subseteq \e(H)\) and equality holds if \(H\) is cocommutative.
	\item If \(H\) is cocommutative, the inclusion \(\e(H) \hookrightarrow H\) induces an algebra isomorphism \(\mathcal U(\e(H)) \cong H\), where \(\mathcal U\) denotes the universal enveloping algebra functor.
	\item If \(H\) is commutative, setting \(H_+ \coloneqq \ker\varepsilon_H\) and \(H_+^{\times 2}\) the space generated
    by products of elements in \(H_+\), \(\ker(\e_H) = H_+^{\times 2}\), and the resulting quotient map splits the short exact sequence
	\[
\begin{tikzcd}
	0 \arrow[r] &H_+^{\times 2} \arrow[r] &H_+ \arrow[r] &H/H_+^{\times 2} \arrow[r] \arrow[l,bend left,start
  anchor={[xshift=1.5ex,yshift=0.5ex]},"\e_H"] &0
\end{tikzcd}.
\]
	\item If \(H\) is commutative, \(\e_H(H) \hookrightarrow H\) induces an algebra isomorphism \(\bigodot (\e_H(H)) \cong H\), where \(\bigodot\) denotes the symmetric tensor algebra functor.
\end{enumerate}
\end{proposition}

We briefly comment on those aspects that are not difficult to show, and the consequences in the context considered here.
Idempotence and orthogonality can be shown using the fact that $\mathrm{id}^{* p} \circ \mathrm{id}^{* q} = \mathrm{id}^{* pq}$ and that $\mathrm{id}^{* p} = \sum_{k \geq 1} p^k \e_k$ in a graded fashion, see \cite{Lod94}. The decomposition of \(H\) as a direct sum involves writing
\[
\mathrm{id} = \exp_*\circ\log_*(\mathrm{id}) = \sum_{n = 0}^\infty\e_n
\]
Assertion (2) is an easy consequence of the definition of \(\e_H\) and duality of the bialgebra operations; we will be
using it for the dual pair \((\hck,\hgl)\), i.e.\ \(\e_\mathrm{CK}^\circ = \e_\mathrm{GL}\).
As for (3), the only non-zero term in \eqref{eq:e} with \(h \in \mathrm{Prim}(H)\) corresponds to \(n = 0\), i.e.\
\(\left.\e_H\right|_{\mathrm{Prim}(H)} = \mathrm{id}_{\mathrm{Prim}(H)}\), and the inclusion follows by idempotence. We stress that equality does not hold for \(\hck\), as illustrated by the following simple example:
\[
\eck\big( \Forest{[\beta[\alpha]]}\big) = \Forest{[\beta[\alpha]]} - \tfrac 12 \Forest{[\alpha]}\Forest{[\beta]} \not\in \cP,
\]
of which \(\pi(\Forest{[\alpha]}\Forest{[\beta]})\) is the symmetrisation. Statement (4) is the celebrated Milnor-Moore theorem \cite{MM65}. To see (5), recall the notion of \emph{infinitesimal character} \cite[Proposition 22]{Man06}, in our case taking values in \(H\) itself, i.e.\ elements \(c \in \cL(H,H)\) s.t.\ 
\[
c(xy) = 1_H(x) c(y) + c(x)1_H(y).
\]
When \(H\) is commutative, \(\exp_*\) and \(\log_*\) define bijections, inverse to each other, between the group of \(H\)-valued characters, i.e.\ algebra morphisms \(H \to H\), and the Lie algebra of \(H\)-valued infinitesimal characters.
Then if \(x,y \in H_+\) it is immediate that \(\e_H(xy)=\log_*(\mathrm{id})(xy) = 0\). The splitting of the short exact
sequence amounts to saying that, for \(h \in H_+^{\times 2}\), \(\e_H(h) = h + k\) with \(k \in H_+^{\times 2}\), which is again evident from \eqref{eq:e}. In the case of \(H = \hck\), the exact sequence is
\[
\begin{tikzcd}
	0 \arrow[r] &\cF \arrow[r] &\cH_+ \arrow[r] &\cT \arrow[r] \arrow[l,bend left,"\eck"] &0 
\end{tikzcd}.
\]
Assertion (5) is known as Leray's theorem, see again \cite[Theorem 7.5]{MM65}, and can be viewed as of a piece with the
statement in the cocommutative case, in the sense that the universal enveloping algebra of the trivial Lie algebra is
the symmetric algebra.

We now proceed with \(H = \hck\), and continue to denote \(\e\) the corresponding Eulerian idempotent. We denote \(\cE
\coloneqq \e(\hck)\) and, as usual, view it as a functor \(\cat{Vec} \to \cat{Vec}\).
We also view \(\e\) as a natural transformation \(\cH\Rightarrow\cE\), which is possible since \(\e\) is defined via the natural operations of \(\hck\), and \Cref{prop:functor}. \Cref{prop:e} implies that the inclusion \(\cE \hookrightarrow \hck\) induces a natural isomorphism \(\hck \cong \bigodot(\cE)\).

In this section we will show how to define a natural isomorphism from \(\hck \cong \bigshuffle (\cP)\). This will be achieved by combining both maps \(\pi\) and \(\e\). Before doing this, however, it is instructive to see that considering the two descriptions of \(\hck\) separately --- one as a tensor algebra and one as a symmetric algebra --- does not yield the desired result.
\begin{proposition}[Two failed attempts at a Hopf isomorphism]\label{prop:failed}\ \\[-2ex]
\begin{enumerate}
	\item \(\top^{-1} \colon \hck \to \bigshuffle(\cP)\) is a coalgebra isomorphism but not an algebra morphism.
	\item The unique algebra morphism \(\varphi \colon \hck \to \bigshuffle(\cP)\) with \(\varphi|_\cE = \top^{-1}|_\cE\) is an algebra isomorphism but not in general a coalgebra morphism.
\end{enumerate}
\begin{proof}
That \(\top^{-1}\) is a coalgebra isomorphism was stated in \Cref{thm:T}. To check that it is not an algebra isomorphism, the following suffices:
\begin{align*}
\top^{-1}(\Forest{[\alpha]} \cdot \Forest{[\beta]}) = \big(\Forest{[\alpha]} \Forest{[\beta]} - \Forest{[\beta[\alpha]]} - \Forest{[\alpha[\beta]]}\big) + \Forest{[\alpha]} \otimes \Forest{[\beta]} +  \Forest{[\beta]} \otimes \Forest{[\alpha]} \neq  \Forest{[\alpha]} \otimes \Forest{[\beta]} +  \Forest{[\beta]} \otimes \Forest{[\alpha]} = \top^{-1}(\Forest{[\alpha]}) \shuffle \top^{-1}(\Forest{[\beta]})
\end{align*}

We now consider the map \(\varphi\). It is clear from free commutativity of \(\cH\) over \(\cE\) that such a map exists and is an algebra morphism. Here and below we use the isomorphism in (1) to identify the underlying space $\hck$ with $\bigotimes(P)$. The \say{shuffle} of elements of $\hck$ is denoted (slight abusing the notation) $\shuffle \coloneqq \top \circ \shuffle \circ \top^{-1}$, which by \Cref{cor:fprod} is comprises the top level term in the forest product
\begin{align*}
(p_1 \tee \ldots \tee p_m) \shuffle (p_{m+1} \tee \ldots \tee p_{m+n}) &= \sum_{\sigma \in \mathrm{Sh}(m,n)}p_{\sigma^{-1}(1)} \tee \ldots \tee p_{\sigma^{-1}(m+n)} \\
&= \pi_{m+n}\big( (p_1 \tee \ldots \tee p_m) \cdot (p_{m+1} \tee \ldots \tee p_{m+n}) \big),
\end{align*}
although we note that this will not be used in this proof (cf. \cite[Lemma 11.1]{Foi2002}).
We are now interested in showing that the algebra morphism \((\hck, \cdot) \to (\hck, \shuffle)\), which we continue to call \(\varphi\), is bijective.
To do this, we observe that it admits the expression \(\varphi = \smash{\exp_*^{\shuffle}} \circ \log_*(\mathrm{id})\),
where the superscript indicates that the convolution product is taken w.r.t.\ the product \(\shuffle\) instead of the
forest product whilst in both cases the coproduct \(\dck\) is the same.
This is because \(\varphi\) maps
\begin{align*}
  \hck \ni h &=\exp_* \circ \log_*(\mathrm{id})(h) \\
             &=\sum_{n = 0}^\infty \frac{1}{n!} \e(h^{(1)}) \cdots  \e(h^{(n)}) \\
             &\mapsto  \sum_{n = 0}^\infty \frac{1}{n!} \e(h^{(1)}) \shuffle \cdots \shuffle \e(h^{(n)}) \\
             &=\exp_*^{\shuffle} \circ \log_*(\mathrm{id})(h).
\end{align*}
From this it follows that for \(h\in \cE\), denoting \(\e_{\shuffle} = \log^{\shuffle}_*(\mathrm{id})\) the Eulerian idempotent of the Hopf algebra \(\bigshuffle(P)\)
\begin{align*}
\e_{\shuffle}(h) &= \sum_{n = 1}^\infty \frac{(-1)^{n-1}}{n} h^{(1)} \shuffle \cdots \shuffle h^{(n)} \\
&= \varphi \left( \sum_{n = 1}^\infty \frac{(-1)^{n-1}}{n}h^{(1)}  \cdots  h^{(n)} \right)\\
&= \varphi \circ \eck(h) \\
&= \varphi(h) \\
&= h
\end{align*}
which by idempotence implies \(\cE \subseteq \e_{\shuffle}(\cH)\), and running through the same argument with products reversed yields the other inclusion, i.e.\ \(\cH\) is free commutative over the same \(\cE\) w.r.t.\ either product. From this it follows immediately that \(\varphi\) is bijective with inverse \(\exp_* \circ \smash{\log_*^{\shuffle}}(\mathrm{id})\).

In order to find a counterexample to the coalgebra morphism property, we must consider degree-3 trees with decorations. We show that
\[
\dck \varphi \bigg( \Forest{[\c[\b[\a]]]}\bigg) -\varphi^{\otimes 2}  \dck  \bigg(\Forest{[\c[\b[\a]]]}\bigg) \neq 0
\]
We refer the reader to \Cref{app:calc} for the full calculation. This completes the proof.
\end{proof}
\end{proposition}
	
The first map in \Cref{prop:failed} did not use Eulerian idempotents, and indeed did not require commutativity, while the second failed to leverage cofreeness. In the following theorem we put both structures to use. We note that the use of \say{\(\Log\)} and \say{\(\Exp\)} is motivated by Hoffman's isomorphism, see \Cref{subsec:hoffman} below. Recall from the above proof that here and below we implicitly identify $\hck$ with $\bigotimes (\cP)$ via the isomorphism $\top$, which results in the identifications $\pi_1 = \pi$, $\otimes = \top$, etc. $\e$ will always refer to $\e_\mathrm{CK}$.
	
\begin{theorem}[Natural isomorphism \(\hck \cong \bigshuffle(\cP)\)] \ \\ \label{thm:iso}
	\noindent There exists a unique isomorphism \(\Log \colon \hck \to \bigshuffle(\cP)\) with the property that 
	\[
	\pi \circ \Log = \pi \circ \e
	\]
	and it is given by
	\[
	\Log = \sum_{n \geq 0} ( \pi \circ \e )^{\otimes n} \circ \dckRed^{(n-1)} \,.
	\]
	Its inverse \(\Exp \colon \bigshuffle(\cP) \to \hck\) is given on $\bigshuffle(\cP)_+$ by the following combinatorial expression, involving a summation over increasing sequences of integers and their compositions
	\[
	\sum_{\substack{1 \leq n^0 < \ldots < n^k \\ m^1_1 + \ldots + m^1_{n^{\scaleto{0}{3pt}}} = n^1 \\ \ldots \\ m^k_1 + \ldots + m^k_{n^{\scaleto{k-1}{3pt}}} = n^k \\ k \geq 0, \ m^i_j \geq 1}}  (-1)^k [(\pi \circ \e \circ \pi_{m_1^1}) \otimes \cdots \otimes (\pi \circ \e \circ \pi_{m_{n^{\scaleto{0}{3pt}}}^1})] \circ \cdots \circ [(\pi \circ \e \circ \pi_{m_1^k}) \otimes \cdots \otimes (\pi \circ \e \circ \pi_{m_{n^{\scaleto{k-1}{3pt}}}^k})]
	\]
	(where the summand corresponding to $k = 0$ is $\mathrm{id}$) or equivalently by the recursion
	\[
	\Exp^m_n = -\sum_{m \leq k < n} \Exp^m_k \circ \Log^k_n,\quad m < n
	\]
	and \(\Exp^n_n = \mathrm{id} \), where \(\Exp^m_n\coloneqq \pi_m\circ\Exp\circ\pi_n\), and similarly for \(\Log\), for $1 \leq m \leq n$. \(\Log\) and \(\Exp\) are natural isomorphisms between the functors \(\hck, \bigshuffle \circ \cP \colon \cat{Vec} \to \bnohat\) restricting to the identity on \(\cP\).
	\begin{proof}
		By cofreeness \eqref{eq:cofree}, there is a unique coalgebra homomorphism \(\Log\) with a specified projection \(\pi
		\circ \Log\) (which in this case equals to \(\pi \circ \e\)). By \Cref{coalgebra_prop} it follows that such a map will have the following expression: for \(p_1,\ldots,p_n \in \cP\)
		\[
		\Log(p_1 \otimes \ldots \otimes p_n) = \sum_{\substack{n_1 + \ldots + n_k = n \\ n_1,\ldots,n_k \geq 1}} (\pi \circ \e)(p_1 \otimes \ldots \otimes p_{n_1}) \otimes \cdots \otimes (\pi \circ \e)(p_{n_1 + \ldots + n_{k-1}+1} \otimes \ldots \otimes p_n)
		\]
		which by \Cref{thm:T} equals the second expression in the statement; the conclusion then follows by linear extension. We must show that this coalgebra morphism is an algebra morphism: denoting \(\times = \cdot\) the forest product for better clarity, this is the case if and only if
		\[
		\Log \circ (\argument \times \argument) =  (\argument \shuffle \argument) \circ \Log^{\otimes 2} \colon \hck^{\otimes 2} \to \bigshuffle(\cP)
		\]
		Since \(\hck\) and \(\bigshuffle(\cP)\) are bialgebras and by the above, both these maps are coalgebra morphisms, and thus again by cofreeness they coincide if and only if
		\[
		\pi \circ \e \circ (\argument \times \argument) = \pi \circ \Log \circ (\argument \times \argument) =  \pi \circ
		(\argument \shuffle \argument) \circ \Log^{\otimes 2} = 0,
		\]
		on $\cH_+ \times \cH_+$ since \(\pi\) vanishes on \(\bigshuffle (\cP)_+^{\shuffle 2}\). But by point 5.\ of \Cref{prop:e}, \(\e_{\mathrm{CK}}\)
		maps \(\cH_+^{\times 2}\) to zero, concluding the proof that \(\Log\) is a bialgebra, and hence Hopf algebra morphism.
		\Cref{coalgebra_prop} already implies that \(\Log\) is invertible, since \(\Log|_\cP = (\pi \circ \e)|_\cP = \mathrm{id}_\cP\), a consequence of point 3.\ of \Cref{prop:e}; this also implies that also its inverse \(\Exp\) restricts to the identity on \(\cP\).
		
		In order to prove the formula for \(\Exp\), it is convenient to think of a coalgebra map between cofree coalgebras as a
		triangular matrix: again by \Cref{coalgebra_prop}, setting \(\Log^m_n \coloneqq \pi_m \circ \Log \circ \pi_n \colon \hck \to \bigshuffle(\cP)\) for $1 \leq m \leq n$, calling $\Log_+ \coloneqq \Log|_{\cH_+}$ we have
		\begin{align}
			&\Log_+ = \underbrace{\sum_{1 \leq m < n} \Log^m_n}_{\eqqcolon L} + \underbrace{\sum_{n \geq 1} \Log^n_n}_{= \mathrm{id}} \nonumber \\
			&\text{with}\quad \Log^m_n = \sum_{\substack{n_1 + \ldots + n_m = n \\ n_1,\ldots, n_m \geq 1}} (\pi \circ \e \circ \pi_{n_1}) \otimes \cdots \otimes (\pi \circ \e \circ \pi_{n_m}) \label{eq:lognm}
		\end{align}
		$\Log_+ = L + \mathrm{id}$ admits the explicit inverse
		\begin{align*}
			(L + \mathrm{id}) \circ \sum_{k \geq 0}(-1)^k L^{\circ k} = \mathrm{id} = \sum_{k \geq 0}(-1)^k L^{\circ k} \circ (L + \mathrm{id})
		\end{align*}
		with $L^{\circ 0} = \mathrm{id}$ and the sum locally finite, since \(L\) is nilpotent of degree \(N\) on \(\bigoplus_{n = 1}^N \cP^{\otimes n}\). This implies
		\begin{align*}
			\Exp_+ &=\Log^{-1}_+ \\
			&= \sum_{k \geq 0}(-1)^k L^{\circ k}\\
			&=\sum_{k \geq 0}\sum_{\substack{1 \leq m^1 < n^1 \\ \ldots \\ 1 \leq m^k < n^k}} (-1)^k \Log^{m^1}_{n^1} \circ \cdots \circ \Log^{m^k}_{n^k} \\
			&=\sum_{k \geq 0} (-1)^k \sum_{1 \leq n^0 < \ldots < n^k} \Log^{n^0}_{n^1} \circ \cdots \circ \Log^{n^{k-1}}_{n^k}
		\end{align*}
		The non-recursive expression is now obtained by substituting in \eqref{eq:lognm} for each term $\Log^{n^{i-1}}_{n^i}$. As for the recursive expression, we have $\Exp^n_n = \mathrm{id}$ and for $1 \leq m < n$
		\begin{align*}
			\Exp^m_n &= \sum_{k \geq 1} (-1)^k \sum_{1 \leq n^0 < \ldots < n^k} \Log^{n^0}_{n^1} \circ \cdots \circ \Log^{n^{k-1}}_{n^k} \\
			&= - \Log^m_n + \sum_{k \geq 2} (-1)^k\sum_{\substack{n^0 < \ldots < n^k \\ m \eqqcolon n^0, \ n \eqqcolon n^k}} \Log^{n^0}_{n^1} \circ \cdots \circ \Log^{n^{k-1}}_{n^k} \\
			&= - \Log^m_n - \sum_{\substack{m < h < n \\ m \eqqcolon n^0,\ h \eqqcolon n^{k-1}}} \sum_{k \geq 2} (-1)^{k-1} \sum_{n^0 < \ldots < n^{k-1}} \Log^{n^0}_{n^1} \circ \cdots \circ \Log^{n^{k-2}}_{n^{k-1}} \circ \Log^h_n \\
			&= - \Log^m_n - \sum_{m < h < n} \Exp^m_h \circ \Log^{h}_{n}  
		\end{align*}
		from which we conclude by expressing $\Log^m_n$ as $\sum_{m = h < n} \Exp^m_h \circ \Log^h_m$. Naturality follows again from \Cref{prop:functor} and the fact that all operations involved in the expressions of \(\Log\) and \(\Exp\) are natural.
	\end{proof}
\end{theorem}

\begin{remark}\label{rem:general}
To keep the presentation concise and concrete we have stated the above theorem in the case which is relevant to the
Itô-Stratonovich transformation of branched rough paths, \Cref{cor:itoStrat} below. However, both the statement and
the proof carry over word for word to the case of arbitrary functors \(\cat{Vec} \to \bnohat\), simply by replacing \(\pi\) with the cofreeness projection and \(\e\) with the Eulerian idempotent for the commutative \(\binf\)-algebra in question. In particular, while in \Cref{sec:ito} we rely on the dual Hopf algebra being pre-Lie, this is not needed for this result, which also applies, for instance, to the Munthe-Kaas Wright Hopf algebra of planar rooted trees \cite{MR2407032}, cf.\ \cite[Proposition 4.3]{ebrahimifard2023primitive} and references therein.
\end{remark}

\begin{remark}\label{rem:FP}
Since the first version of this paper appeared as a preprint, in \cite{FP24} Foissy and Patras have expanded the study of natural isomorphisms between commutative $\binf$-algebras and their shuffle algebra over primitives. They observe that, given any commutative $\binf$-algebra and any Lie idempotent $\ell^*$ of the dual Hopf algebra, $\pi \circ \ell$ yields a natural isomorphism $\Log_\ell$ similar to the one of \Cref{thm:iso}. By taking $\ell = \e$ we recover our result, but different idempotents define different isomorphisms; for example taking $\ell$ to be the Dynkin idempotent, $\mathrm{Dyn} \coloneqq \EuScript{D}/|\, \cdot \,|$ where $\EuScript{D}$ is the Dynkin operator (see \eqref{eq:Dynkin} below), direct calculation yields that, for example, for primitives $p$ and $q$ one obtains
\[
\Log_{\mathrm{Dyn}}(p \tee q) = p \otimes q - \frac{|q|}{|p| + |q|} \pi(pq) \neq p \otimes q - \frac 12 \pi(pq) = \Log_{\e}(p \tee q) .
\]
In particular, in the quasi-shuffle case $\Log_{\mathrm{Dyn}}$ does not coincide with the Hoffman isomorphism, while $\Log = \Log_{\e}$ does (cf.\ \Cref{subsec:hoffman} below). Notice moreover that the grading appears expicitly in $\Log_{\mathrm{Dyn}}$ but not in $\Log$, cf.\ \cite[Section 8]{FP24} which focuses on the non-graded case and introduces a more general form of naturality. We believe these ideas to be relevant to the question of uniqueness, see \Cref{rem:uniq} below; pushing such considerations further requires a separate effort.
\end{remark}

\begin{example}
  We compute \(\Exp\) up to level 3 in the primitiveness grading (which yields more compact formulae than the forest basis).
  The expression in the forest basis can be obtained via \Cref{expl:P3}.
  \begin{align*}
    \Exp\left(p\right) &= p,\qquad\Exp\left(p\otimes q\right) = p\tee q +\frac12\pi(pq) \\
    \Exp\left( p\otimes q\otimes r \right) &=\begin{multlined}[t] p\tee q\tee r+\frac12\left( \pi(pq)\tee r + p\tee\pi(qr) \right) +
    \frac16\pi(pqr)\\
  +\frac14\left( \pi((p\tee q)r) + \pi(p(q\tee r)) - \pi((q\tee p)r) - \pi(p(r\tee q)) \right) \end{multlined}
  \end{align*}
  We note that everything on the first line is already present in Hoffman’s isomorphism for the quasi-shuffle case (see below
  \Cref{subsec:hoffman}), while terms in the second line are specific to the case of Connes-Kreimer.
The following is an example of $\mathrm{Log}$ at level 4 in the undecorated case, in which we use the identity $\pi(\Forest{[]}\Forest{[]}\Forest{[[]]}) = \pi(\Forest{[[]]}\Forest{[[]]})$.
\begin{align*}
      \Log \, \Forest{[[[[]]]]} &=\begin{multlined}[t] \bullet \otimes \bullet \otimes \bullet \otimes \bullet - \frac 12 \big[ \bullet \otimes \bullet \otimes \pi(\Forest{[]}\Forest{[]}) + \bullet \otimes \pi(\Forest{[]}\Forest{[]}) \otimes \bullet + \pi(\Forest{[]}\Forest{[]}) \otimes \bullet \otimes \bullet \big]\\
      	+ \frac 13 \bullet \otimes \pi(\Forest{[]}\Forest{[]}\Forest{[]}) + \frac 14 \pi(\Forest{[]}\Forest{[]}) \otimes \pi(\Forest{[]}\Forest{[]}) + \frac 13 \pi(\Forest{[]}\Forest{[]}\Forest{[]}) \otimes \bullet - \frac 14 \pi(\Forest{[]}\Forest{[]}\Forest{[]}\Forest{[]}) + \frac 12 \pi(\Forest{[[]]}\Forest{[[]]}).\qedhere\end{multlined}
\end{align*}
\end{example}
	
We now discuss the consequences that this has for rough paths. Recall that a \(\cQ\)-valued geometric \(\rho\)-rough path (of inhomogeneous regularity) is defined in complete analogy with \Cref{def:brp}, with the only difference that the Hopf algebra \(\hck\) is replaced with \(\bigshuffle(\cP)\). An important caveat is that the grading on \(\bigshuffle(\cP)\) used in the regularity requirement is the inhomogeneous one, which takes into account the grading on \(\cP\) inherited from \(\hck\)
\begin{equation}\label{eq:weight}
|p_1 \otimes \cdots \otimes p_n| = |p_1| + \ldots + |p_n|
\end{equation}
and \(\rho\) refers to terms of the worst regularity, \(X|_V\). This also affects how the tensor algebra (along with Davie expansions, etc.) is truncated.
In particular from the above examples we clearly see that \(\Exp\) (and \(\Log\)) are not graded maps when using the
standard grading in \(\bigshuffle(\cP)\) since it decreases primitiveness, but it is under \eqref{eq:weight}.

Recall that the graded dual bialgebra to \(\bigshuffle(\cP)\) is \(\bigotimes (\cQ) \coloneqq (\bigotimes (\cQ), \otimes, \Delta_\shuffle)\), where \(\Delta_\shuffle\) denotes \say{unshuffling} (see, for instance, \cite[\S 1.5]{Reu93}), and the dual pairing is the tensor product of the dual pairings \(\cP \otimes \cQ \to \bbR\) inherited from \(\hck \otimes \hgl \to \bbR\) of \Cref{thm:dual}. We can now state the main application of \Cref{thm:iso}.
\begin{corollary}[It\^o-Stratonovich correction formula for RDEs driven by branched rough paths]\label{cor:itoStrat}
Let \(\bX\) be a \(\rho\)-branched rough path controlled by \(\omega\). Then \(\overline \bX \coloneqq \Exp^*(\bX)\) defines a \(\cQ\)-valued geometric \(\rho\)-rough path controlled by \(\omega\) and the \(\bX\)-driven RDE \eqref{eq:rde} is equivalent to the \(\overline \bX\)-driven one
\[
  \dif Y = \widehat{F} \circ \e_{\mathrm{GL}}|_{\cQ}(Y) \dif \overline{\bX} = \sum_{\scrf \in \mathscr{H}} \varsigma(\scrf)^{-1}\sum_{l\in\ell(\scrf)}\langle \widehat{F}(Y),
  \e_\mathrm{GL} \circ \pi^* (\scrf_l)
\rangle\,\dif \overline{\bX}{}^{\pi(\scrf_l)}
\]
where \(\widehat{F}\) is taken as in \Cref{thm:davie.drift}.
\end{corollary}	
\begin{proof}
	The first statement is an immediate consequence of \Cref{thm:iso} and the two definitions of rough path, together with the fact that \(\Exp\) preserves the grading in \eqref{eq:weight}.
  The \(\overline\bX\)-driven equation in the statement has Davie expansion
	\begin{align*}
		Y_{s,t} &\approx \langle \Log (\widehat F(Y_s)), \Exp^*(\bX_{s,t}) \rangle \approx  \langle \Exp \circ \Log (\widehat F(Y_s)), \bX_{s,t} \rangle \approx \langle \widehat F(Y_s), \bX_{s,t} \rangle
	\end{align*}
\end{proof}

\begin{example}[The It\^o-Stratonovich formula for \(\p = 3\)]
  Continuing with \Cref{ex:Davie}, we rewrite the RDE in geometric terms using \Cref{cor:itoStrat}.
  One directly obtains:
  \begin{align*}
	\dif Y_t
	={}&F_{\Forest{[\gamma]}}(Y_t)\,\dif\overline{\bX}{}_t^{\Forest{[\gamma]}}-\frac12\widehat{F}_{\Forest{[\beta[\alpha]]}}(Y_t)
	\dif\overline{\bX}{}_t^{\pi(\Forest{[\alpha]}\Forest{[\beta]})}\\ &+ \frac 16 \bigg( \frac 12 \widehat F_{\Forest{[\gamma[\alpha][\beta]]}} + 2 \widehat F_{\Forest{[\gamma[\beta[\alpha]]]}} \bigg)(Y_t) \dif \overline\bX{}^{\pi(\Forest{[\alpha]}\Forest{[\beta]}\Forest{[\gamma]})}
	\\ &+\frac 12 \bigg(
	\widehat{F}_{\Forest{[\alpha[\beta][\gamma]]}}
	+ \widehat{F}_{\Forest{[\alpha[\beta[\gamma]]]}}
	+ \widehat{F}_{\Forest{[\gamma[\alpha[\beta]]]}}
	\bigg)(Y_t) \dif \overline{\bX}{}^{\pi(\Forest{[\alpha]}\Forest{[\gamma[\beta]]})}_t
\end{align*}
where the sum over labels has been omitted in order to declutter the notation, and the expression has been simplified by taking (anti)symmetric relations into account in the contractions. Note that, truncating at degree 2, since
\[
-\frac 12 \widehat{F}_{\Forest{[\beta[\alpha]]}} = \frac 12 F_{\Forest{[\alpha]}\Forest{[\beta]}} - \frac 12 F_{\Forest{[\beta]}} \rhd F_{\Forest{[\alpha]}}
\]
	we recover the well-known Itô--Stratonovich correction for SDEs driven by semimartingales (cf.\ \eqref{eq:itostrat}, keeping in mind that our drift terms $F_{\Forest{[\gamma]}\Forest{[\gamma]}}$ carry a factor of $\frac 12$). We refer again to \Cref{ex:Davie} for the concrete expression of \(\widehat{F}\) in terms of the original vector fields \(F\).
\end{example}

\cite{K80} first observed that solutions to Stratonovich SDEs with nilpotent vector fields could be written as ODEs with vector fields which are linear functions of the log-signature, extending the Doss-Sussman representation valid in the one-dimensional case. The same principle was used in \cite{CG95} for numerical purposes, dropping the nilpotency assumption, and see \cite{boutaib2013dimensionfree} for the geometric rough path case. Here we limit ourselves to writing down a similar (formal) expansion for the branched case, which is proved in much the same way as it is in the geometric case. See also \cite{kern2023flow} for the case of manifold-valued RDEs (with a different Hopf algebra).

	As a premise to the next result, we note that in analogy to the geometric case, given the signature \Cref{eq:signature}, it is natural to define the \emph{log-signature} by $\log_\star \EuScript{S}(\bX)$, where 
	\[
	\log_\star(h) = \sum_{n \geq 1} \frac{(-1)^{n-1}}{n} (h-1)^{\star n}
	\]
	Using the character property of \(\EuScript{S}(\bX)\), we can express log-signature coordinates in terms of signature coordinates using the Eulerian idempotent (which \emph{is} linear):
	\[
		\begin{split}
			\langle h, \log_\star \EuScript{S}(\bX) \rangle &= \langle \e_{\mathrm{CK}}(h) , \EuScript{S}(\bX) \rangle\\
			&= \sum_{n \geq 1} \frac{(-1)^{n-1}}{n} \langle h^{(1)} \cdots h^{(n)}, \EuScript{S}(\bX)\rangle \\
			&= \sum_{n \geq 1} \frac{(-1)^{n-1}}{n} \langle h^{(1)},
			\EuScript{S}(\bX) \rangle \cdots \langle h^{(n)},
			\EuScript{S}(\bX)\rangle.
		\end{split}
	\]

  \begin{corollary}[log-ODE method for branched rough paths]\label{cor:logode}
		Let $Y$ be as in \eqref{eq:rde}. For \(\varphi\in C^{\infty}(W)\) and $s < t$ we have the formal identity $\varphi(Y_t) = \varphi(Z_1)$, where
		\[
      \dot Z = \langle \widehat F(Z), \log_\star \EuScript{S}(\bX)_{s,t} \rangle = \sum_{\scrf \in \scrT} \varsigma(\scrf)^{-1} \sum_{l \in \ell(\scrf)} \widehat F_{\scrf_l}(Z) \langle \e(\scrf_l), \EuScript{S}(\bX)_{s,t} \rangle  , \quad Z_0 = Y_s
		\]
	\end{corollary}
		\begin{proof}
		For $\varphi \in C^\infty(W)$
		\begin{align*}
		\varphi(Y_t) &= \langle \mathbf{F}\varphi(Y_s), \EuScript{S}(\bX)_{s,t} \rangle \\
		&= \langle \mathbf{F}\varphi(Y_s), \EuScript{S}(\bX)_{s,t} \rangle \\
		&= \sum_{n = 0}^\infty \frac{1}{n!} \langle \mathbf{F}\varphi(Y_s), (\log_\star\EuScript{S}(\bX)_{s,t})^{\star n} \rangle \\
		&= \sum_{n = 0}^\infty \frac{1}{n!}\langle \widehat{F}, \log_\star \EuScript{S}(\bX)_{s,t} \rangle^{\circledast n}\varphi(Y_s) \\
		&= \varphi(Z_1) \, .
		\end{align*}
		since the last expression is the formal Taylor expansion of the ODE in the statement. We have used the fact $\mathbf{F}$ is an algebra morphism and that $\log_\star \EuScript{S}(\bX)_{s,t} \in \cT = \mathrm{Prim}(\hgl)$ (a consequence of the fact that the signature is grouplike). In the coordinate representation, we are summing over trees because $\e$ vanishes on proper forests.
		\end{proof}
	For example, truncating the expansion at level-2 this reads
	\begin{align*}
	\dot Z &= F_\gamma(Z) X^\gamma_{s,t} + (F_\beta \rhd F_\alpha)(Z) \langle \Forest{[\a [\b]]} - \tfrac 12 \Forest{[\a]} \Forest{[\b]} , \bX_{s,t} \rangle + \tfrac 12 F_{\Forest{[\a]}\Forest{[\b]}}(Z) \bX_{s,t}^{\pi({\Forest{[\a]}\Forest{[\b]}})} + \ldots \\
	&= F_\gamma(Z) X^\gamma_{s,t} + \tfrac 12(F_\beta \rhd F_\alpha)(Z)  \langle \Forest{[\a [\b]]} - \Forest{[\b [\a]]} , \overline \bX_{s,t} \rangle + \tfrac 12 [ F_{\Forest{[\a]}\Forest{[\b]}} - F_\beta \rhd F_\alpha](Z) \overline \bX{}^{\pi({\Forest{[\a]}\Forest{[\b]}})}_{s,t} + \ldots
\end{align*}
with the second identity written in terms of the associated Stratonovich rough path, i.e., after applying \Cref{cor:itoStrat} (see also \cite{KMM23} for similar results). Note that the expansions above can be truncated so that only derivatives of $F$ up to order $\lfloor \rho \rfloor$ appear so the smoothness assumption can be relaxed, cf.\ \Cref{rk:RDE}.

	\subsection{The quasi-geometric case: Hoffman's exponential}\label{subsec:hoffman}
In \cite{Hof00}, Hoffman defined a Hopf isomorphism between the quasi shuffle and shuffle bialgebras. In this subsection, we show how this isomorphism can be viewed as a special case of \Cref{thm:iso}. In the context of this paper, and of the notation introduced in \Cref{subsec:quasi}, we call it \(\qExp \colon \bigshuffle (\bigodot(U)) \to \bigqshuffle (U)\) and its inverse \(\qLog \colon \bigqshuffle (U) \to \bigshuffle (\bigodot(U))\). We recall their definition: for \(s_1,\ldots,s_n \in \bigodot(U)\)
\begin{equation}\label{eq:hoffman}
	\begin{split}
	\qExp(s_1 \otimes \cdots \otimes s_n) &\coloneqq \sum_{\substack{n_1 + \ldots + n_k = n \\ n_1 , \ldots , n_k \geq 1}} \frac{1}{n_1! \cdots n_k!} (s_1 \odot \cdots \odot s_{n_1}) \otimes \cdots \otimes (s_{n_1 + \ldots + n_{k-1}+1} \odot \cdots \odot s_n) \\
	\qLog(s_1 \otimes \cdots \otimes s_n) &\coloneqq \sum_{\substack{n_1 + \ldots + n_k = n \\ n_1 , \ldots , n_k \geq 1}} \frac{(-1)^{n-k}}{n_1 \cdots n_k} (s_1 \odot \cdots \odot s_{n_1}) \otimes \cdots \otimes (s_{n_1 + \ldots + n_{k-1}+1} \odot \cdots \odot s_n)
	\end{split}
\end{equation}
The following theorem establishes these isomorphisms as particular cases of \Cref{thm:iso} (defined w.r.t.\ the functor \(\bigqshuffle \colon \cat{Vec} \to \bnohat\), see \Cref{rem:general}); furthermore it shows how our and Hoffman's isomorphisms are related via the arborification quotient maps of \Cref{prop:arb}, and are furthermore related to the recently introduced arborified exponential introduced in \cite{BCE20}, denoted here by \(\widetilde{\mathfrak{a}}_\odot\). It is defined in complete analogy with \(\widetilde{\mathfrak{a}}\) (see \Cref{prop:arb}), the only difference being that forests in the source are decorated with \(\bigodot(U)\), not \(U\), and that decorations are multiplied associatively when taking the shuffle product; it has the alternative right inverse given by \(\gamma_1 \odot \cdots \odot \gamma_n \mapsto \Forest{[\gamma_1 \odot \cdot\cdot\cdot \odot \gamma_n]}\) defined without invoking \(\pi\). The existence of a unique and explicit Hopf algebra morphism \(\mathfrak a \qExp\) making the lower parallelogram commute is proved in \cite[Theorem 2]{BCE20}.
\begin{theorem}
Hoffman's isomorphisms are given by the same formulae of \Cref{thm:iso}, as specified by \Cref{rem:general}, and the following diagram of Hopf morphisms
\[
	\begin{tikzcd}[ampersand replacement=\amp]
	\amp \bigshuffle(\cP(U))  \arrow{rrr}{\Exp} \amp \amp \amp \hck(U) \arrow[dd,twoheadrightarrow,"\widetilde{\mathfrak{a}}"] \\
\hck(\bigodot (U)) \arrow[dr,"\mathfrak a",twoheadrightarrow] \arrow[ur] \arrow[rrr,"\mathfrak a \qExp",near end] \amp \amp \amp \hck(\bigodot(U)) \arrow[ur,hookleftarrow,"\hck(j)",swap] \arrow[dr,"\widetilde {\mathfrak a}_\odot",twoheadrightarrow] \\
	\amp \bigshuffle (\bigodot (U)) \arrow{rrr}{\qExp} \arrow[uu,hookrightarrow,crossing over,near start,"\bigshuffle(i)",swap] \amp \amp \amp \bigqshuffle (U)
\end{tikzcd}
\]
where \(i \coloneqq \bigodot (U) \hookrightarrow \cP(U)\) and \(j \coloneqq U \hookrightarrow \bigodot(U)\) are the natural inclusions (and the unlabelled map is the composition \(\mathfrak a \circ \bigshuffle(i)\)), commutes.

\begin{proof}
We begin with the first statement, which we verify in the case of the logarithm; the analogoues statement for the exponential follows by uniqueness of inverses. By cofreeness, it suffices to check that \(\pi_1 \circ \qLog = \pi_1 \circ \e_{\qshuffle}\), where \(\pi_1 \colon \bigqshuffle (U) \twoheadrightarrow \bigodot(U)\) and \(\e_{\qshuffle}\) is the Eulerian idempotent for the quasi shuffle algebra. We compute, using \Cref{def:euler} and \eqref{eq:hoffman}
\begin{align*}
	\pi_1 \circ \e_{\qshuffle} (s_1 \otimes \cdots \otimes s_n) &= \sum_{m \geq 1}\frac{(-1)^{m-1}}{m} \pi_1 \circ \qshuffle^{(m-1)} \circ \widetilde\Delta_\otimes^{(m-1)} (s_1 \otimes \cdots \otimes s_n) \\
	&= \frac{(-1)^{n-1}}{n}(s_1 \odot \cdots \odot s_n) \\
	&= \pi_1 \circ \qLog (s_1 \otimes \cdots \otimes s_n)
\end{align*}
since the only term to survive application of \(\pi_1\) appears once as a summand in the \((n-1)\)-fold deconcatenation. Although redundant, it is informative to see why Hoffman's exponential has a much simpler closed-form expression than \(\Exp\) for general \(\binf\)-algebras. We check that \(\pi \circ \qExp\) satisfies the same recursion as \(\Exp\) (which is already implied by \Cref{rem:general}):
\begin{align*}
&- \sum_{1 \leq k < n} \qExp^1_k \circ \qLog^k_n (s_1 \otimes \cdots \otimes s_n) \\
={} & - \sum_{1 \leq k < n} \qExp^1_k \circ \sum_{\substack{n_1 + \ldots + n_k = n \\ n_1, \ldots, n_k \geq 1}} \frac{(-1)^{n-k}}{n_1 \cdots n_k} (s_1 \odot \cdots \odot s_{n_1}) \otimes \cdots \otimes (s_{n_1 + \ldots + n_{k-1}+1} \odot \cdots \odot s_{n_k}) \\
={} &\sum_{\substack{1 \leq k < n \\ n_1 + \ldots + n_k = n \\ n_1, \ldots, n_k \geq 1}} \frac{(-1)^{n-k}}{n_1 \cdots n_k} s_1 \odot \cdots \odot s_n
\end{align*}
and the assertion follows from the fact that the coefficient is equal to \(1/n!\), as can be verified by writing out \(t = \exp \circ \log (t+1) - 1\) in terms of power series. From this we see that the simplifications for the closed-form expression for the exponential in the quasi-geometric case are due to the fact that projection onto primitives behaves associatively, i.e.\ \(\pi_1(x\qshuffle y) = \pi_1(\pi_1(x)\qshuffle \pi_1(y))\), which is not true on \(\hck\).

We must show commutativity of two rectangles, two parallelograms, and a rectangle. We begin by considering the alternative rectangle
\begin{equation}\label{eq:altDiag}
\begin{tikzcd}[row sep = large, column sep = large]
\bigshuffle(\cP(U)) \arrow[r,"\Log",leftarrow] \arrow[d,twoheadrightarrow,"\bigshuffle(q)",swap] & \hck(U) \\
\bigshuffle(\bigodot(U)) \arrow[r,"\qLog",leftarrow] & \bigqshuffle(U) \arrow[u,hookrightarrow,"\widetilde\iota",swap]
\end{tikzcd}
\end{equation}
where \(q \coloneqq \widetilde{\mathfrak a}|_{\cP(U)}\) and \(\widetilde \iota\) was defined in \Cref{prop:arb.intr} (recall that this is not a Hopf morphism). \(\bigshuffle(q)\) is a left inverse to \(\bigshuffle(i)\) and \(\widetilde \iota\) is a right inverse to \(\widetilde{\mathfrak a}\) (as shown in \Cref{prop:arb}), the two vertical maps appearing in the statement of the theorem. We show commutativity, using the first statement and again cofreeness
\begin{align*}
\pi_1 \circ \qLog &= \pi_1 \circ \e_{\qshuffle} \\
&= \bigshuffle(q) \circ \pi \circ \e_\mathrm{CK} \circ \widetilde \iota \\
&= \bigshuffle(q) \circ \pi \circ \Log \circ \widetilde \iota \\
&= \pi_1 \circ \bigshuffle(q) \circ \Log \circ \widetilde \iota
\end{align*}
where we have used that
\begin{align*}
\bigshuffle(q) \circ \pi \circ \e_\mathrm{CK} \circ \widetilde \iota (s_1 \otimes \cdots \otimes s_n) &= \frac{(-1)^n}{n} s_1 \odot \cdots \odot s_n \\
&= \pi_1 \circ \e_{\qshuffle} (s_1 \otimes \cdots \otimes s_n)
\end{align*}
by the same argument used in the first part of this proof. Incidentally, note how this proof will show that \(\qExp\) factors as a Hopf isomorphism in terms of \(\Exp\), but that \(\qLog\) only factors as a coalgebra morphism in terms of \(\Log\). Now, returning to the rectangle in the statement, we have
\begin{align*}
\qExp &= \qLog^{-1}\\
&= (\bigshuffle(q) \circ \Log \circ \widetilde \iota\hspace{0.1em})^{-1} \\
&= \widetilde{\mathfrak{a}} \circ \Exp \circ \bigshuffle(i)
\end{align*}
using partial inverses.
	
Returning to the other faces, the triangle on the left commutes by definition. The triangle on the right commutes since \(\widetilde{\mathfrak{a}}_\odot\) restricted to \(U\)-decorated forests is equal to \(\widetilde{\mathfrak{a}}\), since they are defined in the same way. Commutativity of the lower parallelogram is the defining property of \(\mathfrak a \qExp\). Finally, commutativity of the upper parallelogram follows from a diagram chase:
\begin{align*}
\hck(j) \circ \Exp \circ (\bigshuffle(i) \circ \widetilde{\mathfrak{a}}) &=  \widetilde\iota_\odot \circ \qExp \circ \mathfrak a \\
&= \hck(j) \circ \widetilde\iota \circ \qExp \circ \mathfrak a  \\
&= \mathfrak{a} \qExp
\end{align*}
where \(\widetilde\iota_\odot\) is the canonical right inverse to \(\widetilde{\mathfrak{a}}_\odot\) defined in analogy with \(\widetilde \iota\).
\end{proof}
\end{theorem}

\begin{remark}
It is also possible to check commutativity of the rectangle in the main diagram more directly, by showing that \(\qExp\) satisfies the same
\begin{align*}
  \Exp^1_n = - \sum_{1 \leq k < n} {\Exp}^1_k \circ{\Log}.\quad\qedhere
\end{align*}
\end{remark}

	\subsection{Relationship with previous approaches}\label{subsec:relate} 
As mentioned in the introduction, the problem of transforming a branched rough path into a geometric one over a larger
space is not new. It was first tackled by Hairer and Kelly \cite{HK15}, and later by Boediharjo and Chevyrev
\cite{BC2019}. In this subsection we review these approaches and compare them with ours. To summarise, the picture that
will emerge is that our solution is quite distinct from that of \cite{HK15}, which is not canonical and relies on
abstract existence results. \cite{BC2019} provides an algebraic solution, but still one that is not canonical, as it depends on a choice of a basis. The dual of our \(\Exp\) can be identified as a specific choice of the maps considered therein; however, there is no guarantee that a generic basis will produce a \emph{natural} isomorphism, as simple examples demonstrate. In particular any sort of uniqueness is out of the question without the constraint of naturality; in fact this is what guarantees uniqueness in the context of quasi-shuffle.
Indeed, all natural Hopf morphisms \(\Phi\colon\bigshuffle(\bigodot(U))\to\bigqshuffle(U)\) are induced by linear maps
\(\varphi\colon\bigotimes(\bigodot(U))\to\bigodot(U)\) of the form
\[
  \varphi(s_1\otimes\dotsm\otimes s_n) = c_n(s_1\odot\dotsm\odot s_n)
\]
for a sequence of constants \(c_1, c_2,\dotsc\in\mathbb{R}\). The induced map is an isomorphism if and only if \(c_1=1\).
This is a highly non-trivial, see \cite{FPT16}.

\begin{proposition}\label{prop:hoffUniq}
  Hoffman's exponential \(\qExp\), defined in \eqref{eq:hoffman}, is the unique natural Hopf transformation \(\bigshuffle\circ\bigodot\Rightarrow\bigqshuffle\).
\end{proposition}
\begin{proof}
  We recall that \(\varphi\) induces a bialgebra morphism if and only if
  \[
    \varphi(x\shuffle y) = \varphi(x)\odot\varphi(y).
  \]
  In particular, considering any element of \(s\in\bigotimes(\bigodot(U))\) we see that
  \[
    c_n s^{\odot n} = \varphi(s^{\otimes n}) = \frac{1}{n!}\varphi(s^{\shuffle n}) = \frac{1}{n!}s^{\odot n}.\qedhere
  \]
\end{proof}

\begin{remark}\label{rem:uniq}
  We conjecture that an appropriate naturality constraint also yields uniqueness of \(\Exp\) in the case of \(\hck\). However, a similar argument to the proof of the previous proposition does not immediately go through, since the structure of the maps is more involved.
\end{remark}

The approach of \cite{HK15} is, in fact, partly algebraic. The authors consider a natural map
\[
\psi \colon \hck \to \bigshuffle (\cT), \qquad \psi(\scrt) = \scrt + \psi(\scrt') \otimes \scrt'', \ \scrt \in \cT
\]
which is fixed by the further requirement that it be an algebra morphism (this is because \(\hck\) is free commutative
over trees \(\cT\)). The map \(\psi\) plays the role of our \(\Log\), since it maps \(\hck\) to a shuffle algebra over a larger
space, and is injective, since it is injective on \(\cT\) as seen from the leading term. 
But, unlike \(\Log\), dimensionality (along the grading) immediately implies it is not surjective.
Then \(\cT\)-valued geometric rough paths \(\overline \bX\) satisfying \(\psi^*(\overline \bX) = \bX\) are considered,
and the authors prove a theorem analogous to \Cref{cor:itoStrat} for such rough paths. Since \(\psi\) is not
invertible, a recursive procedure for proving the existence of such a rough path is employed, and this is achieved via
the non-constructive existence theorem of \cite{LV07}. The lack of uniqueness of this rough path is expressed dually by
the fact that there is redundancy in the choice of vector fields in the geometric RDE \cite[Remark 5.9]{HK15}. One may
try to avert these difficulties by asking whether \(\psi\) admits a Hopf left inverse, which one could then use to
define \(\overline \bX\) directly, but it is highly unclear whether such a map should exist, especially with the extra requirement of naturality.

In \cite{BC2019}, the authors use the result, proved independently (in different contexts) in \cite[Theorem 8.4]{Foi2002} and \cite[Corollary 6.3]{Ch10}, that \(\cT(U)\) is free as a Lie algebra over some space \(\mathcal B(U)\). Specifically, the latter can be summarised by stating that the upper triangle in the following diagram
\begin{equation}
\begin{tikzcd}[column sep = large,row sep = large]
	&\cat{pre-Lie} \arrow[d] \arrow[dd,"(\bigodot(\, \cdot \, ){,} \circledast)",bend left=45] \\
	\cat{Vec} \arrow[ru, "(\cT{,} \curvearrowleft)"] \arrow[r, "(\cT{,}{[} \, \cdot \, {,} \, \cdot \, {]}_\star)"] \arrow[d,"\mathcal B"] &\cat{Lie} \arrow[d,"\mathcal U",swap] \\
	\cat{Vec} \arrow[r,"\bigotimes",swap] \arrow[ru,"\mathcal L",swap] &\cat{Alg}
\end{tikzcd}
\end{equation}
commutes, where \(\bigotimes\) is the free algebra functor and \(\mathcal U\) is the universal enveloping algebra
functor. Commutativity of the
lower triangle is the Poincar\'e-Birkhoff-Witt theorem, i.e.\ that the universal enveloping algebra of a free Lie
algebra is free. Note that the cited article is
about free pre-Lie algebras, and does not consider the Connes-Kreimer or Grossman-Larson Hopf algebras. The connection
is made thanks to the fact that \(\cT\) with tree grafting \(\curvearrowleft\) coincides with the free pre-Lie algebra
functor and the functor of Grossman-Larson primitives, compatibly with their Lie bracket given as antisymmetrisation of
\(\star\): by Milnor-Moore it then follows that \(\hgl = \mathcal U\circ\cT\) which is thus equal to \(\bigotimes \circ \mathcal B\) by commutativity of the above diagram.
The Oudom-Guin theorem provides an explicit expression for the product in \(\hgl\).

This point of view is related to ours as follows: \(\mathcal B\) (or at least a particular natural choice for it) coincides with \(\e_\mathrm{GL}(\cQ)\). It is enough to show that \(\hgl\) is free over \(\e_\mathrm{GL}(\cQ)\), since \(\e_\mathrm{GL}\) is valued in \(\mathrm{Prim}(\hgl) = \cT\) (the Lie algebra associated to the free associative algebra is the canonical model for the free Lie algebra, so the upper triangle commutes). Then \Cref{thm:iso} can be read as stating that \(\hck\) is cofree over the projection \(\pi \circ \eck\), and thus dually \(\hgl\) is free over the injection \(\egl \circ \pi^*\), i.e.\ over \(\e_\mathrm{GL}(\cQ)\).

We wish to stress the following point, which is central to the scope of this paper. There are uncountably many spaces \(\cB\) over which a free algebra, in our case \(\hgl\), is free. Requiring that the isomorphism \(\hgl(U) \cong \bigotimes(\cB(U))\) be natural in \(U\) is a much more restrictive property. In \cite[Theorem 2.3]{BC2019} the authors focus on finding such a basis composed of trees (as opposed to linear combinations of them), and in \cite[\S 6]{BC2019} compute the explicit free basis of the truncated \(\hgl^{2}\): having fixed a frame (ordered basis) of \(U\), take \(\cB(U) = \big\{1, \ \Forest{[\gamma]}, \ \Forest{[\beta[\alpha]]}\big\}_{\gamma, \alpha \leq \beta}\). The resulting isomorphism \(\hgl \to \bigotimes(\cQ)\), however, is not natural, since any coordinate permutation changes the order of the basis. Moreover, order-2 truncation makes it possible to compare its dual with Hoffman's exponential, and the two do not coincide
\[
\hck(U) \ni \Forest{[\beta[\alpha]]} \mapsto \Forest{[\alpha]} \otimes \Forest{[\beta]} + \deltaup_{\alpha \leq \beta} \pi(\Forest{[\alpha]} \Forest{[\beta]}) \neq \Forest{[\alpha]} \otimes \Forest{[\beta]} + \tfrac 12 \pi(\Forest{[\alpha]} \Forest{[\beta]}) = \qExp\big(  \Forest{[\beta[\alpha]]}  \big)
\]
In terms of stochastic integration and rough paths, this means that the geometric rough path associated to the It\^o rough path via this isomorphism will not coincide with the Stratonovich rough path, rather its second order terms \(\int_{s < u < v < t}\dif X^\alpha_u \dif X^\beta_v + \deltaup_{\alpha \leq \beta}[X]^{\alpha\beta}\) will depend on the particular order of the coordinates chosen. In fact, this example, to be compared with the next calculation, shows that bases constituted of single trees will never be natural (as this will immediately fail at order 2).
\begin{example}[\(\egl(\cQ)^3\)]
We compute the free basis of $\hgl$ up to level 3: this is a straightforward rearrangement after applying definitions of $\e_\mathrm{GL}$ and $\pi^*$ (or \Cref{ex:dual})
\begin{align*}
	\egl(\cQ)^3 ={}\mathrm{span} \bigg\{ &\Forest{[\c]}, \ \Forest{[\a[\b]]} + \Forest{[\b[\a]]} , \ \Forest{[\a[\b][\c]]} + \Forest{[\b[\a][\c]]} + \Forest{[\c[\a][\b]]} + 2 \sum_{\sigma \in \mathbb S_3} \Forest{[\sigma(\gamma)[\sigma(\beta)[\sigma(\alpha)]]]} ,\\
	&\Forest{[\alpha[\beta][\gamma]]}
	-\Forest{[\gamma[\alpha][\beta]]}
	+\Forest{[\alpha[\beta[\gamma]]]}
	-\Forest{[\gamma[\beta[\alpha]]]}
	+\Forest{[\gamma[\alpha[\beta]]]}
	-\Forest{[\alpha[\gamma[\beta]]]}
	\bigg\}
\end{align*}
We note that this is a strict subspace of $\cT^3$.
\end{example}

The following two considerations further motivate our interest in considering natural isomorphisms as they relate to rough analysis: the first involves probability theory, while the second comes from differential geometry.
\begin{remark}[Preservation of coordinate permutation-invariance in law]
Many rough paths considered in the literature (see \cite[Part III]{FV10}) involve a process with i.i.d.\ components w.r.t.\ a chosen basis \(B\) of the underlying space. If constructed \say{naturally}, one may expect the resulting rough path terms to be invariant, in law, under coordinate permutation. By this we mean, if the rough path is branched, that \(\bX^\scrf \sim \bX^{\cH(\sigma)(\scrf)}\) where \(\cH(\sigma) \in \mathrm{Aut}(\hck)\) is the map induced by some \(\sigma \in \mathbb S_B\) (\emph{not} for \(\sigma\) a permutation of the vertices of \(\scrf\): the corresponding statement will of course not hold in general unless \(\sigma \in \mathbb S_\scrf\)). Then, when transforming it to a geometric rough path using an isomorphism \(\Phi \colon \hgl \to \tens(\cB)\), we may wish for this property to continue to hold for the transformed rough path. This will hold if the transformation is natural, since for a word \(w \in \bigotimes (\cP)\)
\[
\Phi(\bX)^{\bigotimes \circ \cP(\sigma)(w)} = \langle \Phi^*(\tens\circ \cP(\sigma)(w)), \bX \rangle = \langle \hck(\sigma)(\Phi^*(w)), \bX \rangle \sim \Phi(\bX)^\scrf
\]
but may not in the opposite case  (as illustrated by the order-2 example above). 
\end{remark}

\begin{remark}[It\^o-Stratonovich transformations on manifolds]
In \cite{F22} one of the authors introduced a definition of branched rough path on a manifold, and used it to study integrals and differential equations (the latter in the quasi-geometric case), using the framework of bracket extensions. Using the framework of this article, this definition would be available without invoking bracket extensions, cf.\ \Cref{rem:kelly}: it involves considering a collection of branched rough paths \(\bX^i\) and an atlas indexed by the same set, requiring the compatibility condition \((\varphi_j^{-1} \circ \varphi_i)_*\bX^i = \bX^j\) on the overlap of the two charts, which hinges on defining the \emph{pushforward} of a rough path through a smooth map, done earlier in the paper using a natural construction involving tree grafting (and done in \cite{LCL07} for geometric rough paths using ordered shuffles). If such a manifold-valued branched rough path is defined, we may ask about transforming it to a geometric one using an isomorphism \(\Phi\) as above. This, however, will only be possible if such an isomorphism commutes with pushforwards
\[
f_* \Phi(\bX) = \Phi(f_*\bX) \quad \implies \quad \Phi(\bX^j) = (\varphi_j^{-1} \circ \varphi_i)_* \Phi(\bX^i)
\]
While we leave an in-depth discussion of these topics for further work, we point out that this can only be expected to hold if \(\Phi\) is natural. If it is not, the isomorphism will depend on the system of local coordinates chosen, i.e.\ it will not produce a manifold-valued geometric rough path.
\end{remark}

\section{The one-dimensional case}\label{sec:one}

\subsection{Kailath-Segall polynomials}
In this section we develop a generalization of the Kailath-Segall polynomials \cite{Segall1976OrthogonalFO} to branched rough paths.
The goal is to write the quantities representing the iterated integrals 
\[
\int_{s < u_1 < \ldots < u_n < t} \dif X^p_{u_1} \cdots \dif X^p_{u_n}
\]
with \(p \in \cP\), as polynomials in increments of \(X\) and its \say{higher-order variations}.
In fact, we will show that these polynomials are universal in some sense, because they arise as the image of the Eulerian idempotent in the free Hopf algebra over a sequence of symbols behaving like divided powers.

Let us consider a sequence of formal symbols \(\mathbf{x}=(x_n:n\ge 1)\), and let \(\mathbb{R}[\mathbf{x}]\) be the free polynomial algebra generated by them.
We define a coproduct by setting
\[
	\Delta x_n=\sum_{j=0}^n x_j\otimes x_{n-j}
\]
with the convention \(x_0\coloneqq 1\). It is clear that it is a commutative, cocommutative bialgebra.
We give it a grading by \(|x_n|=n\), therefore its also connected.
Hence, it is a Hopf algebra.

It is immediate that for any \(k>1\) we have
\[
	\Delta^{(k-1)}x_n=\sum_{i_1+\dotsb+i_k=n}x_{i_1}\otimes\dotsm\otimes x_{i_k},
\]
so the associated Eulerian idempotent (see \Cref{def:euler}) satisfies
\[
	\mathrm{e}(x_n)=\sum_{k=1}^n\frac{(-1)^{k-1}}{k}\sum_{i_1+\dotsb+i_k=n}x_{i_1}\dotsm x_{i_k}.
\]
Given the product is commutative, we may rearrange terms and obtain the perhaps better-known expression
	\begin{equation}\label{eq:Pn}
		P_n \coloneqq \mathrm{e}(x_n) = \sum_{a_1+2a_2+\dotsb+na_n=n}(-1)^{a_1+\dotsb+a_n-1}\frac{(a_1+\dotsb+a_n-1)!}{a_1! \dotsm a_n!}x_1^{a_1} \dotsm x_n^{a_n}
	\end{equation}
Since we are in a cocommutative Hopf algebra, by \Cref{prop:e}, \(P_n\) is primitive and since it is also commutative, again by \Cref{prop:e} the collection \(P_n\) is a set of algebra generators. In particular, we can rewrite \(x_n\) as a polynomial in these variables.
In fact, since \(\mathrm{id}=\exp_*(\e)\) we may write

\[
x_n = \sum_{k=1}^n\frac{1}{k!}\sum_{i_1+\dotsb+i_k=n}P_{i_1}\dotsb P_{i_k}.
\]
Again, applying the same reordering as before we may write
\begin{equation}\label{eq:KSpolynomials}
x_n = \sum_{a_1+2a_2+\dotsb+na_n=n}\frac{P_{1}^{a_1}}{a_1!}\dotsb \frac{P_n^{a_n}}{a_n!}.
\end{equation}
Expressions of this form in \(\hck\) were already known to Foissy \cite{Foi2002}.

We transfer these identities to \(\hck\).
Working over an arbitrary vector space \(U\) as in the previous sections, and given \(p \in \cP\), let \(\cH_{\top(p)}\) denote the free commutative algebra generated by the elements \(p^{\top n}\), \(n \in \mathbb N\) in \(\hck\). A special case is \(\cH_{\top(\Forest{[\c]})}\) is the sub-Hopf algebra generated by all ladders with nodes decorated with \(\c\). It is immediate that the mapping
\begin{equation}\label{eq:isoDivided}
  \mathbb{R}[\mathbf{x}] \to \cH_{\top(p)}, \quad x_n \mapsto p^{\top n}
\end{equation}
is a Hopf isomorphism. Combining all of these simple observations immediately yields the following generalisation of \cite[\S 9.1]{Foi2002}, in which ordinary ladders are replaced by the elements \(p^{\top n}\):
\begin{proposition}[Branched Kailath-Segall polynomials]\label{thm:KSidentity}
	For \(p \in \cP\), \(p^{\top n}\) satisfies the expression for \(x_n\) \eqref{eq:KSpolynomials} where \(P_n\) is given by \eqref{eq:Pn} with \(p^{\top k}\) replacing \(x_k\). Evaluating a branched rough path \(\bX\) on this identity expresses \(\langle p^{\top n}, \bX_{s,t} \rangle\) as a polynomial in the path increments \(\langle P_k, \bX_{s,t} \rangle\).
\end{proposition}

Continuing to use \(P_n\) as in the above proposition, we are able to obtain a classical representation formula for
iterated integrals of semimartingales as a special case. We note that this formula continues to hold for discontinuous
semimartingales, i.e.\ processes not explicitly considered in the rest of this article; nonetheless, these same
algebraic identities continue to hold. See for example \cite{avram1986symmetric,Segall1976OrthogonalFO,solut08}.
Discontinuous semimartingales, such as
general Lévy processes, are important examples in this regard, as they allow us to attain the identities in the fullest
generality.

The following result is a generalization of \cite[Theorem 4.2]{J2011}, which corresponds to the case \(m=1\).
\begin{corollary}[Classical Kailath-Segall polynomials for quasi-geometric rough paths]\label{cor:classicalKS}
Let \(\bX\) be a one-dimensional quasi-geometric rough path. Then 
\[
\langle \mathscr \pi(\mathscr r_m)^{\top n}, \bX_{s,t} \rangle = (-1)^n\sum_{a_1+2a_2+\dotsb+na_n=n} \frac{(-1)^{a_1 + \ldots + a_n}}{a_1! \cdot 2^{a_1} a_2! \cdots n^{a_n} a_n!} \langle \pi(\scrr_m), \bX_{s,t} \rangle^{a_1} \cdots \langle \pi(\scrr_{mn}) , \bX_{s,t} \rangle^{a_n}
\]
\begin{proof}
By \Cref{Thm_quasi_geo}, defining 
\begin{equation}\label{eq:Ptilde}
\widetilde P_k \coloneqq \frac{(-1)^{k-1}}{k}\pi(\scrr_k)
\end{equation}
we have that \(\langle P_k, \bX \rangle = \langle \widetilde P_k, \bX \rangle\) since \(P_k = \pi(P_k)\) and by \Cref{Thm_quasi_geo}. We conclude by replacing \(P_k\) with \(\widetilde P_k\) in \eqref{eq:KSpolynomials}.
\end{proof}
\end{corollary}

\begin{remark}
The classical Kailath-Segall polynomials are written in terms of the variables \(\pi(\scrr_k)\). One might ask for the
branched Kailath-Segall polynomials for \(p^{\top n}\) to be expressed explicitly in terms of images of \(\pi\). Following
\eqref{eq:ibp}, this can be interpreted as saying that for geometric rough paths, \(\bX^{\scrl_n} = X^n/n!\), for
quasi-geometric rough paths \Cref{cor:classicalKS} holds, and for more general branched rough paths, \(p^{\top n}\) is a polynomial in further images of \(\pi\). These can be obtained as in \Cref{cor:classicalKS}, i.e.\ by observing that \(P_n \in \cP\), \(P_n = \pi(P_n)\), i.e.\
\begin{equation}\label{eq:piPn}
	P_n	= \sum_{a_1+2a_2+\dotsb+na_n=n}(-1)^{a_1+\dotsb+a_n-1}\frac{(a_1+\dotsb+a_n-1)!}{a_1! \dotsm a_n!}\pi((p^{\top 1})^{a_1} \dotsm (p^{\top n})^{a_n}).
\end{equation}
and substituting this expression into \Cref{thm:KSidentity}. There are a couple of issues with this substitution, however. First of all, if one wishes for the variables of these polynomials to be \(\pi(\scrf)\) with \(\scrf \in \scrF\), the above expression needs to be re-worked and the polynomials are no longer universal in these variables, i.e.\ they will depend on \(p\) (or rather its \say{shape}). Moreover, since \(\{\pi(\scrf) \mid \scrf \in \scrF\}\) is not an independent set, such polynomials are not unique (although there will still be a natural choice that emerges from the above substitution) as they are instead in the quasi-geometric case, thanks to \Cref{lemma_independence}. Nevertheless, such polynomials can in principle be computed, in particular when \(p = \Forest{[\c]}\), in which case they represent the branched corrections to the formula \(\bX^{\scrl_n} = X^n/n!\) valid in the geometric case.
\end{remark}

We now return to the general setting to obtain a branched analogue of what is called the Kailath-Segall formula, i.e.\ a
recursive expression for \eqref{eq:KSpolynomials}.
Let us recall the definition of the Dynkin operator \(\EuScript{D}\colon H\to H\) over a graded Hopf algebra \(H\):
\begin{equation}\label{eq:Dynkin}
	\EuScript{D}(h)\coloneq |h_{(2)}|\EuScript{S}(h_{(1)})h_{(2)}
\end{equation}
where \(\EuScript{S}\colon H\to H\) is the antipode, which, we recall, satisfies Takeuchi's formula \cite{tak71}: $\EuScript{S}(\mathbf{1}) = \mathbf{1}$ and
\begin{equation}\label{eq:Takeuchi}
	\EuScript{S}(h) = \sum_{n \geq 1} (-1)^n h^{(1)} \cdots h^{(n)}
\end{equation}
The following lemma relates the Dynkin and Eulerian idempotents in a very special case.
\begin{lemma}\label{lem:De}
	If \(H\) is commutative and cocommutative, \(\EuScript{D}(h) = |h|\e(h)\). 
	\begin{proof}
		Combining \eqref{eq:Dynkin} and \eqref{eq:Takeuchi} we have
		\begin{align*}
    \EuScript{D}(h) &= \sum_{n \geq 1} (-1)^n |h_{(2)}| (h_{(1)})^{(1)} \cdots (h_{(1)})^{(n)} h_{(2)} \\
			&= \sum_{n \geq 1} (-1)^{n-1} |h^{(n)}| h^{(1)} \cdots h^{(n)} \\
			&= \sum_{n \geq 1} (-1)^{n-1} \frac{|h^{(1)}| + ... +|h^{(n)}|}{n} h^{(1)} \cdots h^{(n)} \\
			&= |h|\e(h)
		\end{align*}
		where we have used that for each \(n\) the expression \(\sum_{(h)} h^{(\sigma(1))} \cdots h^{(\sigma(n))} = h^{(1)} \cdots h^{(n)}\) does not depend on \(\sigma \in \mathbb S_n\), which follows from commutativity and cocommutativity combined.
	\end{proof}
\end{lemma}

Formula \eqref{eq:Dynkin} can actually be inverted, since the antipode is the convolutional inverse of the identity map
in \(H\), that is, \(\EuScript{S} * \mathrm{id} = \mathrm{id} * \EuScript{S} = \mathbf{1} \circ \varepsilon\) for \(|h| > 0\). We readily obtain
\begin{align*}
	\mathrm{id}*\EuScript{D}&=\mathrm{id}*\EuScript{S}*|\,\cdot\,|=|\,\cdot\,|
\end{align*}
In the particular case of \(\mathbb{R}[\mathbf{x}]\), after applying this identity and \Cref{lem:De}, we get the following recursion:
\[
	nx_n=\sum_{k=1}^{n}kP_kx_{n-k}.
\]
Taking again the isomorphism \eqref{eq:isoDivided} (and reasoning as in \eqref{eq:Ptilde}) we have thus proved the following result.
\begin{theorem}[Branched Kailath-Segall formula]
Let \(\bX\) be a one-dimensional branched rough path
\begin{equation}\label{eq:KS_formula}
	\langle p^{\top n}, \bX_{s,t}\rangle=\sum_{k=1}^n\frac{k}{n}\langle P_{k}, \bX_{s,t}\rangle\langle p^{\top (n-k)}, \bX_{s,t}\rangle\,.
\end{equation}
As special case when \(\bX\) is quasi-geometric we recover the usual Kailath-Segall formula \cite[Theorem 1]{Segall1976OrthogonalFO}
\begin{equation}\label{eq:KS_formula.qs}
	\langle \pi(\scrr_m)^{\top n}, \bX_{s,t}\rangle=\sum_{k=1}^n\frac{(-1)^{k-1}}{n}\langle \pi(\scrr_{mk}), \bX_{s,t}\rangle\langle \pi(\scrr_m)^{\top (n-k)}, \bX_{s,t}\rangle\,.
\end{equation}
\end{theorem}

\begin{remark}
It is clear by symmetrising that the results of this chapter extend from \(\cH_{\top(p)}\) to the Hopf algebra \(\cH_\top\) generated as an algebra by the elements
\[
\frac{1}{n!} \sum_{\sigma \in \mathscr S_n} p_{\sigma(1)} \tee \ldots \tee p_{\sigma(n)}
\]
as \(p_1,\ldots,p_n\) range in \(\cP\). \(\cH_\top\) is the largest cocommutative sub-Hopf algebra of \(\hck\) \cite[Lemma 26]{F2017}, and hence the largest whose elements can be written as polynomials in primitives.
\end{remark}

One common application of Kailath-Segall polynomials is to expand the so-called Doléans-Dade exponential of a process.
Here, we consider its branched version, denoted by \(\mathcal{E}(X)\), and defined as the solution of the RDE
\begin{equation}\label{eq:ddexp}
  \mathrm{d}\mathcal{E}(X)_t = \mathcal{E}(X)_t \,\mathrm{d}X_t,\quad \mathcal{E}(X)_0=1.
\end{equation}
In the notation of previous section, this RDE is defined by the vector field \(F_{\Forest{[]}}\equiv\mathrm{id}\).
An immediate consequence, we obtain the following polynomial expansion.

\begin{proposition}\label{prp:ddexp}
  Let \(\bX\in\mathscr{C}^{\rho}\) be a branched rough path, and consider the stochastic exponential of its trace \(X\), defined by \eqref{eq:ddexp}.
  Then,
  \[
    \mathcal{E}(X)_t = \exp\left(\sum_{k=1}^{\infty}\langle P_k,\bX_{0,t}\rangle\right).
  \]
\end{proposition}
\begin{proof}
  We begin by noting that due to the so-called neo-classical inequality \cite[Definition 7.2]{Gub10} (see also \cite[Theorem 4]{Boe18}), there exists a (potentially large) constant \(K>0\), depending only on \(\rho\), such that \(|\bX_{s,t}^{\mathscr{f}}|\le \frac{K}{(\mathscr{f}!)^{1/\rho}}|t-s|^{|\mathscr{f}|/\rho}\) for every forest \(\mathscr{f}\in\mathcal{F}\).

  We note that since \(F_{\Forest{[]}}(y)\), it is immediate that the map \(\widehat{F}\) given by \Cref{prp:prelie.vf} satisfies \(\widehat{F}_{\ell_n}(x)=x\) for all \(n\ge 1\) and vanishes otherwise.
  Davie expansion yields
  \[
    \mathcal{E}(X)_{s,t}=\mathcal{E}(X)_s\sum_{n=1}^{\infty}\bX_{s,t}^{\ell_n},
  \]
  and we remark that the sum on the right-hand side is convergent since
  \[
    \sum_{n=1}^{\infty}|\bX_{s,t}^{\ell_n}|\le\sum_{n=1}^{\infty}\frac{K^n}{(n!)^{1/\rho}}|t-s|^{n/\rho}<\infty.
  \]
  Setting \(s=0\) we see that
  \[
    \mathcal{E}(X)_{t}=1+\sum_{n=1}^{\infty}\bX_{0,t}^{\ell_n}.
  \]
  By \eqref{eq:KSpolynomials} we obtain
  \begin{align*}
    \mathcal{E}(X)_t &= 1+\sum_{n=1}^{\infty}\left\langle\sum_{k=1}^n\frac{1}{k!}\sum_{i_1+\dotsb+i_k=n}P_{i_1}\dotsm P_{i_k},\bX_{0,t}\right\rangle\\
                     &= 1+\sum_{k=1}^{\infty}\frac{1}{k!}\left\langle\left( \sum_{n\ge 1}P_n \right)^k ,\bX_{0,t}\right\rangle\\
                     &= 1 + \sum_{k=1}^{\infty}\frac{1}{k!}\left(\sum_{n\ge 1}\langle P_n,\bX_{0,t}\rangle\right)^k.\qedhere
  \end{align*}

  \begin{remark}
    We note that this result can also be proven by appealing to the log-ODE method detailed in \Cref{cor:logode}.
    Indeed, by choosing the test function to be \(\varphi\equiv\mathrm{id}\) we obtain that \(Y_t=Z_1\) where
    \[
      \frac{\mathrm{d}}{\mathrm{d}\theta}Z_\theta=\langle\widehat{F}(Z_\theta),\log_\star S(\bX)_{0,t}\rangle.
    \]
    Since \(\widehat{F}\) vanishes on non-ladder trees and is the identity map otherwise, the right-hand side simplifies to
    \[
      \frac{\mathrm{d}}{\mathrm{d}\theta}Z_\theta=\left(\sum_{n=1}^{\infty}\langle \operatorname{e}(\ell_n),\bX_{0,t}\rangle\right)Z_\theta
    \]
    whence the result follows, after applying \eqref{eq:Pn}.
  \end{remark}

  \begin{remark}
    In the case \(\bX\) is quasi-geometric, \Cref{prp:ddexp,eq:KS_formula.qs} yield
    \[
      \mathcal{E}(X)_t=\exp\left(\sum_{k=1}^{\infty}\frac{(-1)^{k+1}}{k}\langle\pi(\mathscr{r}_k),\bX_{0,t}\rangle\right) 
    \]
    which is a well-known formula for the solution of \eqref{eq:ddexp} and may be obtained by a purely Itô calculus argument (cf. \cite[Proposition 2.6]{J2011}).
  \end{remark}
\end{proof}

\subsection{Construction of branched rough paths: stochastic examples at order 4}\label{subsec:four}
In this last subsection, we propose a general method for constructing and classifying examples of branched rough paths, and provide stochastic examples. The recipe is as follows:
\begin{enumerate}
	\item Depending on regularity and number of labels, identify a free commutative generating set \(\mathcal A\) of \(\hck\) constituted of single words, under the delimiter \(\top\), in the alphabet \(\cP\) (by this we actually mean \emph{single words}, not their linear combinations; this is somewhat like a Lyndon basis---such a set always exists, e.g.\ by picking a linear basis of $\cT$ constituted of single words). It may be that not all primitives in a basis of $\cP$ are needed as letters in $\mathcal A$, and the choice of which letters to use in which slots should be motivated by the next two steps.
	\item Assign to each primitive element \(p\) with \(|p| \leq \p\), and which appears as a letter in one of the elements of \(\mathcal A\), a path \(X^p\): this should be motivated by the specifics of the problem at hand.
	\item Using some kind of integration theory which satisfies the Chen identity and preserves regularity in the expected way, postulate the values of \(\bX^a\) for \(a \in \mathcal A\): since $\mathcal A$ generates $\hck$ as an algebra, this uniquely defines a branched rough path $\bX$.
	\item We can now check, using \Cref{Thm_quasi_geo}, whether \(\bX\) is geometric or quasi-geometric. This may involve computing \(\bX\) on elements not considered in its definition.
\end{enumerate}
The point is that cofreeness makes it possible to restrict the definition of \(\bX\) to a
certain number of linear iterated integrals, for which the Chen identity is simpler to state and verify than its
branched counterpart. Following these three points, we will be able to exhibit examples of rough paths which are
\emph{truly branched}, i.e.\ fail the quasi-geometricity test. We use this terminology in part to echo the notion of
\say{truly rough} \cite[Definition 6.3]{FH20}: the two share the feature of it being necessary to consider the rough path as member of a space of \say{worse-behaved} objects.

We proceed to put the above plan in action in the regime of one-dimensional branched rough paths of regularity \(\rho \in (4,5)\). The choice to consider this case is due to three simple
reasons: first of all, explicit computations for undecorated forests of order up to $4$ are tractable, and we will rely on those of \cite[Appendix A]{VRST22}. Secondly, we will focus on fractional Brownian motion with Hurst parameter $H = 1/4$, for which rough path lifts are problematic to define in multidimensional case \cite{CQ02}, but not in the one-dimensional case. Finally, at level four, as
described in \Cref{Thm_quasi_geo} branched rough paths strictly contain quasi-geometric rough paths, which makes the change-of-variable formulae more interesting and of a form not yet considered. 

Let \(X\) be a one-dimensional
fractional Brownian motion (fBm) of Hurst parameter \(H=1/4\), the centered Gaussian process $X$ with covariance function $\mathbb E[X_s X_t] = \frac 12 (\sqrt{s} + \sqrt{t} - \sqrt{|t-s|})$, defined over a complete probability space \((\Omega,\mathcal{F}, \mathbb{P})\). It is well-known (e.g.\ see \cite{rogers97}) that \(X\) has a.s. infinite quadratic variation, that is, for
any \(t\in [0,T]\) and any sequence of partitions of \(\Pi_n\) of \([0,t]\) such that \(|\Pi_n|\to 0\) the random variable
\[[X]^n_t:=\sum_{[u,v]\in \Pi_n}(X_{v}- X_u)^2\]
diverges as \(n\to \infty\) in \(L^1(\mathbb{P})\). This property is shared with fractional Brownian motion of any Hurst parameter \(0<H<1/2\). However, in the critical case \(H=1/4\), taking a uniform partition of \([0,t]\) and considering the renormalised approximations to the quadratic variation
\begin{equation}\label{corrected_quadratic}
	\overline{[X]}_t^n\coloneqq [X]^n_t-\mathbb{E}[X]^n_t = \sqrt{\frac{t}{n}}\sum_{k=0}^{n-1}\left(\sqrt{\frac{n}{t}}(X_{(k+1)t/n}- X_{kt/n})^2-1 \right)
\end{equation}
it can be shown \cite[Theorem 1.1]{Nourdin09} that
\[\overline{[X]}_t^n\overset{(d)}{\longrightarrow}C_{1/4} W_t \quad (n\to +\infty)\]
for an independent standard Brownian motion \(W\) and \(C_{1/4}\approx 1535\) a fixed constant, which we drop from now on for simplicity.

Considering that in It\^o calculus with $X$ a semimartingale, $\langle \pi(\bullet \bullet), \bX_{0,t} \rangle$ is the ordinary quadratic variation, it is interesting to see what the resulting objects are in our case (with $X$ a $\frac 14$-fBm) if we postulate this term to be the weak limit of the renormalised approximations. We thus set
\[
\langle \pi(\bullet \bullet), \bX_{0,t} \rangle \coloneqq W_t
\]
with $W \perp \!\!\!\! \perp X$. Requiring that all additional terms be postulated in the It\^o or Stratonovich sense leads to the following classification.
\begin{proposition}\label{thm_stochastic}
	There exist (after having fixed $W$, up to modifications) four distinct \(\rho\in(4,5) \)-branched rough paths over \(X\), call them \(\bX^{S,S}\),  \(\bX^{I,I}\),  \(\bX^{S,I}\),  \(\bX^{I,S}\), satisfying the conditions 
	\begin{equation}
		\begin{split}\label{canonical_conditions}
			\langle \pi(\Forest{[]}\Forest{[]}), \bX_{s,t}\rangle  &= W_t-W_s, \quad \langle \Forest{[]}\tee\pi(\Forest{[]}\Forest{[]}), \bX_{s,t}\rangle  = \int_s^t(X_r-X_s)\dif W_r,\\
			\langle\pi(\Forest{[]}\Forest{[]})\tee\Forest{[]}, \bX_{s,t}\rangle  &= (X_t-X_s)(W_t-W_s)-\int_s^t(X_r-X_s)\dif W_r, \\
			\langle\pi(\Forest{[]}\Forest{[]}\Forest{[]}), \bX_{s,t}\rangle &= \langle\pi(\Forest{[]}\Forest{[]}\Forest{[]})\tee\Forest{[]}, \bX_{s,t}\rangle = \langle\Forest{[]}\tee\pi(\Forest{[]}\Forest{[]}\Forest{[]}), \bX_{s,t}\rangle = 0 \\[0.5ex]
		\end{split}
	\end{equation}
	and defining \(\langle\, \Forest{[]}\tee\Forest{[]}\tee\pi(\Forest{[]}\Forest{[]}), \bX_{s,t}\rangle\),  \(\langle\, \pi(\Forest{[]}\Forest{[]})\tee\pi(\Forest{[]}\Forest{[]}), \bX_{s,t}\rangle\) either as It\^o or Stratonovich integrals, i.e.\:
		\begin{align*}
		\langle\,\Forest{[]}\tee\Forest{[]} \tee  \pi(\Forest{[]}\Forest{[]}),
		\bX^{I,I}_{s,t}\rangle&=\langle\,\Forest{[]}\tee\Forest{[]} \tee  \pi(\Forest{[]}\Forest{[]}),
		\bX^{I,S}_{s,t}\rangle= \int_s^t  \frac{(X_r- X_s)^2- (W_r-W_s)}{2}\,\dif W_r,\\
		\langle\,\Forest{[]}\tee\Forest{[]} \tee  \pi(\Forest{[]}\Forest{[]}),
		\bX^{S,I}_{s,t}\rangle&=\langle\,\Forest{[]}\tee\Forest{[]} \tee  \pi(\Forest{[]}\Forest{[]}),
		\bX^{S,S}_{s,t}\rangle= \int_s^t  \frac{(X_r- X_s)^2- (W_r-W_s)}{2}\,{\circ\dif}W_r,\\
		\langle\,\pi(\Forest{[]}\Forest{[]})\tee  \pi(\Forest{[]}\Forest{[]}),
		\bX^{I,S}_{s,t}\rangle&=\langle\,\pi(\Forest{[]}\Forest{[]})\tee  \pi(\Forest{[]}\Forest{[]}),
		\bX^{S,S}_{s,t}\rangle= \int_s^t (W_r-W_s)\,{\circ\dif} W_r = \frac{(W_t-W_s )^2}{2},\\
		\langle\,\pi(\Forest{[]}\Forest{[]})\tee  \pi(\Forest{[]}\Forest{[]}),
		\bX^{S,I}_{s,t}\rangle&=\langle\,\pi(\Forest{[]}\Forest{[]})\tee  \pi(\Forest{[]}\Forest{[]}),
		\bX^{I,I}_{s,t}\rangle= \int_s^t (W_r-W_s)\,\dif W_r.
	\end{align*}
\end{proposition}
We refer the reader to \Cref{app:four} for the proof, and elaborate on the conditions in the above theorem. By independence (zero covariation), the It\^o and Stratonovich integrals $\int_s^t(X_r-X_s)\dif W_r = \int_s^t(X_r-X_s)\circ\dif W_r$ agree, so there is not much extra choice to be made here. The identity on the second line of \eqref{canonical_conditions} is integration by parts, which again by independence is the same as taking the limit in $L^2$ of the Riemann sums $\sum_{[s,t]} (W_r-W_s)(X_t - X_s)$. From the algebraic identity
\[\pi(\Forest{[]}\Forest{[]})\;\Forest{[]}= \pi(\Forest{[]}\Forest{[]})\tee\Forest{[]}+ \Forest{[]}\tee\pi(\Forest{[]}\Forest{[]})+\pi(\Forest{[]}\Forest{[]}\Forest{[]}) \]
we deduce that \(\langle\pi(\Forest{[]}\Forest{[]}\Forest{[]}), \bX_{s,t}\rangle=0\), which motivates setting all terms involving $\pi(\Forest{[]}\Forest{[]}\Forest{[]})$ to zero.

We can now bring the plan set out earlier to its conclusion. Recall from \eqref{eq:fourQuasi} that (in the current setting) $\langle \pi(\Forest{[[]]}\Forest{[[]]}) , \bX \rangle$ is the unique obstruction to $\bX$ being quasi-geometric. The algebraic identity 
\[\pi(\Forest{[[]]}\Forest{[[]]})= \frac{1}{3}\left(\pi(\Forest{[]} \Forest{[]}\Forest{[]})\tee \Forest{[]}+\Forest{[]}\tee\pi(\Forest{[]}\Forest{[]}\Forest{[]}) - \pi(\Forest{[]}\Forest{[]}\Forest{[]}) \Forest{[]} +\pi(\Forest{[]}\Forest{[]})\;\pi(\Forest{[]}\Forest{[]})-2 \pi(\Forest{[]}\Forest{[]})\tee\pi(\Forest{[]}\Forest{[]}) \right)\]
can be interpreted as explaining the obstruction as a sum of two defects of integration-by-parts. Since in all cases $\bX$ evaluated in terms containing $\pi(\Forest{[]}\Forest{[]}\Forest{[]})$ vanishes, it follows that
\[
\langle \pi(\Forest{[[]]}\Forest{[[]]}), \bX \rangle = \frac 13 \langle \pi(\Forest{[]}\Forest{[]})\pi(\Forest{[]}\Forest{[]}) - 2 \pi(\Forest{[]}\Forest{[]}) \tee \pi(\Forest{[]}\Forest{[]}), \bX \rangle
\]
for all four choices. This reinforces the intuition that lack of quasi-geometricity corresponds to higher-order failure of integration-by-parts identities. Therefore, while $\bX^{I,S}$ and $\bX^{S,S}$ are quasi-geometric, \(\bX^{I,I}\) and \(\bX^{S,I}\) satisfy
\[
\langle\pi(\Forest{[[]]}\Forest{[[]]}),\bX^{S,I}_{s,t}\rangle=\langle\pi(\Forest{[[]]}\Forest{[[]]}),\bX^{I,I}_{s,t}\rangle=
\frac{1}{3}\left((W_t-W_s)^2- 2\int_s^t(W_r-W_s)\,\dif W_r\right) = \frac{t-s}{3}.\]
\Cref{Thm_quasi_geo} then implies \(\bX^{I,I}\) and \(\bX^{S,I}\) do
not lie in the image of $\widetilde{\mathfrak a}{}^*$, making them examples of truly branched rough paths, nevertheless defined using tools from stochastic analysis in a rather straightforward way. We leave it as an interesting question for future work to relate the expansions for $\int \! f(X) \dif \bX$ for the four choices of $\bX$ (computed in \cref{prop:intfBm} for generic integrand) to the weakly converging discrete approximations to stochastic integrals of \cite[Theorem 1.1, Theorem 4.1]{NR09}.

\appendix 
\section{Calculation in \Cref{prop:failed}}\label{app:calc}
By definition of $\varphi$ one has
\begin{align*}
	\varphi (\Forest{[\b[\a]]})=\Forest{[\b[\a]]}- \frac{1}{2}\Forest{[\a]}\Forest{[\b]} + \frac{1}{2} \Forest{[\a]}\shuffle\Forest{[\b]} = \Forest{[\b[\a]]}- \frac{1}{2}\pi(\Forest{[\a]}\Forest{[\b]}) \,,
\end{align*}
Thanks to this identity and the trivial one $ \varphi(\Forest{[\c]})=\Forest{[\c]}\,$, we compute
\begin{align*}
	\varphi^{\otimes 2}  \dck  \bigg(\Forest{[\c[\b[\a]]]} \bigg)&=\varphi^{\otimes 2}  \left(\mathbf{1} \otimes \Forest{[\c[\b[\a]]]}+ \Forest{[\b[\a]]}\otimes \Forest{[\c]}+ \Forest{[\a]}\otimes \Forest{[\c[\b]]}+  \Forest{[\c[\b[\a]]]}\otimes \mathbf{1}\right)\\& =\mathbf{1} \otimes \varphi\left(\Forest{[\c[\b[\a]]]}\right)+  \varphi\left(\Forest{[\c[\b[\a]]]}\right)\otimes \mathbf{1}+\Forest{[\a]}\otimes \varphi(\Forest{[\c[\b]]})+ \varphi(\Forest{[\b[\a]]})\otimes \Forest{[\c]}\\&=\mathbf{1} \otimes \varphi\left(\Forest{[\c[\b[\a]]]}\right)+  \varphi\left(\Forest{[\c[\b[\a]]]}\right)\otimes \mathbf{1}+\Forest{[\a]}\otimes (\Forest{[\c[\b]]}- \frac{1}{2}\pi(\Forest{[\b]}\Forest{[\c]})) +(\Forest{[\b[\a]]}- \frac{1}{2}\pi(\Forest{[\a]}\Forest{[\b]}))\otimes\Forest{[\c]} \,.
\end{align*}
On the other hand, we compute $\varphi$ on the same tree via the expression for $\pi(\Forest{[\a]}\Forest{[\b]}\Forest{[\c]} )$ in \Cref{expl:P3}:
\begin{align*}
	\varphi \left(\Forest{[\c[\b[\a]]]}\right)&= \Forest{[\c[\b[\a]]]}-\frac{1}{2}(\Forest{[\c]} \Forest{[\b[\a]]}+\Forest{[\a]} \Forest{[\c[\b]]}) +  \frac{1}{3}\Forest{[\a]}\Forest{[\b]}\Forest{[\c]}  +\frac{1}{2} \Forest{[\c]}\shuffle(\Forest{[\b[\a]]}- \frac{1}{2}\Forest{[\a]}\Forest{[\b]})+ \frac{1}{2} (\Forest{[\c[\b]]}- \frac{1}{2}\Forest{[\b]}\Forest{[\c]})\shuffle\Forest{[\a]}+ \frac{1}{6} \Forest{[\a]}\shuffle\Forest{[\b]}\shuffle\Forest{[\c]}\\&=\Forest{[\c[\b[\a]]]}-\frac{1}{2}(\Forest{[\c]} \Forest{[\b[\a]]}+\Forest{[\a]} \Forest{[\c[\b]]})+  \frac{1}{3}\Forest{[\a]}\Forest{[\b]}\Forest{[\c]}+ \frac{1}{6}\sum_{\sigma\in \mathbb{S}_3}\Forest{[\sigma(\c)[\sigma(\b)[\sigma(\a)]]]}- \frac{1}{4}\left(\Forest{[\c]}\shuffle\Forest{[\a]}\Forest{[\b]}+\Forest{[\b]}\Forest{[\c]}\shuffle\Forest{[\a]} \right)+ \frac{1}{2}\left(\Forest{[\c]}\shuffle \Forest{[\b[\a]]} + \Forest{[\c[\b]]}\shuffle \Forest{[\a]}\right)\\&= 2\Forest{[\c[\b[\a]]]}+ \frac{1}{2}\left(\Forest{[\b[\a[\c]]]} +\Forest{[\b[\c[\a]]]} +\Forest{[\c[\a[\b]]]}+  \Forest{[\a[\c[\b]]]}\right)-\frac{1}{2}(\Forest{[\c]} \Forest{[\b[\a]]}+\Forest{[\a]} \Forest{[\c[\b]]})+  \frac{1}{3}\pi(\Forest{[\a]}\Forest{[\b]}\Forest{[\c]}) \\&\phantom{=}{}+ \frac{1}{12}(\Forest{[\c]}\tee\pi(\Forest{[\a]}\Forest{[\b]})+\pi(\Forest{[\b]}\Forest{[\c]})\tee\Forest{[\a]} + \pi(\Forest{[\a]}\Forest{[\b]})\tee\Forest{[\c]}+\Forest{[\a]}\tee\pi(\Forest{[\b]}\Forest{[\c]}))+ \frac{1}{3}\left(\Forest{[\b]}\tee\pi(\Forest{[\a]}\Forest{[\c]})+ \pi(\Forest{[\a]}\Forest{[\c]})\tee \Forest{[\b]}\right) \,.
\end{align*}

\noindent Applying the coproduct we derive
\begin{align*}
	\dck \varphi \bigg( \Forest{[\c[\b[\a]]]} \bigg) 
	&=\mathbf{1} \otimes \varphi\bigg(\Forest{[\c[\b[\a]]]}\bigg)+  \varphi\bigg(\Forest{[\c[\b[\a]]]}\bigg)\otimes \mathbf{1}+ \frac{1}{2}\wdck \, \left(4\Forest{[\c[\b[\a]]]}+ \Forest{[\b[\a[\c]]]} +\Forest{[\b[\c[\a]]]} +\Forest{[\c[\a[\b]]]}+  \Forest{[\a[\c[\b]]]}\right)\\&\phantom{=} -\frac{1}{2} \wdck\left(\Forest{[\c]} \Forest{[\b[\a]]}+\Forest{[\a]} \Forest{[\c[\b]]}\right)+ \frac{1}{3}\wdck\left(\Forest{[\b]}\tee\pi(\Forest{[\a]}\Forest{[\c]})+ \pi(\Forest{[\a]}\Forest{[\c]})\tee\Forest{[\b]}\right) \\&\phantom{=}+ \frac{1}{12}\wdck(\Forest{[\c]}\tee\pi(\Forest{[\a]}\Forest{[\b]})+\pi(\Forest{[\b]}\Forest{[\c]})\tee\Forest{[\a]} + \pi(\Forest{[\a]}\Forest{[\b]})\tee\Forest{[\c]}+\Forest{[\a]}\tee\pi(\Forest{[\b]}\Forest{[\c]}))\\&=\mathbf{1}\otimes \varphi\bigg(\Forest{[\c[\b[\a]]]}\bigg) +\varphi\bigg(\Forest{[\c[\b[\a]]]}\bigg)\otimes \mathbf{1}- \frac{1}{6}\left(\Forest{[\b]}\otimes \pi(\Forest{[\a]}\Forest{[\c]})+\pi(\Forest{[\a]}\Forest{[\c]})\otimes \Forest{[\b]} \right)+ \frac{1}{12}(\Forest{[\c]}\otimes \pi(\Forest{[\a]}\Forest{[\b]})+\pi(\Forest{[\b]}\Forest{[\c]})\otimes \Forest{[\a]})\\&\phantom{=}+  \left(\frac{3}{2}\Forest{[\b[\a]]}+\frac{1}{2} \Forest{[\a[\b]]}- \frac{1}{2}\Forest{[\a]}\Forest{[\b]}+ \frac{1}{12}\pi(\Forest{[\a]}\Forest{[\b]})\right)\otimes \Forest{[\c]}+ \Forest{[\a]}\otimes \left(\frac{3}{2}\Forest{[\c[\b]]}+\frac{1}{2}\Forest{[\b[\c]]}-\frac{1}{2}\Forest{[\b]}\Forest{[\c]} + \frac{1}{12} \pi(\Forest{[\b]}\Forest{[\c]})\right)\,.
\end{align*}	

\noindent Subtracting, we obtain
\begin{align*}
	&\dck \varphi \bigg( \Forest{[\c[\b[\a]]]}\bigg) -\varphi^{\otimes 2}  \dck  \bigg(\Forest{[\c[\b[\a]]]}\bigg) \\
	={}& \frac{1}{12}(\Forest{[\c]}\otimes \pi(\Forest{[\a]}\Forest{[\b]})+\pi(\Forest{[\a]}\Forest{[\b]})\otimes\Forest{[\c]}+\pi(\Forest{[\b]}\Forest{[\c]})\otimes \Forest{[\a]} + \Forest{[\a]}\otimes\pi(\Forest{[\b]}\Forest{[\c]})) -\frac{1}{6} (\Forest{[\b]}\otimes (\pi(\Forest{[\a]}\Forest{[\c]})+ \pi(\Forest{[\a]}\Forest{[\c]})\otimes \Forest{[\b]})
\end{align*}
which does not vanish when $\c$, $\b$ and $\a$ are all distinct.

\section{Proofs and details of \Cref{subsec:four}}\label{app:four}
	We begin by an explicit description of the truncated space of primitives.
	Note we truncate \(\mathcal{P}\) and \(\mathcal{Q}\) at different degrees since these are the minimal degrees required
	for the controlledness condition to hold, and thus enough to define rough integration.
	\begin{proposition}\label{prof:P4Q3}
		The vector spaces  \(\mathcal{P}^{4}= \bigoplus_{n=1}^4 \cP^{(n)}\) and \(\mathcal{Q}^3= \bigoplus_{n=1}^3\cQ^{(n)}\) are generated respectively by the elements
		\[\left\{  \Forest{[]}, \, \pi(\Forest{[]} \Forest{[]})\,, \pi(\Forest{[]} \Forest{[]}\Forest{[]})\,, \,\pi(\Forest{[]} \Forest{[]} \Forest{[]} \Forest{[]})\,,\pi(\Forest{[[]]} \Forest{[[]]}) \right\}\,,\quad \left\{  \Forest{[]}, \, \Forest{[]} \Forest{[]}\,,\Forest{[]} \Forest{[]}\Forest{[]}\, \right\}\,,\]
		which constitute  a basis. In the remainder of this section we will write $\hck$ for $\hck(\mathbb R)$ and similar.
	\end{proposition}
	\begin{proof}
		From the definition of the \(\pi\) operator, we can indeed compute \(\pi(\scrf)\) for any forest \(\scrf\) of
		cardinality smaller than \(4\), obtaining that \(\pi\) is non zero over the forests (c.f. \cite[Appendix A]{VRST22})
		\[\left\{  \Forest{[]}, \, \Forest{[]} \Forest{[]}\,, \Forest{[]} \Forest{[]}\Forest{[]}\,, \,\Forest{[]} \Forest{[]}
		\Forest{[]} \Forest{[]}\,,\Forest{[[]]} \Forest{[[]]} \,, \Forest{[]} \Forest{[]}\Forest{[[]]}\right\}.\]
		Using the result of \Cref{lemma_independence}, the identity \(\pi(\Forest{[[]]} \Forest{[[]]}) =  \pi(\Forest{[]}
		\Forest{[]}\Forest{[[]]})\) and remarking that \(\pi(\Forest{[[]]} \Forest{[[]]})\) is linearly independent from
		\(\pi(\Forest{[]} \Forest{[]} \Forest{[]}  \Forest{[]})\), we can easily conclude. The basis of \(\mathcal{Q}\) was
		already computed in \Cref{ex:dual}.
	\end{proof}

	Thus
	\begin{align*}
		\bigg\{\Forest{[]}, \, \pi(\Forest{[]} \Forest{[]}) \,,
		\Forest{[[]]}\,, \Forest{[[[]]]}\,,  \pi(\Forest{[]} \Forest{[]}) \tee \Forest{[]}\,,  \Forest{[]}\tee
		\pi(\Forest{[]} \Forest{[]})\,,  \pi(\Forest{[]} \Forest{[]}\Forest{[]})\,,\Forest{[[[[]]]]}\,, \pi(\Forest{[]}
		\Forest{[]}) \tee \Forest{[]} \tee \Forest{[]}\,,  \Forest{[]}\tee \pi(\Forest{[]} \Forest{[]}) \tee \Forest{[]}\,,\\
		\Forest{[]} \tee  \Forest{[]}\tee  \pi(\Forest{[]} \Forest{[]})\,,\pi(\Forest{[]} \Forest{[]}) \tee \pi(\Forest{[]}
		\Forest{[]}) \,,\pi(\Forest{[]} \Forest{[]}\Forest{[]})\tee \Forest{[]}\,,\Forest{[]}\tee\pi(\Forest{[]}
		\Forest{[]}\Forest{[]}) \,,\pi(\Forest{[]} \Forest{[]} \Forest{[]} \Forest{[]})\,,\pi(\Forest{[[]]}
		\Forest{[[]]})\bigg\}
	\end{align*}
	is a linear basis of $\hck^4$ adapted to the coalgebra structure, but still contains polynomial relations. The following proposition extracts an algebraic basis.

	\begin{proposition}\label{thm_generating_basis}
		Any element of \(\hck^{4}\) can be uniquely written as a polynomial w.r.t.\ the forest product with respect to the new variables
		\begin{equation}\label{new_variables}
			\begin{split}
				\mathcal{A}\coloneqq
				\bigg\{&  \Forest{[]}, \, \pi(\Forest{[]} \Forest{[]})\,, \pi(\Forest{[]} \Forest{[]}\Forest{[]})\,, \,\Forest{[]}\tee \pi(\Forest{[]} \Forest{[]}),\, \pi(\Forest{[]} \Forest{[]})\tee \pi(\Forest{[]} \Forest{[]})\,, \\& \pi(\Forest{[]} \Forest{[]}\Forest{[]})\tee\Forest{[]} \,,  \Forest{[]}\tee \pi(\Forest{[]} \Forest{[]}\Forest{[]})\,, \Forest{[]}\tee\Forest{[]}\tee \pi(\Forest{[]} \Forest{[]}) \bigg\}\,.
			\end{split}
		\end{equation}
	\end{proposition}

	\begin{proof}
		Denoting by \(\cT^4\) the set of trees of weight smaller than \(4\), to prove the result it is sufficient to show that the free commutative algebra over \(\cT^4\) coincides with the free commutative algebra over \(\mathcal{A}\). Since \(\cT^4\) and  \(\mathcal{A}\) have the same cardinality, it is sufficient to show that we can write every element of \(\cT^4\) as a polynomial in elements of \(\mathcal{A}\) w.r.t.\ the forest product. By using the product formula  \eqref{eq:bInfProd} combined with the operation \(\tee\) we can check the identities

		\begin{align*}
			\Forest{[[]]} &= \frac{1}{2}(\Forest{[]}\Forest{[]}- \pi(\Forest{[]}\Forest{[]}))\\
			\Forest{[[[]]]} &= \frac{1}{6}(\Forest{[]}\Forest{[]}\Forest{[]}- 3\,\pi(\Forest{[]}\Forest{[]})\Forest{[]}+ 2\,\pi (\Forest{[]}\Forest{[]}\Forest{[]} ))\\
			\Forest{[[][]]} &= \frac{1}{3}( \Forest{[]}\Forest{[]}\Forest{[]}- 3\,\Forest{[]}\tee \pi(\Forest{[]}\Forest{[]})- \pi (\Forest{[]}\Forest{[]}\Forest{[]}))\\
			\Forest{[[[[]]]]} &= \begin{multlined}[t]\frac{1}{24}( \Forest{[]}\Forest{[]}\Forest{[]}\Forest{[]}+
				4\,\pi(\Forest{[]}\Forest{[]}\Forest{[]})\Forest{[]}+\pi(\Forest{[]}\Forest{[]})\pi(\Forest{[]}\Forest{[]})+4\,\pi(\Forest{[]}\Forest{[]})\tee\pi(\Forest{[]}\Forest{[]})\\-6\,\Forest{[]}\Forest{[]}\pi(\Forest{[]}\Forest{[]})+ 4\,\pi(\Forest{[]}\Forest{[]}\Forest{[]})\tee\Forest{[]}+4\,\Forest{[]}\tee\pi(\Forest{[]}\Forest{[]}\Forest{[]}))\end{multlined}\\
			\Forest{[[[]][]]} &=\begin{multlined}[t]\frac{1}{24}(3 \Forest{[]}\Forest{[]}\Forest{[]}\Forest{[]}-6\,\Forest{[]} \Forest{[]}\pi(\Forest{[]}\Forest{[]}) -\pi(\Forest{[]}\Forest{[]})\pi(\Forest{[]}\Forest{[]})-4\, \pi(\Forest{[]}\Forest{[]})\tee\pi(\Forest{[]}\Forest{[]})\\-4\,\pi(\Forest{[]}\Forest{[]}\Forest{[]})\tee\Forest{[]}-4\,\Forest{[]}\tee\pi(\Forest{[]}\Forest{[]}\Forest{[]})+4\,\pi(\Forest{[]}\Forest{[]}\Forest{[]})\Forest{[]}-24\,\Forest{[]} \tee  \Forest{[]} \tee \pi(\Forest{[]}\Forest{[]}))\end{multlined}\\
			\Forest{[[][][]]}&=\frac{1}{4}( \Forest{[]}\Forest{[]}\Forest{[]}\Forest{[]}-\pi(\Forest{[]}\Forest{[]})\pi(\Forest{[]}\Forest{[]})-4\,\pi(\Forest{[]}\Forest{[]})\tee\pi(\Forest{[]}\Forest{[]}) - 12\,\Forest{[]}\tee \Forest{[]} \tee \pi(\Forest{[]}\Forest{[]})- 4\,\Forest{[]}\tee  \pi(\Forest{[]}\Forest{[]}\Forest{[]}))\\
			\Forest{[[[][]]]}&=\frac{1}{12}( \Forest{[]}\Forest{[]}\Forest{[]}\Forest{[]}-4\,\pi(\Forest{[]}\Forest{[]}\Forest{[]})\Forest{[]}+3\pi(\Forest{[]}\Forest{[]})\pi(\Forest{[]}\Forest{[]})\\ &\phantom{=} +12\,\Forest{[]}\tee \Forest{[]} \tee \pi(\Forest{[]}\Forest{[]})- 12\,\Forest{[]}\tee \pi(\Forest{[]}\Forest{[]})\Forest{[]}+ 12\,\Forest{[]}\tee \pi (\Forest{[]}\Forest{[]}\Forest{[]}))\,,
		\end{align*}
		concluding the argument.
	\end{proof}
	\begin{corollary}\label{cor_det_uniqueness}
		Any given function \(\mathbf{X}\colon \Delta_T\to \mathbb{R}^\mathcal{A}\) satisfying the Chen property and the regularity
		property in Definition \ref{def:brp} for some control \(\omega\) and \(\rho\in (4, 5)\) uniquely extends to a \(\rho\)-branched rough path \(\mathbf{X}\) controlled by \(\omega\).
	\end{corollary}
	\begin{proof}
		The regularity condition is obviously satisfied by the linear extension  \(\bX\colon\Delta_T\to\hgl^4\), and the character property follows from the fact that \(\mathcal{A}\) is an  algebraically independent set.  Thus, we only need to check that Chen's property is preserved. This follows from the fact that \(\bX_{s,t}\) and  \(\bX_{s,u}*\bX_{u,t}\) are two characters agreeing on \(\mathcal{A}\).
	\end{proof}

	We proceed with the proof of the main result of \Cref{subsec:four}.

	\begin{proof}[Proof of \Cref{thm_stochastic}]
		Keeping in mind that $\top$ is bracketed from the left, i.e.\ $\Forest{[]}\tee\Forest{[]} \tee  \pi(\Forest{[]}\Forest{[]}) = (\Forest{[]}\tee\Forest{[]}) \tee  \pi(\Forest{[]}\Forest{[]})$, and using that
		\(
		\langle \Forest{[]} \tee \Forest{[]}, \bX_{s,t} \rangle = \langle \Forest{[]}\Forest{[]} - \pi(\Forest{[]}\Forest{[]}) , \bX_{s,t}  \rangle/2
		\) the four possible choices are those in the statement of the proposition:
		\begin{align*}
			\langle\,\Forest{[]}\tee\Forest{[]} \tee  \pi(\Forest{[]}\Forest{[]}),
			\bX^{I,I}_{s,t}\rangle&=\langle\,\Forest{[]}\tee\Forest{[]} \tee  \pi(\Forest{[]}\Forest{[]}),
			\bX^{I,S}_{s,t}\rangle= \int_s^t  \frac{(X_r- X_s)^2- (W_r-W_s)}{2}\,\dif W_r,\\
			\langle\,\Forest{[]}\tee\Forest{[]} \tee  \pi(\Forest{[]}\Forest{[]}),
			\bX^{S,I}_{s,t}\rangle&=\langle\,\Forest{[]}\tee\Forest{[]} \tee  \pi(\Forest{[]}\Forest{[]}),
			\bX^{S,S}_{s,t}\rangle= \int_s^t  \frac{(X_r- X_s)^2- (W_r-W_s)}{2}\,{\circ\dif}W_r,\\
			\langle\,\pi(\Forest{[]}\Forest{[]})\tee  \pi(\Forest{[]}\Forest{[]}),
			\bX^{I,S}_{s,t}\rangle&=\langle\,\pi(\Forest{[]}\Forest{[]})\tee  \pi(\Forest{[]}\Forest{[]}),
			\bX^{S,S}_{s,t}\rangle= \int_s^t (W_r-W_s)\,{\circ\dif} W_r = \frac{(W_t-W_s )^2}{2},\\
			\langle\,\pi(\Forest{[]}\Forest{[]})\tee  \pi(\Forest{[]}\Forest{[]}),
			\bX^{S,I}_{s,t}\rangle&=\langle\,\pi(\Forest{[]}\Forest{[]})\tee  \pi(\Forest{[]}\Forest{[]}),
			\bX^{I,I}_{s,t}\rangle= \int_s^t (W_r-W_s)\,\dif W_r.
		\end{align*}
		Note that all differences are attributable to the classical It\^o-Stratonovich correction $\int_s^t(W_r-W_s)\circ \dif W_r = \int_s^t(W_r-W_s) \dif W_r + (t-s)/2$ occurring in different places. The conditions in the statement thus unambiguously define four distinct a maps \(\Delta_T\to \mathbb{R}^\mathcal{A}\), and by \ref{cor_det_uniqueness} it is sufficient to show that \(\bX\) satisfies the regularity assumption a.s.\ and the Chen property. For Chen, observe that It\^o and Stratonovich integration are additive on consecutive time intervals, which implies that the Chen identity holds on $\mathcal A$, and the statement follows by extension thanks to the fact that $\mathcal A$ is a free commutative generating set using the bialgebra properties, cf.\ \cite[Theorem 1.5]{BFGJ24}.

		Concerning the regularity in \(\rho\)-variation of all four possible models, it is sufficient to prove that the two-parameter process
		\[ (s,t)\to\left(X_t-X_s,W_t-W_s, \int_s^t(X_r-X_s)\,\dif W_r, \int_s^t(X_r-X_s)^2\,\dif W_r\right)\]
		admits a modification whose components have respectively a.s. finite \(1/\rho\), \(2/\rho\), \(3/\rho\) and
		\(4/\rho\) H\"older norm with \(\rho\in (4,5)\). For the first two components, we apply Kolmogorov's continuity theorem to the processes \(X\) and \(W\), thanks to the fact that
		\[\mathbb{E}|X_t-X_s|^p=C |t-s|^{\frac{p}{4}}, \quad \mathbb{E}|W_t-W_s|^p=C' |t-s|^{\frac{p}{2}} \]
		for any \(p\geq 1\) and some constants \(C,C'>0\). We obtain modifications of \(X\) and \(W\) (which we will denote in
		the same way) with the desired regularity properties and, in addition their Hölder seminorm are also \(q\)-integrable for any \(q\geq 4\).
		Passing to the stochastic integrals we only show how to prove that the third component is \(3/\rho\)-Hölder, since the
		last component follows similar reasoning.
		By applying the Burkholder-Davis-Gundy inequality we obtain for any \(p \geq 1\) the inequality
		\[ \mathbb{E} \left|\int_s^t(X_r-X_s)\,\dif W_r\right |^{2p}\lesssim \left|\int_s^t  \mathbb{E}|X_r-X_s|^2\,\dif r\right|^p\lesssim |t-s|^{\frac{3p}{2}},\]
		(we adopt the notation \(\lesssim\) to denote inequality up to a constant).
		We now apply \cite[Lemma 2.1]{nua_tindel11} to \(R_{s,t}\coloneqq\int_s^t(X_r-X_s)\,\dif W_r\) to get an estimate for its Hölder regularity.
		Indeed, we note that \(\delta R_{sut}\coloneqq R_{s,t}-R_{s,u}-R_{u,t}=(X_u-X_s)(W_t-W_u),\) so that
		\[
		\sup_{(s,t)\in\Delta_T}\frac{|\delta R_{sut}|}{|t-u|^{3/\rho}}\lesssim\|X\|_{1/\rho}\|W\|_{2/\rho}.
		\]
		We then obtain for any \((s,t)\in \Delta_T\), \(s\neq t\)
		\[\frac{|R_{s,t}|}{| t-s|^{3/\rho}}\lesssim
		\left(\int_{\Delta_T}\frac{|R_{uv}|^{2p}}{|u-v|^{6p/\rho+ 4}}\,\dif u \dif v\right)^{1/2p}+ \Vert X\Vert_{1/\rho} \Vert W\Vert_{2/\rho}\]
		where \(\Vert\cdot\Vert_{\alpha}\) denotes the \(\alpha\)-Hölder seminorm of a path. Since the right-hand side of this
		inequality does not depend on \((s,t)\), the estimate still holds when taking supremum of the left-hand side over the
		set \(D_N\coloneqq 2^{-N}\mathbb{N}\times 2^{-N}\mathbb{N}\cap\Delta_T\).
		By taking \(p\) sufficiently large so that
		\[\mathbb{E}\int_{\Delta_T}\frac{|R_{uv}|^{2p}}{|u-v|^{6p/\rho+ 4}}\,\dif u \dif v<\infty, \]
		which is possible since \(1/4-1/\rho>0\), we obtain that the sequence of random variables
		\[ Z_N=\sup_{(s,t)\in D_N}\frac{|R_{s,t}|}{| t-s|^{3/\rho}}\]
		is not only integrable but also uniformly integrable, and converges to a random variable \(Z_{\infty}\) with finite
		moment of order \(2p\). Since \(D_N\) is a monotone sequence of increasing sets whose limit \(D_{\infty}\) is a dense
		and countable subset of \(\Delta_T\), up to modifications the process \((s,t)\to \int_s^t(X_r-X_s)\,\dif W_r\) admits a modification which is a.s. \(3/\rho\) H\"older with also integrable H\"older seminorm.
	\end{proof}

	Combining the definition of the different branched rough paths in \Cref{thm_stochastic} we can also identify some rough integrals with other classical objects. In particular it follows from their definition that
	\[\int_s^t \mathbf{H}_r\,\dif \bX^{\pi(\bullet \bullet \bullet)}_r=0\]
	for \(\bX= \bX^{I,S},   \bX^{S,S}, \bX^{I,I}\) and \(\bX^{S,I}\).
	Moreover by using the properties of Riemann sums one has  also
	\[\int_s^t \mathbf{H}_r\,\dif \bX^{\pi(\bullet \bullet \bullet \bullet)}_r=\int_s^t \langle \mathbf{H}_r,\mathbf{1}^*
	\rangle\,\dif r,\quad \int_s^t \mathbf{H}_r\,\dif\bX^{\pi(\Forest{[[]]}\Forest{[[]]})}_r= \frac{1}{3}\int_s^t \langle
	\mathbf{H}_r,\mathbf{1}^* \rangle\,\dif r \]
	when \(\bX\) equals \(\bX^{I,I}\) or \(\bX^{S,I}\), and zero in the other two cases. Concerning the integral with
	respect to \(\pi(\bullet \bullet)\)  we have a specific identification of the rough integrals with some particular linear combination of stochastic integrals and classical integrals.
	\begin{proposition}\label{prop:intfBm}
		Let \(\rho\in(1/5,1/4)\). Suppose that \(\mathbf{H}\in\mathscr{D}_{\bX}^{4}\) a.s. is adapted to the joint  filtration of \(X\) and \(W\). The following identities holds a.s.:
		\begin{itemize}
			\item when \(\bX=\bX^{I,I}\)
			\begin{equation}\label{Ito-Ito-integral}
				\int_0^t \mathbf{H}_r\,\dif \bX^{\pi(\bullet \bullet)}_r= \int_0^t \langle \mathbf{H}_r, \mathbf{1}^*\rangle\,\dif W_r\,,
			\end{equation}
			\item when \(\bX= \bX^{I,S}\)
			\begin{equation}\label{Ito-strato-integral}
				\int_0^t \mathbf{H}_r\,\dif \bX^{\pi(\bullet \bullet)}_r= \int_0^t \langle \mathbf{H}_r, \mathbf{1}^*\rangle\,\dif
				W_r+\frac{1}{4}\int_0^t \langle \mathbf{H}_r,\Forest{[]} \Forest{[]} \rangle\,\dif r,
			\end{equation}
			\item when \(\bX= \bX^{S,I}\)
			\begin{equation}\label{strato-ito-integral}
				\int_0^t \mathbf{H}_r\,\dif \bX^{\pi(\bullet \bullet)}_r= \int_0^t \langle \mathbf{H}_r,\mathbf{1}^* \rangle\,\dif
				W_r-\frac{1}{4} \int_0^t \langle \mathbf{H}_r ,\Forest{[]}\star \Forest{[]} \rangle\,\dif r,
			\end{equation}
			\item when \(\bX= \bX^{S,S}\)
			\begin{equation}\label{strato-strato-integral}
				\int_0^t \mathbf{H}_r\,\dif \bX^{\pi(\bullet \bullet)}_r= \int_0^t \langle\mathbf{H}_r, \mathbf{1}^*\rangle\,\dif
				W_r- \frac{1}{4}\int_0^t \langle  \mathbf{H}_r, \Forest{[[]]} \rangle\,\dif r.
			\end{equation}
		\end{itemize}

	\end{proposition}
	\begin{proof}
		It is sufficient to show the result when \(t=1\). We begin by proving the identification \eqref{Ito-Ito-integral}
		using an argument similar to \cite[Prop 5.1]{FH20}, and using the expansion
		\begin{align*}
			\int_s^t \mathbf{H}_r\,\dif\mathbf{X}^{\pi(\bullet \bullet)}_r\approx{}&\langle	\mathbf{H}_s,\mathbf{1}^* \rangle  \langle \pi(\Forest{[]}\Forest{[]}),\mathbf{X}_{s,t}\rangle+ \mathbf{H}_s,\Forest{[]}\rangle  \langle \Forest{[]}\tee \pi(\Forest{[]}\Forest{[]}),\mathbf{X}_{s,t}\rangle\\&+\frac{1}{2} \langle\mathbf{H}_s, \Forest{[]} \Forest{[]} \rangle  \langle \pi(\Forest{[]}\Forest{[]})\tee \langle \pi(\Forest{[]}\Forest{[]}),\mathbf{X}_{s,t}\rangle+\langle\mathbf{H}_s, \Forest{[]} \star \Forest{[]}\rangle\langle \Forest{[]}\tee\Forest{[]}\tee \pi(\Forest{[]}\Forest{[]}),\mathbf{X}_{s,t}\rangle,
		\end{align*}
		implied by \Cref{lem:intCoords}. By definition of rough integral when  \(\bX= \bX^{I,I}\) this
		equals the a.s. limit
		\[\int_0^1 \mathbf{H}_r\,\dif \bX^{\pi(\bullet \bullet)}_r= \lim_{n\to +\infty}\sum_{[s,t]\in \Pi_n} \mathbf{Z}_{s,t}^{I,I}\]
		where  \(\mathbf{Z}_{s,t}^{I,I}\) is explicitly given by
		\[
		\mathbf{Z}_{s,t}^{I,I} =\begin{multlined}[t]\langle \mathbf{H}_s,\mathbf{1}^* \rangle (W_t-W_s)+ \langle
			\mathbf{H}_s, \Forest{[]}\rangle  \int_s^t(X_r-X_s)\,\dif W_r\\ + \frac{1}{2}\langle  \mathbf{H}_s,\Forest{[]}
			\Forest{[]}- \Forest{[]} \star \Forest{[]}\rangle  \int_s^t (W_r-W_s)\,\dif W_r+ \frac{1}{2} \langle \mathbf{H}_s,
			\Forest{[]} \star \Forest{[]} \rangle \int_s^t  (X_r- X_s)^2\,\dif W_r\end{multlined}
		\]
		and \(\Pi_n\) is a sequence of partitions of \([0,1]\) whose size \(|\Pi_n|\) converges to \(0\). Using the result in
		\cite[Chapter 4]{Yor99} we know that
		\[\int_0^t \langle \mathbf{1}^*, \mathbf{H}_r\rangle \dif W_r =\lim_{n\to +\infty}^{\mathbb{P}}\sum_{[s,t]\in
			\Pi_n}\langle \mathbf{H}_s,\mathbf{1}^* \rangle (W_t-W_s).\]
		Therefore the result follows by showing that
		\[
		\begin{split}
			\lim_{n\to +\infty}^{\mathbb{P}}\sum_{[s,t]\in \Pi_n}\mathbf{Z}_{s,t}^{I,I}=\langle \mathbf{H}_s, \mathbf{1}^*\rangle (W_t-W_s)
		\end{split}
		\]
		We will prove this convergence into two steps: First, we will show that all components of \(\mathbf{H}\) are uniformly
		bounded by the sequence of random variables
		\begin{equation}\label{eq_def_Zn}
			\begin{split}
				Z_n\coloneqq\sum_{[s,t]\in \Pi_n} \langle \mathbf{H}_s,\Forest{[]} \rangle  \int_s^t(X_r-X_s)\dif W_r &+
				\frac{1}{2}\langle\mathbf{H}_s, \Forest{[]} \Forest{[]}- \Forest{[]} \star \Forest{[]} \rangle  \int_s^t (W_r-W_s)
				\dif W_r\\&+ \frac{1}{2}\langle \mathbf{H}_s , \Forest{[]} \star \Forest{[]}\rangle \int_s^t  (X_r- X_s)^2\dif
				W_rm
			\end{split}
		\end{equation}
		which converges to \(0\) in \(L^2\).
		If \(\Pi_n=\{0=t_0<t_1^n<\ldots t_N^n=1\}\) one has the estimate
		\[
		\begin{split}
			&\Vert Z_n\Vert_{L^2}^2\lesssim \\& \bigg\Vert\sum_{k=0}^{N-1}\int_{t^n_k}^{t^n_{k+1}}(X_r-X_{t^n_k})\dif
			W_r\bigg\Vert_{L^2}^2+\bigg\Vert\sum_{k=0}^{N-1}\int_{t^n_k}^{t^n_{k+1}}(X_r-X_{t^n_k})^2\dif W_r\bigg\Vert_{L^2}^2+
			\bigg\Vert\sum_{k=0}^{N-1}\int_{t^n_k}^{t^n_{k+1}}(W_r-W_{t^n_k})\dif W_r\bigg\Vert_{L^2}^2.
		\end{split}
		\]
		Thanks to the It\^o isometry and elementary computations on the variance of \(X\) and \(W\) the random variables
		inside each sum are centred and uncorrelated, therefore the above sum coincides with
		\[\sum_{k=0}^{N-1}\int_{t^n_k}^{t^n_{k+1}}\mathbb{E}(X_r-X_{t^n_k})^2+ \mathbb{E}(X_r-X_{t^n_k})^4+
		\mathbb{E}(W_r-W_{t^n_k})^2 \dif r \lesssim \sum_{k=0}^{N-1}|t^n_{k+1}- t^n_k|^{3/2}\lesssim |\Pi_n|^{1/2}\,,\]
		thereby obtaining the desired convergence. To treat the general case we simply introduce the stopping time
		\[ \tau_M=\inf\{t\in  [0,1] \colon |\langle \mathbf{H}_{t},\Forest{[]} \rangle|\wedge |\langle \mathbf{H}_{t}, \Forest{[]} \Forest{[]}\rangle| \wedge |\langle \mathbf{H}_{t} , \Forest{[]} \star \Forest{[]}\rangle| \geq M\}\]
		and consider the sequence \(Z_n^M\) defined  as \(Z_n\) in \eqref{eq_def_Zn} where the components of
		\(\mathbf{H}_{t}\) are stopped with respect to the stopping time \(\tau_M\). Applying the previous part to the process
		\(Z_n^M\) we can use the fact that \(\tau_M\to + \infty\) a.s. as \(M \to + \infty\) to obtain that \(Z_n\) converges
		to \(0\) in probability. Passing to the other other three cases, namely \eqref{Ito-strato-integral}, \eqref{strato-ito-integral} and \eqref{strato-strato-integral}, the result follows easily from the first one by simply checking that one has the identities
		\begin{align*}
			\mathbf{Z}_{s,t}^{I,S}&= \mathbf{Z}_{s,t}^{I,I}+ \frac{1}{4} \langle \mathbf{H}_s,\Forest{[]} \Forest{[]}
			\rangle(t-s)\\
			\mathbf{Z}_{s,t}^{S,I}&= \mathbf{Z}_{s,t}^{I,I}- \frac{1}{4} \langle \mathbf{H}_s ,\Forest{[]}\star \Forest{[]} \rangle
			(t-s)\\
			\mathbf{Z}_{s,t}^{S,S}&= \mathbf{Z}_{s,t}^{I,I}- \frac{1}{4}\langle  \mathbf{H}_s, \Forest{[[]]} \rangle(t-s)
		\end{align*}
		as direct consequence of the It\^o formula applied to \((W_t-W_s)^2\).
	\end{proof}

	\begin{remark}
		Using \Cref{prop:intfBm} we obtain that the simple It\^o formula \eqref{eq:itoSimple} with $X^{S,S}$ and \(\mathbf{H}=\varphi(X)\) equals
		\begin{equation}\label{first_rough_identity}
			\varphi(X_{t})- \varphi(X_{s})=  \int_s^t D\varphi(X_r)\,\dif \bX^{S,S}_r
			+ \frac 12 \int_s^t\varphi''(X_r)\,\dif W_r,
		\end{equation}
		as a consequence of the fact that \(\langle  \D^2\varphi(X_t),\Forest{[[]]}\rangle= 0\). This identity is formally similar to a previous result on fractional Brownian motion with Hurst parameter $H=1/4$, see \cite[Theorem 4.1]{Nourdin09}. Indeed, using  Malliavin-calculus techniques the authors proved the convergence in law of the midpoint Riemann sums
		\begin{equation}\label{eq:rieSu}
		S^M_n\coloneqq\sum_{k=0}^{\lfloor n/2 \rfloor}\varphi'(X_{\frac{2k-1}{n}})( X_{\frac{2k}{n}}- X_{\frac{2k-2}{n}} )
		\end{equation}
		to the following limit in law
		\begin{equation}\label{riemann_sums}
			S^M_n\overset{(d)}{\to} \varphi(X_{1})- \varphi(X_{0})- \frac{\kappa}{2}\int_0^1 \varphi''(X_r) \dif W_r
		\end{equation}
		for some universal constant \(\kappa\approx 1290\). Similar It\^o formulas in law were deduced for a stochastic process with similar properties \cite{Swa2010}.

		From these probabilistic results we deduce that, starting from a condition $\langle\pi(\Forest{[]} \Forest{[]}), \bX^{S,S}_{s,t} \rangle= \kappa (W_t-W_s)$, we can match the liming object and obtain
		\[
			S^M_n\overset{(d)}{\to} \int_0^1 \varphi(X_r) \dif\bX^{S,S}_r
		\]
		Starting from this identity in law, it would be interesting to understand how Riemann sums like \eqref{eq:rieSu} could be used to define a convergence in law for full rough paths. More generally, it would be also interesting to understand how the remaining  branched models $\bX^{I,I} $, $\bX^{S,I}$ and $\bX^{I,S}$ might arise as limit in law  of proper discretisations.
	\end{remark}

\bibliographystyle{arxivalpha}
\bibliography{refs}

\end{document}